\documentstyle[12pt]{article}
\def\sq{\hbox {\rlap{$\sqcap$}$\sqcup$}}
\def\sq{\hbox {\rlap{$\sqcap$}$\sqcup$}}
\def\R{ {\rm R \kern -.31cm I \kern .15cm}}
\def\C{ {\rm C \kern -.15cm \vrule width.5pt \kern .12cm}}
\def\Z{ {\rm Z \kern -.27cm \angle \kern .02cm}}
\def\N{ {\rm N \kern -.26cm \vrule width.4pt \kern .10cm}}
\def\1{{\rm 1\mskip-4.5mu l} }
\def\lsim{\raise0.3ex\hbox{$<$\kern-0.75em\raise-1.1ex\hbox{$\sim$}}}
\def\gsim{\raise0.3ex\hbox{$>$\kern-0.75em\raise-1.1ex\hbox{$\sim$}}}
\def\noi{\noindent}

\def\beq{\begin{equation}}   \def\eeq{\end{equation}}
\def\bea{\begin{eqnarray}}  \def\eea{\end{eqnarray}}
\def\nn{\nonumber}
\def\noi{\noindent}
\def\beeq{\begin{eqnarray}} \def\eeeq{\end{eqnarray}}
\newcommand\mysection{\setcounter{equation}{0}\section}
\renewcommand{\theequation}{\thesection.\arabic{equation}}
\newcounter{hran} \renewcommand{\thehran}{\thesection.\arabic{hran}}

\def\bmini{\setcounter{hran}{\value{equation}}
   \refstepcounter{hran}\setcounter{equation}{0}
   \renewcommand{\theequation}{\thehran\alph{equation}}\begin{eqnarray}}

\def\bminiG#1{\setcounter{hran}{\value{equation}}
\refstepcounter{hran}\setcounter{equation}{-1}
\renewcommand{\theequation}{\thehran\alph{equation}}
\refstepcounter{equation}\label{#1}\begin{eqnarray}}

\def\emini{\end{eqnarray}\relax\setcounter{equation}{\value{hran}}\ren 
ewcommand{\theequation}{\thesection.\arabic{equation}}}

\begin{document}
\centerline{\Large\bf Long Range Scattering and Modified }
  \vskip 3 truemm \centerline{\Large\bf  Wave Operators for the 
Wave-Schr\"odinger
system\footnote{Work supported in part by NATO Collaborative Linkage 
Grant 976047}} \vskip 0.5
truecm

\centerline{\bf J. Ginibre}
\centerline{Laboratoire de Physique Th\'eorique\footnote{Unit\'e Mixte de
Recherche (CNRS) UMR 8627}}  \centerline{Universit\'e de Paris XI, B\^atiment
210, F-91405 ORSAY Cedex, France}
\vskip 3 truemm
\centerline{\bf G. Velo}
\centerline{Dipartimento di Fisica, Universit\`a di Bologna} 
\centerline{and INFN, Sezione di
Bologna, Italy}

\vskip 1 truecm

\begin{abstract}
We study the theory of scattering for the system consisting of a 
Schr\"odinger equation and a
wave equation with a Yukawa type coupling in space dimension 3. We 
prove in particular the
existence of modified wave operators for that system with no size 
restriction on the data and
we determine the asymptotic behaviour in time of solutions in the 
range of the wave operators.
The method consists in solving the wave equation, substituting the 
result into the
Schr\"odinger equation, which then becomes both nonlinear and 
nonlocal in time, and treating
the latter by the method previously used for a family of generalized 
Hartree equations with
long range interactions.  \end{abstract}

\vskip 3 truecm
\noi AMS Classification : Primary 35P25. Secondary 35B40, 35Q40, 81U99.  \par
\noi Key words : Long range scattering, modified wave operators, 
Wave-Schr\"odinger system.
\vskip 1 truecm

\noindent LPT Orsay 01-48 \par
\noindent June 2001 \par

\newpage
\pagestyle{plain}
\baselineskip 18pt

\mysection{Introduction}
\hspace*{\parindent} This paper is devoted to the theory of 
scattering and more precisely to
the existence of modified wave operators for the Wave-Schr\"odinger (WS) system
$$\hskip 3 truecm \left \{ \begin{array}{ll} i\partial_t u = - {1 
\over 2} \Delta u - Au &\hskip
6.8 truecm(1.1) \\ & \\ \sq A = |u|^2 &\hskip 6.8 truecm(1.2) 
\end{array} \right .$$
\noi where $u$ and $A$ are respectively a complex valued and a real 
valued function defined in
space time ${I\hskip-1truemm R}^{3+1}$, $\Delta$ is the Laplacian in 
${I\hskip-1truemm R}^3$
and $\sq = \partial_t^2 - \Delta$ is the d'Alembertian in 
${I\hskip-1truemm R}^{3+1}$. That
system is Lagrangian with Lagrangian density
$${\cal L} = i \left ( \bar{u} \ \partial_t \ u - u \ \partial_t\ 
\bar{u}\right ) - {1 \over
2} |\nabla u|^2 + {1 \over 2} (\partial_t A)^2 -
{1 \over 2} |\nabla A|^2 + A|u|^2 \ .
\eqno(1.3)$$ \noi Formally, the $L^2$ norm of $u$ is conserved, as 
well as the energy
$$E(u,A) = \int dx \Big \{ {1 \over 2} \left ( |\nabla u|^2 + 
(\partial_t A)^2 + |\nabla A|^2
\right ) - A|u|^2 \Big \} \ . \eqno(1.4)$$
\noi The Cauchy problem for the WS system (1.1) (1.2) is known to be 
globally well posed in
the energy space $X_e = H^1 \oplus \dot{H}^1 \oplus L^2$ for $(u, A, 
\partial_t A)$ \cite{1r}
\cite{2r} \cite{4r} \cite{15r}. \par

A large amount of work has been devoted to the theory of scattering 
for nonlinear equations
and systems centering on the Schr\"odinger equation, in particular 
for nonlinear Schr\"odinger
(NLS) equations, Hartree equations, Klein-Gordon Schr\"odinger (KGS) 
and Maxwell-Schr\"odinger
(MS) systems. As in the case of the linear Schr\"odinger equation, 
one must distinguish the
short range case from the long range case. In the former case, 
ordinary wave operators are
expected and in a number of cases proved to exist, describing 
solutions where the Schr\"odinger
function behaves asymptotically like a solution of the free 
Schr\"odinger equation. In the
latter case, ordinary wave operators do not exist and have to be 
replaced by modified wave
operators including a suitable phase in their definition. In that 
respect, the WS system (1.1)
(1.2) in ${I\hskip-1truemm R}^{3+1}$ belongs to the borderline 
(Coulomb) long range case,
because of the $t^{-1}$ decay in $L^{\infty}$ norm of solutions of 
the wave equation. Such is
the case also for the Hartree equation with $|x|^{-1}$ potential. 
Both are simplified models
for the more complicated Maxwell-Schr\"odinger system in 
${I\hskip-1truemm R}^{3+1}$, which
belongs to the same case, as well as the KGS system in 
${I\hskip-1truemm R}^{2+1}$. \par

Whereas a well developed theory of long range scattering exists for 
the linear Schr\"odinger
equation (see \cite{3r} for a recent treatment and for an extensive 
bibliography), there exist
only few results on nonlinear long range scattering. The existence of 
modified wave operators
in the borderline Coulomb case has been proved for the NLS equation 
in space dimension $n =
1$ \cite{19r}. That result has been extended to the NLS equation in 
dimensions $n = 2,3$ and
to the Hartree equation in dimension $n \geq 2$ \cite{5r}, to the 
derivative NLS equation in
dimension $n = 1$ \cite{14r}, to the KGS system in dimension 2 
\cite{20r} and to the MS system
in dimension 3 \cite{22r}. All those results are restricted to the 
case of small data. \par

In a recent series of papers, \cite{6r} \cite{7r} \cite{8r}, we 
proved the existence of
modified wave operators for a family of Hartree type equations with 
general (not only Coulomb)
long range interactions and without any size restriction on the data. 
The method is strongly
inspired by a previous series of papers by Hayashi et al \cite{9r} 
\cite{10r} \cite{11r}
\cite{12r} \cite{13r} on the Hartree equation. In the latter papers 
it is proved first in the
borderline Coulomb case and then in the whole long range case, that 
the global solutions of the
Hartree equation with small initial data exhibit an asymptotic 
behaviour for large time that is
typical of long range scattering and includes in particular the 
expected relevant phase factor.
\par

The present paper is devoted to the extension of the results of 
\cite{6r} \cite{7r} \cite{8r}
to the WS system and in particular to the proof of the existence of 
modified wave operators for
that system without any size restriction on the data. The method 
consists in eliminating the
wave equation by solving it for $A$ in terms of $u$ and substituting 
the result into the
Schr\"odinger equation, thereby obtaining a new Schr\"odinger 
equation which is both nonlinear
and nonlocal in time. The latter is then treated as the Hartree 
equation in \cite{6r}
\cite{7r} \cite{8r}, namely $u$ is expressed in terms of an amplitude 
$w$ and a phase $\varphi$
satisfying an auxiliary system similar to that introduced in 
\cite{11r}. Wave operators are
constructed first for that auxiliary system, and then used to 
construct modified wave
operators for the original system (1.1). The detailed construction is 
too complicated to allow
for a more precise description at this stage, and will be described 
in heuristic terms in
Section 2 below. In subsequent papers, the results of the present one 
will be extended to the
case of the MS system. \par

We now give a brief outline of the contents of this paper. A more 
detailed description of the
technical parts will be given at the end of Section 2. After 
collecting some notation and
preliminary estimates in Section 3, we study the asymptotic dynamics 
for the auxiliary system
in Section 4 and uncover some difficulties due to the different 
propagation properties of
solutions of the wave and Schr\"odinger equations. As a preparation 
for the general case, we
construct in Section 5 the wave operators associated with the 
simplified linear system
obtained by replacing (1.2) by the free wave equation $\sq A = 0$. We 
then solve the local
Cauchy problem at infinity for the auxiliary system in Sections 6 and 
7, which contain the
main technical results of this paper.  We finally come back from the 
auxiliary system to the
original one (1.1) (1.2) and construct the modified wave operators 
for the latter in Section
8, where the final result is stated in Proposition 8.1. \par

We conclude this section with some general notation which will be 
used freely throughout this
paper. We denote by $\parallel \cdot \parallel_r$ the norm in $L^r 
\equiv L^r({I\hskip-1truemm
R}^3)$ and we define $\delta (r) = 3/2 - 3/r$. For any interval $I$ 
and any Banach space $X$,
we denote by ${\cal C}(I,X)$ (resp. ${\cal C}_w(I,X))$ the space of 
strongly (resp. weakly)
continuous functions from $I$ to $X$ and by $L^{\infty}(I,X)$ (resp. 
$L^{\infty}_{loc}(I,X))$
the space of measurable essentially bounded (resp. locally 
essentially bounded) functions from
$I$ to $X$. For real numbers $a$ and $b$, we use the notation $a \vee 
b = {\rm Max}(a,b)$, and
$a \wedge b = {\rm Min} (a, b)$. Furthermore, we define

$$\begin{array}{llll}
&[a \vee b] &= a \vee b &\hbox{if \ $a \not= b$} \\
& & & \\
&&= a + \varepsilon &\hbox{for some \ $\varepsilon > 0$ \ if \ $a = b$} \ , \\
& & &\\
&[a \wedge b] &= a + b - [a \vee b] &\hbox{and} \ \ [a]_+ = [a \vee 0] \ .
\end{array}$$

\noi For any interval $I \subset {I\hskip-1truemm R}^+$, we denote by 
$\bar{I}$ the closure of
$I$ in ${I\hskip-1truemm R}^+ \cup \{\infty\}$ and for any interval 
$I = [a,b)$ we denote by
$I_+$ the interval $I_+ = [a, \infty )$.  In the estimates of 
solutions of the relevant
equations, we shall use the letter $C$ to denote constants, possibly 
different from an estimate
to the next, depending on various parameters, but not on the 
solutions themselves or on their
initial data. We shall use the notation $C(a_1,a_2, \cdots )$ for 
estimating functions, also
possibly different from an estimate to the next, depending in 
addition on suitable norms $a_1,
a_2, \cdots $ of the solutions or of their initial data. Additional 
notation will be given in
Section 3.

\mysection{Heuristics} \hspace*{\parindent} In this section, we 
discuss in heuristic terms the
construction of the modifed wave operators for the system (1.1) 
(1.2), as it will be performed
in this paper. We refer to Section 2 of \cite{6r} \cite{7r} for 
general background and for a
similar discussion adapted to the case of the Hartree equation. The 
problem that we want to
address is that of classifying the possible asymptotic behaviours in 
time of the solutions of
(1.1) (1.2) by relating them to a set of model functions ${\cal V} = 
\{ v= v(v_+)\}$ parametrized
by some data $v_+$ and with suitably chosen and preferably simple 
asymptotic behaviour in
time. For each $v \in {\cal V}$, one tries to construct a solution 
$(u, A)$ of (1.1) (1.2)
such that $(u, A)(t)$ behaves as $v(t)$ when $t \to \infty$ in a 
suitable sense. We then
define the wave operator as the map $\Omega : v_+ \to (u, A)$ thereby 
obtained. A similar
question can be asked for $t \to - \infty$. We restrict our attention to positive time. The
more standard definition of the wave operator is to defin 
e it as the map $v_+ \to (u,A)(0)$,
but what really matters is the solution $(u, A)$ in the neighborhood 
of infinity in time,
namely in some interval $[T, \infty )$, and continuing such a 
solution down to $t = 0$ is
a somewhat different question which we shall not touch here. \par

In cases such as (1.1) (1.2) where the system of interest is a 
perturbation of a simple
linear system, hereafter called the free system, a natural candidate 
for ${\cal V}$ is the
set of solutions of the free system, parametrized by the initial data 
$v_+$ at time $t =
0$ for the Cauchy problem for that system. In the case of the system 
(1.1) (1.2) one is
therefore tempted to consider the Cauchy problem
\beq
\label{2.1e}
\left \{ \begin{array}{ll} i\partial_t u = - \displaystyle{{1 \over 
2}} \Delta u &\qquad u(0) =
u_+ \\ &\\
\sq A = 0 &A(0) = A_+\ , \ \partial_t A(0) = \dot{A}_+ \ ,\end{array}
\right . \eeq
\noi to take $v_+ = (u_+, A_+, \dot{A}_+)$ and to take for $v(v_+)$ 
the solution $(u, A)$ of
(2.1). Cases where such a procedure yields an adequate set ${\cal V}$ 
are called short range
cases. They require that the perturbation has sufficient decay in 
time or equivalently in
space. This is the case for instance for the linear Schr\"odinger 
equation or for the Hartree
equation with potential $V(x) = |x|^{-\gamma}$ for $\gamma > 1$. Such 
is not the case however
for the system (1.1) (1.2). This shows up through the fact that the 
solution $A$ of the wave
equation $\sq A = 0$ decays at best as $t^{-1}$ (in $L^{\infty}$ 
norm), which is the
borderline case of nonintegrability in time. That situation 
corresponds to the limiting case
$\gamma = 1$ (the Coulomb case in space dimension $n = 3$) for the 
linear Schr\"odinger and
for the Hartree equation. A similar situation prevails for the KGS 
system in space dimension 2
and for the MS system in space dimension 3. In the present case, 
which is the borderline long
range case, the set of solutions of the Cauchy problem (2.1) is 
inadequate, and one of the
tasks that will be performed in this paper (see especially Sections 7 
and 8) will be to
construct a better set ${\cal V}$ of model asymptotic functions. \par

Constructing the wave operators essentially amounts to solving the 
Cauchy problem with
infinite initial time. The system (1.1) (1.2) in this form is not 
well suited for that purpose
and we shall now perform a number of transformations leading to an 
auxiliary system for which
that problem can be handled. For additional flexibility we shall 
first of all allow for
imposing initial data at two different initial times $t_0$ and $t_1$ 
for the Schr\"odinger
and wave equations respectively. With the aim of letting $t_1$ and 
$t_0$ tend to infinity in
that order, we shall take $t_0 \leq t_1$. We shall then eliminate the 
wave equation by
solving it and substituting the result into the Schr\"odinger 
equation. We define
$$\omega = (-\Delta )^{1/2} \quad , \quad K(t) = \omega^{-1} \sin 
\omega t \quad , \quad
\dot{K}(t) = \cos \omega t$$
\noi and we replace the wave equation (1.2) by its solution
\beq
\label{2.2e}
A = A_0 + A_1^{t_1} (|u|^2)
\eeq
\noi where
\bea
\label{2.3e}
&&A_0 = \dot{K}(t) \ A_+ + K(t) \ \dot{A}_+ \\
&&A_1^{t_1} (|u|^2) = \int_{t_1}^t dt'\ K(t - t')\ |u(t')|^2 \ .
\label{2.4e}
\eea
\noi Here $A_0$ is a solution of the free wave equation, with initial 
data $(A_+, \dot{A}_+)$ at
time $t = 0$. For $t_1 = \infty$, $(A_+, \dot{A}_+)$ is naturally 
interpreted as the asymptotic
state for $A$, in keeping with the previous discussion. \par

The Cauchy problem for the system (1.1) (2.2) with initial data 
$u(t_0) = u_0$ is no longer a
usual PDE Cauchy problem because $A_1$ depends on $u$ nonlocally in 
time. A convenient way to
handle that difficulty is to first replace that problem by a partly 
linearized form thereof,
namely
\beq
\label{2.5e}
\left \{ \begin{array}{ll} i\partial_t u' = - \displaystyle{{1 \over 
2}} \Delta u' - Au' \ ,
&\quad u'(t_0) = u_0 \\ &\\
  A = A_0 + A_1 (|u|^2) \ .&\end{array}
\right . \eeq
\noi For given $u$, (\ref{2.5e}) is an ordinary (linear) Cauchy 
problem for $u'$. Solving that
problem for $u'$ defines a map $\Gamma : u \to u'$, and solving the 
original problem then
reduces to finding a fixed point of $\Gamma$, which in favourable 
cases can be done for
instance by contraction. We shall make use of that linearization 
method, not for the equation
for $u$, but for the auxiliary system to be defined below. \par

Aside from the nonlocality in time of the nonlinear interaction term, 
which can be handled by
the previous linearization, the system (1.1) (2.2) is rather similar 
to the Hartree type
equations considered in \cite{6r} \cite{7r} \cite{8r}, and we next 
perform the same change of
variables, which is well adapted to the study of the asymptotic 
behaviour in time. The unitary
group
\beq
\label{2.6e}
U(t) = \exp (i(t/2)\Delta )
\eeq
\noi which solves the free Schr\"odinger equation can be written as
\beq
\label{2.7e}
U(t) = M(t) \ D(t) \ F \ M(t)
\eeq
\noi where $M(t)$ is the operator of multiplication by the function
\beq
\label{2.8e}
M(t) = \exp  \left ( i x^2/2t \right ) \ ,
\eeq
\noi $F$ is the Fourier transform and $D(t)$ is the dilation operator
\beq
\label{2.9e}
(D(t)\ f)(x) = (it)^{-n/2} \ f(x/t)
\eeq
\noi normalized to be unitary in $L^2$. We shall also need the 
operator $D_0(t)$ defined by
\beq
\label{2.10e}
\left ( D_0(t) f\right )(x) = f(x/t) \ .
\eeq
\noi We now parametrize $u$ in terms of an amplitude $w$ and of a 
real phase $\varphi$ as
\beq
\label{2.11e}
u(t) = M(t) \ D(t) \exp [- i \varphi (t)] w(t) \ .
\eeq
\noi Substituting (\ref{2.11e}) into (1.1) yields an evolution 
equation for $(w , \varphi )$,
namely  \beq
\label{2.12e}
\left \{ i \partial_t + (2t^2)^{-1} \Delta - i(2t^2)^{-1} (2 \nabla 
\varphi \cdot \nabla +
\Delta \varphi ) + t^{-1} B + \partial_t \varphi - (2t^2)^{-1} 
|\nabla \varphi |^2 \right \} w
= 0 \eeq
\noi where we have expressed $A$ in terms of a new function $B$ by
\beq
\label{2.13e}
A = t^{-1} \ D_0\ B \ .
\eeq
\noi Corresponding to the decomposition (\ref{2.2e}) of $A$, we decompose
\beq
\label{2.14e}
B = B_0 + B_1^{t_1} (w , w)
\eeq
\noi where
\beq
\label{2.15e}
B_1^{t_1}(w_1,w_2) = \int_1^{t_1/t} d \nu \ \nu^{-3} \ \omega^{-1} 
\sin ((\nu - 1) \omega)
D_0(\nu) ({\rm Re} \ \bar{w}_1w_2)(\nu t) \ .\eeq
\indent As in the case of the Hartree equation, we have only one 
evolution equation for two
functions $(w, \varphi )$ and we arbitrarily split that equation into 
two equations for $w$
and $\varphi$. For that purpose, we split $B_1^{t_1}$ into a short 
range and a long range
parts
\beq
\label{2.16e}
B_1^{t_1} = B_S^{t_1} + B_L^{t_1} \ .
\eeq
\noi We take $0 < \beta < 1$ and we define
\beq
\label{2.17e}
\left \{ \begin{array}{l} \left ( FB_S^{t_1} \right )(t, \xi ) = \chi 
(|\xi | > t^{\beta})
FB_1^{t_1}(t, \xi ) \\ \\
\left ( FB_L^{t_1} \right )(t, \xi ) = \chi (|\xi | \leq t^{\beta})
FB_1^{t_1}(t, \xi )   \end{array}
\right . \eeq
\noi where $\chi (|\xi| \displaystyle{\mathrel{\mathop <_{>}}} 
t^{\beta})$ is the characteristic
function of the set $\{(t, \xi ) : |\xi | 
\displaystyle{\mathrel{\mathop <_{>}}} t^{\beta}\}$.
The parameter $\beta$ will satisfy various conditions which will 
appear later, all of which
will be compatible with $\beta = 1/2$. We now split the equation 
(\ref{2.12e}) into the
following system of two equations
\beq
\label{2.18e}
\left \{ \begin{array}{l} \partial_t w = i(2t^2)^{-1} \Delta w + 
t^{-2} Q (\nabla \varphi , w)
+ it^{-1} (B_0 + B_S^{t_1}(w, w)) w \\ \\
  \partial_t \varphi = (2t^2)^{-1} |\nabla \varphi |^2 - t^{-1} \ 
B_L^{t_1}(w,w) \end{array}
\right . \eeq
\noi where we have defined
\beq  \label{2.19e}
Q(s, w) = s \cdot \nabla w + (1/2) (\nabla \cdot s) w
\eeq
\noi for any vector field $s$. Since the right-hand sides of 
(\ref{2.18e}) contain $\varphi$
only through its gradient, we can obtain from (\ref{2.18e}) a closed 
system for $w$ and $s =
\nabla \varphi$ by taking the gradient of the second equation, namely
\beq
\label{2.20e}
\left \{ \begin{array}{l} \partial_t w = i(2t^2)^{-1} \Delta w + 
t^{-2} Q (s , w)
+ it^{-1} (B_0 + B_S^{t_1}(w, w)) w \\ \\
  \partial_t s = t^{-2} s \cdot \nabla s - t^{-1} \nabla 
B_L^{t_1}(w,w) \ .\end{array}
\right .  \eeq
\indent Once the system (\ref{2.20e}) is solved for $(w, s)$, one 
recovers $\varphi$ easily by
integrating the second equation of (\ref{2.18e}) over time. We refer 
to \cite{6r} for details.
The system (\ref{2.20e}) will be referred to as the auxiliary system 
and will play an essential
role in this paper. For the same reason as was explained for the 
partly resolved system (1.1)
(\ref{2.2e}), we  shall use at intermediate stages a partly 
linearized version of the system
(\ref{2.20e}), namely
\beq
\label{2.21e}
\left \{ \begin{array}{l} \partial_t w' = i(2t^2)^{-1} \Delta w' + 
t^{-2} Q (s , w')
+ it^{-1} (B_0 + B_S^{t_1}(w, w)) w' \\ \\
  \partial_t s' = t^{-2} s \cdot \nabla s' - t^{-1} \nabla 
B_L^{t_1}(w,w) \end{array}
\right . \eeq
\noi to be considered as a system of equations for $(w', s')$ for 
given $(w, s)$. The first
question to be considered is whether the auxiliary system 
(\ref{2.20e}) defines a dynamics for
large time, namely whether the Cauchy problem for that system is 
locally well posed in a
neighborhood of infinity in time, more precisely has a unique 
solution defined up to infinity
in time for sufficiently large $t_1$ and sufficiently large initial 
time $t_0$, possibly
depending on the size of the initial data. This property was 
satisfied by the corresponding
auxiliary system associated with the Hartree equation and considered 
in \cite{6r} \cite{7r}.
Here however we encounter serious difficulties associated with the 
difference of propagation
properties of solutions of the Schr\"odinger and wave equations. In 
fact a typical solution of
the free Schr\"odinger equation behaves asymptotically in time as
$$(U(t) u_+) (x) \sim (MDFu_+) (x) = \exp (ix^2/2t)(it)^{-3/2} \ Fu_+(x/t)$$
\noi namely spreads by dilation by $t$ in all directions in the 
support of $Fu_+$, while by the
Huyghens principle $A_0$ remains concentrated in a neighborhood of 
the light cone, more
precisely within a distance $R$ of the latter if the initial data 
$(A_+, \dot{A}_+)$ are
supported in a ball of radius $R$. When switching to the new 
variables $(w, B)$, $w$ tends to
a limit when $t \to \infty$ whereas $B_0$ concentrates in a 
neighborhood of the unit sphere,
within a distance $R/t$ of the latter in the previous case of 
compactly supported data. Note
however that for $t_1 = \infty$, $B_1^{\infty}$ is expected to tend 
to a limit like $w$ and
not to concentrate like $B_0$, as can be guessed from (\ref{2.15e}). \par

We shall treat the auxiliary system (\ref{2.20e}) by energy methods, 
and in particular look
for $w$ in spaces of the type ${\cal C}([T, \infty ), H^k)$ where 
$H^k$ is the usual Sobolev
space based on $L^2$. In order to treat the nonlinear term 
$B_1(w,w)$, we shall need a minimal
regularity, in practice $k > 1$. However, when taking $H^k$ norms of 
$B_0$, the previous
concentration phenomenon implies
$$\parallel B_0 ; H^k\parallel \ \sim O \left ( t^{k-1/2} \right )$$
\noi which has worse and worse asympotic behaviour in time as $k$ 
increases. This difficulty
manifests itself in the following way~: \par

(i) If $t_0 = t_1 < \infty$, the available estimates for the sytem 
(\ref{2.20e}) do not
prevent finite time blow up after $t_0$, even if $A_0 = 0$. \par

This encourages us to take $t_1 > t_0$, and actually the situation 
becomes slightly better in
that case. Nevertheless \par

(ii) the available estimates do not prevent finite time blow up after 
$t_1$, which is the
same fact as (i) with $t_0$ replaced by $t_1$, and \par

(iii) if $A_0 \not= 0$ and if $t_1$ is sufficiently large, the 
available estimates do not
prevent blow up before $t_1$. \par

A definite improvement occurs however if $A_0 = 0$. \par

(iv) If $A_0 = 0$, the available estimates allow for a proof of 
existence of solutions in
$[t_0, t_1]$ for $t_0$ sufficiently large and arbitrary $t_1 > t_0$, 
possibly $t_1 = \infty$.
In particular for $t_1 = \infty$, the solutions are defined up to 
infinity in time.
Furthermore, for those solutions, $w(t)$ has a limit $w_+$ as $t \to 
\infty$. \par

The last case brings us in the same situation as that encountered for 
the Hartree equation in
\cite{6r} \cite{7r} and could be taken as the starting point for the 
construction of partial
modified wave operators (restricted to the case of vanishing $(A_+, 
\dot{A}_+))$ by the same
method as in \cite{6r} \cite{7r}. We shall however refrain from 
performing that construction
and turn directly to the case of nonvanishing $(A_+, \dot{A}_+)$. In 
that case, the need to use
$H^k$ norms with $k > 1$ for $w$ makes the treatment of $A_0$ 
nontrivial, even if one drops the
interaction term $A_1$. As a preparation for the general case, we 
shall therefore first
construct the wave operators at the same level of regularity for the 
simplified system
\beq
\label{2.22e}
\left \{ \begin{array}{l} i \partial_t u = - (1/2) \Delta u - A_0 u \\
\\
  \sq A_0 = 0 \end{array}
\right . \eeq
\noi namely for a linear Schr\"odinger equation with time dependent 
potential $A_0$
satisfying the free wave equation. After the appropriate change of variables
\beq  \label{2.23e}
u = MD w \qquad , \qquad A_0 = t^{-1} \ D_0 \ B_0
\eeq
\noi the Schr\"odinger equation becomes
\beq
\label{2.24e}
R(w) \equiv \partial_t w - i(2t^2)^{-1} \Delta w - i t^{-1} B_0 w = 0 \ .
\eeq
\indent The construction of the wave operators for that equation in 
$L^2$, either in the form
(\ref{2.22e}) or (\ref{2.24e}) can be easily performed by a simple 
variant of Cook's method,
and the construction of the wave operators at the level of $H^k$ 
becomes a regularity problem
for the previous wave operators. Solving that problem for $k \geq 1$ 
(in fact for $k > 1/2$)
requires special assumptions on the asymptotic states $(w_+, A_+, 
\dot{A}_+)$, to the effect
that the product $B_0w_+$ decays faster in time in the relevant norms 
than what would
naturally follow from factorized estimates. Those assumptions can be 
ensured for instance by
imposing support properties of $w_+$, to the effect that $w_+ = 0$ on 
the unit sphere, and
suitable decay of $(A_+,\dot{A}_+)$ at infinity in space. They will 
be needed again in the
treatment of the general problem. \par

The construction of the modified wave operators in the general case 
follows the same pattern as
for the Hartree equation. The aim is to construct solutions of the 
auxiliary system
(\ref{2.20e}) with suitably prescribed asymptotic behaviour at 
infinity, and in particular with
$w(t)$ tending to a limit $w_+$ as $t \to \infty$. For that purpose, 
we choose a suitable
asymptotic pair $(W, \phi )$, and therefore a pair $(W, S)$ with $S = 
\nabla \phi$ with $W(t)$
tending to $w_+$ as $t \to \infty$, and we look for a solution $(w, 
\varphi )$ of the system
(\ref{2.18e}) in the form
\beq
\label{2.25e}
(w, \varphi ) = (W, \phi ) + (q, \psi)
\eeq
\noi or equivalently for a solution $(w, s)$ of the system 
(\ref{2.20e}) in the form
\beq
\label{2.26e}
(w, s ) = (W, S) + (q, \sigma)
\eeq
\noi where $(q, \psi , \sigma )$ tends to zero as $t \to \infty$. 
Actually for technical
reasons, we shall need to modify the auxiliary system slightly, in 
the following way. When
expanding $w = W + q$ in $B_1^{t_1}(w, w)$, we shall replace that quantity by
  \beq
\label{2.27e}
B_1^{t_1, \infty }(w, w) \equiv B_1^{\infty}(W,W) + 2 B_1^{t_1}(W,q) 
+ B_1^{t_1}(q,q) \ .
\eeq
\noi We furthermore define the remainders
\bea
\label{2.28e}
&&R_1(W,S) = \partial_t W - i(2t^2)^{-1} \Delta W - t^{-2} Q(S,W) - 
it^{-1} (B_0 + B_S^{\infty}
(W, W)) W \nn \\
&& \\
&&R_2(W,S) = \partial_tS - t^{-2}S \cdot \nabla S + t^{-1} \nabla 
B_L^{\infty}(W,W) \ .
\label{2.29e}
\eea
\noi The modified auxiliary system for $(q, \sigma )$ then takes the form
\bea
\label{2.30e}
\left \{ \begin{array}{l} \partial_t q = i (2t^2)^{-1} \Delta q + 
t^{-2} (Q(s, q) + Q(\sigma,
W)) + it^{-1} B_0 q  \\
\\
+ i t^{-1} B_S^{t_1, \infty} (w, w) q + it^{-1} \left ( 2 
B_S^{t_1}(W,q) + B_S^{t_1}(q,q)
\right ) W - R_1(W,S)\\
\\
\partial_t \sigma = t^{-2} ( s \cdot \nabla \sigma + \sigma \cdot 
\nabla S) - t^{-1}\nabla
\left ( 2 B_L^{t_1}(W,q) + B_L^{t_1}(q,q)\right ) - R_2 (W, S) \ .
  \end{array}
\right .
\eea
\noi Note that changing $B^{t_1}(w,w)$ to $B^{t_1, \infty}(w,w)$ 
changes $A$ by a solution of
the free wave equation, so that we are still solving the original 
system (1.1) (1.2), with
however a slightly different $A_0$ as compared with (\ref{2.2e}). \par

For the same reason as for the partly resolved system (1.1) (2.2) and 
for the auxiliary system
(\ref{2.20e}), we shall use at intermediate stages a partly 
linearized version of the system
(\ref{2.30e}), namely
\bea
\label{2.31e}
\left \{ \begin{array}{l} \partial_t q' = i (2t^2)^{-1} \Delta q' + 
t^{-2} (Q(s, q') + Q(\sigma,
W)) + it^{-1} B_0 q'  \\
\\
+ i t^{-1} B_S^{t_1, \infty} (w, w) q' + it^{-1} \left ( 2 
B_S^{t_1}(W,q) + B_S^{t_1}(q,q)
\right ) W - R_1(W,S)\\
\\
\partial_t \sigma ' = t^{-2} ( s \cdot \nabla \sigma ' + \sigma \cdot 
\nabla S) - t^{-1}\nabla
\left ( 2 B_L^{t_1}(W,q) + B_L^{t_1}(q,q)\right ) - R_2 (W, S) \ .
  \end{array}
\right .
\eea

The construction of solutions $(q, \sigma )$ tending to zero at 
infinity for the system
(\ref{2.30e}) with $t_1 = \infty$ proceeds in several steps. We 
assume first that $(W, S)$
satisfy suitable boundedness properties and that the remainders 
$R_1(W, S)$ and $R_2(W, S)$
satisfy suitable decay in time. We solve the linearized system 
(\ref{2.31e}) for $(q', \sigma
')$ for given $(q, \sigma )$, both with finite and infinite time 
$t_1$ and initial time $t_0$.
We then solve (\ref{2.30e}) by proving that the map $\Gamma : (q , 
\sigma ) \to (q', \sigma
')$ is a contraction in suitable norms. We also prove that the 
solution of (\ref{2.30e}) with
$t_0 = t_1 < \infty$ converges to the solution with $t_0 = t_1 = 
\infty$ when $t_0 \to
\infty$, a property which is natural in the framework of scattering 
theory. We finally
construct $(W, S)$ satisfying the required boundedness and decay 
properties by solving the
auxiliary system (\ref{2.20e}) with $t_1 = \infty$ approximately by 
iteration, restricting
our attention to the second iteration, which is sufficient to cover 
the range $ 1 < k < 2$.
The pair $(W, S)$ or equivalently $(W, \phi )$ depends only on the 
limit $w_+$ of $w$ and
$W$ as $t \to \infty$. Solving the auxiliary system (\ref{2.30e}) 
with that $(W, S)$ and with
$t_0 = t_1 = \infty$ yields a solution $(w, s)$ of the system 
(\ref{2.20e}) and therefore a
solution $(w, \varphi )$ of the system (\ref{2.18e}) with prescribed 
asymptotic behaviour
characterized by $(W, S)$ or $(W, \phi )$. That solution depends on 
$(w_+, A_+, \dot{A}_+)$.
Plugging that solution into (\ref{2.11e}) (\ref{2.13e}) and taking 
$w_+ = Fu_+$ yields a map
$\Omega : (u_+, A_+, \dot{A}_+) \to (u, A)$ which is the required 
wave operator for the
system (1.1) (1.2). \par

The auxiliary system (\ref{2.18e}) satisfies a gauge invariance 
property similar to that of the
corresponding system for the Hartree equation used in \cite{6r} 
\cite{7r}, and the construction
of the intermediate wave operator for that system can be made in a 
gauge covariant way. For
brevity we shall refrain from discussing that question in this paper. \par

We now describe the contents of the technical parts of this paper, 
namely Sections 3-8. In
Section 3, we introduce some notation, define the relevant function 
spaces and collect a number
of preliminary estimates. In Section 4, we study the Cauchy problem 
for large time for the
auxiliary system (\ref{2.20e}). We solve the Cauchy problem with 
finite initial time for the
linearized system (\ref{2.21e}) (Proposition 4.1), we prove a number 
of uniqueness results
for the system (\ref{2.20e}) (Proposition 4.2), we prove the 
existence of a limit $w(t)$ of
$w_+$ for suitably bounded solutions of the system (\ref{2.20e}) 
(Proposition 4.3), we
discuss in more quantitative terms the possible occurrence of blow up 
mentioned above, and we
finally solve the Cauchy problem for the system (\ref{2.20e}) with 
$t_1 = \infty$ and large
$t_0$ in the special case $A_0 = 0$ (Proposition 4.4). \par

In Section 5, as a preparation for the construction of the wave 
operators for the system
(\ref{2.20e}) with $A_0 \not= 0$, we study the existence of wave 
operators for the linear
problem (\ref{2.22e}) in the form (\ref{2.24e}). In particular we 
prove the existence of
$L^2$-wave operators by a variant of Cook's method (Proposition 5.2), 
we prove the $H^k$
regularity of those wave operators under suitable decay assumptions 
of $R(W)$ for the model
function $W$ (Proposition 5.3) and we finally reduce those decay 
properties to conditions on
the asymptotic state $(w_+, A_+, \dot{A}_+)$. In Section 6 and 7, we 
study the Cauchy problem
at infinity in the general case $A_0 \not= 0$ for the auxiliary 
system (\ref{2.20e}) in the
difference form (\ref{2.30e}). Under suitable boundedness assumptions 
on $(W,S)$ and decay
assumptions on $R_1(W,S)$ andf $R_2(W,S)$ we prove the existence of 
solutions for $t_0$ and
$t_1$ finite and infinite, first for the linearized system 
(\ref{2.31e}) (Propositions 6.1
and 6.2) and then for the nonlinear system (\ref{2.30e}) (Proposition 
6.3). We then choose
appropriate $(W,S)$, prove that they satisfy the required assumptions 
(Lemmas 7.1 and 7.2)
and finally state the result on the Cauchy problem at infinity for 
the system (\ref{2.30e})
in $H^k$ for $1 < k < 2$ (Proposition 7.1). Finally in Section 8, we 
construct the wave
operators for the system (1.1) (1.2) from the results previously 
obtained for the system
(\ref{2.30e}) and we derive the asymptotic estimates for the 
solutions $(u, A)$ in their
range that follow from the previous estimates (Proposition 8.1).

\mysection{Notation and preliminary estimates}
\hspace*{\parindent} In this section we introduce some additional 
notation and we collect a
number of estimates which will be used throughout this paper. We 
shall use the Sobolev spaces
$H_r^k$ defined for $1 \leq r \leq \infty$ by
$$H_r^k = \left \{ u : \parallel u ; H_r^k\parallel \ \equiv \ 
\parallel <\omega >^k u
\parallel_r \ < \infty \right \}$$
\noi where $< \cdot > = (1 + |\cdot |^2)^{1/2}$. The subscript $r$ 
will be omitted if $r = 2$.
\par

We shall look for solutions of the auxiliary system (\ref{2.20e}) in 
spaces of the type ${\cal
C}(I, X^{k, \ell})$ where $I$ is an interval and
$$X^{k, \ell} = H^k \oplus \omega^{-1} \ H^{\ell}$$
\noi namely
\beq
\label{3.1e}
X^{k, \ell} = \left \{ (w, s) : w \in H^k \ , \ \nabla s \in H^{\ell} \right \}
\eeq
\noi and we shall use the notation
\beq
\label{3.2e}
\parallel w; H^k\parallel\ = |w|_k \ , \ \parallel \omega s ; 
H^{\ell}\parallel \ \equiv \
\parallel \nabla s ; H^{\ell} \parallel \ = |s|_{\ell}^{\dot{}} \ .  \eeq
\noi Here it is understood that $\nabla s \in L^2$ includes the fact 
that $s \in L^6$. \par

We shall use extensively the following Sobolev inequalities, stated 
here in ${I\hskip-1truemm
R}^n$, but to be used only for $n = 3$. \\

\noi {\bf Lemma 3.1.} {\it Let $1 < q$, $r < \infty$, $1 < p \leq 
\infty$ and $0 \leq j < k$.
If $p = \infty$, assume that $k - j > n/r$. Let $\sigma$ satisfy $j/k 
\leq \sigma \leq 1$ and
$$n/p - j = (1 - \sigma )n/q + \sigma (n/r - k) \ .$$
\noi Then the following inequality holds
\beq
\label{3.3e}
\parallel \omega^j u \parallel_p \ \leq C \parallel u \parallel_q^{1 
- \sigma} \ \parallel
\omega^ku \parallel_r^{\sigma} \ .  \eeq} \\
\indent The proof follows from the Hardy-Littlewood-Sobolev (HLS) 
inequality (\cite{16r},
p.~117) (from the Young inequality if $p = \infty$), from 
Paley-Littlewood theory and
interpolation.\\

We shall also use extensively the following Leibnitz and commutator 
estimates.\\

\noi {\bf Lemma 3.2.} {\it Let $1 < r, r_1, r_3 < \infty$ and
$$1/r = 1/r_1 + 1/r_2 = 1/r_3 + 1/r_4 \ .$$
\noi Then the following estimates hold
\beq
\label{3.4e}
\parallel \omega^m (uv) \parallel_r \ \leq C \left ( \parallel 
\omega^m u \parallel_{r_1}
  \ \parallel v \parallel_{r_2} + \parallel  \omega^m v \parallel_{r_3} \
  \parallel u \parallel_{r_4} \right ) \eeq
\noi for $m \geq 0$, and
\beq
\label{3.5e}
\parallel [\omega^m , u]v\parallel_r \ \leq C \left ( \parallel 
\omega^m u \parallel_{r_1} \
  \parallel v \parallel_{r_2} + \parallel  \omega^{m -1} v \parallel_{r_3} \
  \parallel \nabla u \parallel_{r_4} \right ) \eeq
\noi for $m \geq 1$, where $[\ ,\ ]$ denotes the commutator.
}\\

The proof of those estimates is given in \cite{17r} \cite{18r} with 
$\omega$ replaced by
$<\omega >$ and follows therefrom by a scaling argument. \\

We shall also need the following consequence of Lemma 3.2. \\

\noi {\bf Lemma 3.3.} {\it Let $m \geq 0$ and $1 < r < \infty$. Then 
the following estimate
holds}
\beq
\label{3.6e}
\parallel \omega^m (e^{\varphi} - 1) \parallel_r \ \leq \ \parallel 
\omega^m \varphi
\parallel_{r} \exp \left ( C \parallel \varphi  \parallel_{\infty} 
\right ) \ . \eeq

\vskip 3 truemm
\noi {\bf Proof.} For any integer $n \geq 2$, we estimate
\bea
\label{3.7e}
a_n &\equiv& \parallel \omega^m \ \varphi^n\parallel_r \ \leq C \left 
( \parallel  \omega^m
\varphi \parallel_{r} \
  \parallel \varphi \parallel_{\infty}^{n-1} + \parallel  \omega^m 
\varphi^{n-1} \parallel_{r} \
  \parallel \varphi \parallel_{\infty} \right )  \nn \\
&=& C \left ( a_1 \ b^{n-1} + a_{n-1} \ b \right )
\eea
\noi by (3.4), where $b = \parallel \varphi \parallel_{\infty}$ and 
we can assume $C \geq 1$
without loss of generality. It follows easily from (\ref{3.7e}) that
$$a_n \leq n(Cb)^{n-1} \ a_1$$
\noi for all $n \geq 1$, from which (\ref{3.6e}) follows by expanding 
the exponential.\par
\hfill $\sq$\par

We next give some estimates of $B_1^{t_1}$, $B_S^{t_1}$ and 
$B_L^{t_1}$ defined by
(\ref{2.15e}) (\ref{2.17e}). It follows immediately from (\ref{2.17e}) that
\beq
\label{3.8e}
\parallel \omega^m B_S^{t_1} \parallel_2 \ \leq t^{\beta (m- p)} 
\parallel  \omega^p B_S^{t_1}
\parallel_2 \ \leq t^{\beta (m - p)} \parallel  \omega^p B_1^{t_1} \parallel_2
  \eeq
\noi for $m \leq p$ and similarly
\beq
\label{3.9e}
\parallel \omega^m B_L^{t_1} \parallel_2 \ \leq t^{\beta (m- p)} 
\parallel  \omega^p B_L^{t_1}
\parallel_2 \ \leq t^{\beta (m - p)} \parallel  \omega^p B_1^{t_1} \parallel_2
  \eeq
\noi for $m \geq p$. On the other hand it follows from (\ref{2.15e}) that
\beq
\label{3.10e}
\parallel \omega^{m +1} B_1^{t_1} (w_1, w_2)\parallel_2 \ \leq 
I_m^{t_1} \left ( \parallel
\omega^m (w_1 \bar{w}_2) \parallel_2 \right )
\eeq
\noi where $I_m^{t_1}$ is defined by
\beq
\label{3.11e}
\left ( I_m^{t_1}(f) \right ) (t) = \left | \int_1^{t_1/t} d\nu \ 
\nu^{-m-3/2} \ f(\nu t)
\right | \eeq
\noi or equivalently
$$\left ( I_m^{t_1}(f) \right ) (t) = t^{m+1/2} \left | \int_t^{t_1} 
dt' \ t'^{-m-3/2} \ f(t')
\right | $$
\noi for $t > 0$, $t_1 > 0$. Most of the subsequent estimates of 
$B_1^{t_1}$ will follow from
(\ref{3.8e}) (\ref{3.9e}) (\ref{3.10e}) and from an estimate of $\parallel
\omega^m(w_1\bar{w}_2)\parallel_2$. The latter follows from the HLS 
inequality if $- 3/2 < m <
0$ and from (\ref{3.4e}) if $m \geq 0$. For future reference, we 
quote the following special
case, which will occur repeatedly
\beq
\label{3.12e}
\parallel \omega^{2k - 1/2} B_1^{t_1} (w, w)\parallel_2 \ \leq C \ 
I_{2k-3/2}^{t_1} \left (
\parallel  \omega^k w\parallel_2^2 \right )
  \eeq
\noi and which holds for $0 < k < 3/2$. The required estimate
\beq
\label{3.13e}
\parallel \omega^{2k - 3/2} |w|^2 \parallel_2 \ \leq C  \parallel 
\omega^k w\parallel_2^2
  \eeq
\noi follows from the HLS inequality if $2k < 3/2$ and from 
(\ref{3.4e}) if $2k \geq 3/2$, as
mentioned above, and from Sobolev inequalities. \par

We next give a special estimate of the long range part $B_L^{t_1}$ of 
$B_1^{t_1}$. \\

\noi {\bf Lemma 3.4.} {\it Let $m > - 3/2$. Then}
\beq
\label{3.14e}
\parallel \omega^{m+1}  B_L^{t_1} (w_1, w_2)\parallel_2 \ \leq C \ 
t^{\beta (m + 3/2)}
I_{-3/2}^{t_1} \left ( \parallel  w_1 \parallel_2 \ \parallel  w_2 
\parallel_2\right
) \ .
\eeq
\vskip 3 truemm

\noi {\bf Proof.} Let $f = D_0(\nu ) {\rm Re} \ w_1 \bar{w}_2$. From 
(\ref{2.15e})
(\ref{2.17e}), we estimate
\begin{eqnarray*}
&&\parallel \omega^{m+1}  B_L^{t_1} (w_1, w_2)\parallel_2 \ \leq 
\left | \int_1^{t_1/t} d\nu \
\nu^{-3} \ \parallel  \chi (|\xi | \leq t^{\beta} ) |\xi|^m \ Ff(\xi 
) \parallel_2 \right |
\\
&&\leq \left | \int_1^{t_1/t} d\nu \
\nu^{-3} \ \parallel  \chi (|\xi | \leq t^{\beta} ) |\xi|^m 
\parallel_2 \ \parallel Ff
\parallel_{\infty} \right |
  \\
&&\leq C \left | \int_1^{t_1/t} d\nu \ t^{\beta (m + 3/2)} \parallel 
( w_1 \ \bar{w}_2) (\nu t)
\parallel_{1} \right |
\end{eqnarray*}
\noi which implies (\ref{3.14e}). \par
\hfill $\sq$\par

We finally collect some estimates of the solution of the free wave 
equation $\sq A_0 = 0$ with
initial data $(A_+, \dot{A}_+)$ at time zero, given by (\ref{2.3e}). \\

\noi {\bf Lemma 3.5.} {\it Let $m \geq 0$. Let $\omega^mA_+ \in L^2$, 
$\omega^{m-1} \dot{A}_+
\in L^2$, $\nabla^2 \omega^m A_+ \in L^1$ and $\nabla \omega^m 
\dot{A}_+ \in L^1$. Then the
following estimate holds
\beq
\label{3.15e}
\parallel \omega^{m}  A_0 \parallel_r \ \leq b_0 \ t^{-1+2/r} \quad 
\hbox{for $2 \leq r \leq
\infty$}  \eeq
\noi for all $t > 0$.} \\

\noi {\bf Proof.} It suffices to prove (\ref{3.15e}) for $m = 0$, for 
$r = 2$ and $r = \infty$.
For $r = 2$, it follows from (\ref{2.3e}) that
\beq
\label{3.16e}
\parallel   A_0 \parallel_2 \ \leq \  \parallel   A_+ \parallel_2 + 
\parallel   \omega^{-1}
\dot{A}_+ \parallel_2
\eeq
\noi for all $t \in {I\hskip-1truemm R}$. For $ r = \infty$, the 
result follows from the
divergence theorem applied to the standard representation of 
solutions of the free wave
equation in terms of spherical means \cite{21r}. \par
\hfill $\sq$\par

The time decay expressed by (\ref{3.15e}) is known to be optimal, and 
we shall always consider
solutions $A_0$ of the free wave equation satisfying those estimates 
for suitable $m$. In the
applications, we shall use the estimates (\ref{3.15e}) in the 
equivalent form expressed in
terms of $B_0$ defined by (\ref{2.13e}), namely
\beq
\label{3.17e}
\parallel \omega^m \   B_0 \parallel_r \ \leq \ b_0 \ t^{m-1/r} \quad 
\hbox{for $2 \leq r \leq
\infty$} \ .
\eeq

\mysection{Cauchy problem and preliminary asymptotics for the auxiliary system}
\hspace*{\parindent} In this section, we study the Cauchy problem for 
the auxiliary system
(\ref{2.20e}) and we derive some preliminary asymptotic properties of 
its solutions. This
section illustrates both the method of solution with the help of the 
linearized version
(\ref{2.21e}) of that system and the difficulties arising from the 
different propagation
properties of the Schr\"odinger and wave equations. In particular we 
are able to prove the
existence of solutions up to infinity in time only if $A_0 = 0$. This 
section could be the
starting point for the construction of partial wave operators with 
vanishing asymptotic
states for the field $A$, a construction which would be very similar 
to that of the wave
operators for the Hartree equation performed in \cite{6r} \cite{7r}, 
but which we shall
refrain from performing here. The general case of non-vanishing 
asymptotic states for $A$
will be treated by a similar but more complicated method in Section 6 
below. \par

The basic tool of this section consists of a priori estimates for 
suitably regular
solutions of the linearized system (\ref{2.21e}). Those estimates can 
be proved by a
regularisation and limiting procedure and hold in the integrated form 
at the available
level of regularity. For brevity, we shall state them in differential 
form and we shall
restrict the proof to the formal computation. \par

We first estimate a single solution of the linearized system 
(\ref{2.21e}) at the level of
regularity where we shall eventually solve the auxiliary system 
(\ref{2.20e}). \\

\noi {\bf Lemma 4.1.} {\it Let $1 < k \leq \ell$, $\ell > 3/2$ and 
$\beta > 0$. Let $I \subset
[1, \infty )$ be an interval and let $t_1 \in \bar{I}$. Let $B_0$ 
satisfy the estimates
(\ref{3.17e}) for $0 \leq m \leq k$. Let $(w, s), (w', s') \in {\cal 
C} (I, X^{k,\ell})$ with
$w \in L^{\infty}(I, H^k)$ and let $(w', s')$ be a solution of the 
system (\ref{2.21e}) in $I$.
Then the following estimates hold for all $t \in I$~: \par
$\parallel w'\parallel_2 =$ const.
\bea
\label{4.1e}
\left | \partial_t| w'|_k \right | &\leq& C\ b_0 \left \{ \parallel 
w'\parallel_2^{1/k} \
|w'|_k^{1-1/k} + t^{k-1-\delta /3} \  \parallel w'\parallel_2^{1-\delta /k}\ |w'|_k^{\delta
/k} \right \}\nn \\
&+& C \left \{ t^{-2} |s|_{\ell}^{\dot{}} + t^{-1-\beta_1} \ I_{m_1}^{t_1} (|w|_k^2) \right \}
|w'|_k
\eea
\beq
\label{4.2e}
\left | \partial_t| s'|_{\ell}^{\dot{}} \right | \leq C\ t^{-2} 
|s|_{\ell}^{\dot{}}
|s'|_{\ell}^{\dot{}} + C\ t^{-1+\beta_2} \Big ( I_{m_1-k}^{t_1} (\parallel
w\parallel_2 |w|_k) + I_{m_1}^{t_1}(|w|_k^2 ) \Big )
\eeq
\noi where $0 < \delta \leq [k \wedge 3/2]$,
\bea
\label{4.3e}
&&\beta_1 = \beta [1 \wedge 2(k - 1)] = \beta (1 - 2 [3/2 - k]_+) \ , \\
\label{4.4e}
&&m_1 = [k \wedge (2k - 3/2)] = k - [3/2 - k]_+ \ , \\
&&\beta_2  = \beta (\ell + 1 - k + [3/2 - k]_+) \ .
\label{4.5e}
\eea
}
\vskip 3 truemm

\noi {\bf Proof.} We omit the superscript $t_1$ in all the proof. We 
first estimate $w'$. It is
clear from (\ref{2.21e}) that $\parallel w'\parallel_2$ = const. We 
next estimate
\bea
\label{4.6e}
&&\left | \partial_t \parallel \omega^k w'\parallel_2 \right | \leq 
t^{-1} \parallel
[\omega^k , B_0] w' \parallel_2 + t^{-2} \Big \{  \parallel [\omega^k, s]
\cdot \nabla w'\parallel_2 + \parallel (\nabla \cdot s) \omega^k w' 
\parallel_2 \nn \\
&&+\parallel \omega^k ((\nabla \cdot s)w')\parallel_2 \Big \} + 
t^{-1} \parallel [\omega^k ,
B_s(w, w)]w' \parallel_2  \ .
\eea
\noi The contribution of $B_0$ is estimated by Lemma 3.2 and (\ref{3.17e}) as
  \bea
&&\parallel [\omega^k , B_0] w' \parallel_2\ \leq C \left ( 
\parallel \nabla B_0
\parallel_{\infty} \ \parallel \omega^{k-1} w' \parallel_2 + 
\parallel \omega^k  B_0
\parallel_{3/\delta} \ \parallel  w' \parallel_r \right ) \nn \\
&&\leq C \ b_0 \left ( t \parallel \omega^{k-1} w' \parallel_2 + t^{k 
- \delta /3} \parallel  w'
\parallel_r \right )
\label{4.7e}
\eea
\noi with $0 < \delta = \delta (r) < k \wedge 3/2$. This yields the 
first term in the RHS
of (\ref{4.1e}) by Sobolev inequalities and interpolation. We next 
estimate by Lemma 3.2
\begin{eqnarray*}
&&\parallel [\omega^k , s] \cdot \nabla w' \parallel_2\ + \ \parallel 
(\nabla \cdot s)
\omega^k w'  \parallel_{2} \ + \ \parallel \omega^{k} (( \nabla \cdot 
s) w') \parallel_2 \\
&&\leq C \left (  \parallel \nabla s \parallel_{\infty} \ \parallel 
\omega^{k} w' \parallel_2
+ \parallel \omega^k  s  \parallel_{3/\delta} \ \parallel \nabla w' 
\parallel_r + \parallel \omega^k
(\nabla \cdot s)  \parallel_{3/\delta ' }\ \parallel  w'
\parallel_{r'}\right )
\end{eqnarray*}
\noi where $0 < \delta = \delta (r) \leq [(k - 1) \wedge 3/2]$ and $0 
< \delta ' = \delta (r')
\leq [k \wedge 3/2]$. Choosing $\delta = [(k - 1) \wedge 1/2]$ and 
$\delta ' = [k \wedge 3/2]$
and using Sobolev inequalities, we continue the previous estimate by
\bea
\label{4.8e}
&&\cdots \leq C \Big (  \parallel \nabla s \parallel_{\infty} \ 
\parallel \omega^{k} w'
\parallel_2 + \parallel \omega^{[k \vee 3/2]} \nabla  s 
\parallel_{2} \ \parallel \omega^{[k
\wedge 3/2]} w' \parallel_{2}\nn \\
&&+ \chi (k > 3/2) \parallel \omega^{k} \nabla  s  \parallel_{2}\ \parallel
w'\parallel_{\infty} \Big ) \leq C \ |s|_{\ell}^{\dot{}} \ |w'|_k \ .   \eea
\noi We next estimate the contribution of $B_S$ to (\ref{4.6e}). By 
Lemma 3.2 and Sobolev
inequalities, we estimate
\bea
&&\parallel [\omega^k , B_S(w, w)] w' \parallel_2 \nn \\
&&\leq C \left \{ \parallel \nabla B_S(w,w) \parallel_3 \ \parallel 
\omega^{k-1} w' \parallel_6
+ \parallel \omega^k B_S(w,w)\parallel_{3/ \delta} \ \parallel 
w'\parallel_r \right \} \nn \\
&&\leq C \left \{ \parallel \omega^{3/2} B_S(w,w) \parallel_2 \ 
\parallel \omega^{k} w'
\parallel_2 + \parallel \omega^{k+ 3/2 - \delta} B_S(w,w)\parallel_2 
\ \parallel
w'\parallel_r \right \}   \label{4.9e}  \eea

\noi where $0 < \delta = \delta (r) \leq 3/2$. We choose $\delta = [k 
\wedge 3/2]$ and continue
(\ref{4.9e}) as follows~: \par

If $k < 3/2$, so that $\delta = k$,
\bea
\label{4.10e}
\cdots &\leq& C \parallel \omega^{3/2} B_S(w,w) \parallel_2\ \parallel
\omega^k w'  \parallel_{2} \nn \\
&\leq& C \ t^{-2\beta (k - 1)} \parallel \omega^{2k-1/2} B_1(w,w) 
\parallel_2\ \parallel
\omega^k w'  \parallel_{2} \nn \\
&\leq& C \ t^{-2\beta (k - 1)} I_{2k-3/2} \left ( \parallel \omega^k 
w \parallel_2^2\right )
\parallel \omega^k w'  \parallel_{2} \nn \\
&\leq& C \ t^{-\beta_1} I_{m_1} \left ( |w|_k^2 \right ) |w'|_k
   \eea
\noi by Sobolev inequalities, by (\ref{3.8e}) (\ref{3.12e}) and by 
the definitions (\ref{4.3e})
(\ref{4.4e}). \par

If $k = 3/2$, so that $\delta = 3/2 - \varepsilon$,
\bea
\label{4.11e}
\cdots &\leq& C \parallel \omega^{3/2+ \varepsilon} B_S(w,w) 
\parallel_2\ \parallel
\omega^{3/2 - \varepsilon} w'  \parallel_{2} \nn \\
&\leq& C \ t^{-\beta (1 - 2\varepsilon)} \parallel \omega^{5/2 - 
\varepsilon} B_1(w,w)
\parallel_2\ \parallel \omega^{3/2 - \varepsilon} w'  \parallel_{2} \nn \\
&\leq& C \ t^{-\beta (1 - 2\varepsilon)} I_{3/2- \varepsilon} \Big ( 
\parallel \omega^{(3-
\varepsilon)/2} w \parallel_2^2 \Big ) \parallel \omega^{3/2 - 
\varepsilon} w'  \parallel_{2}
\nn \\ &\leq& C \ t^{-\beta_1} I_{m_1} \left ( |w|_k^2 \right ) |w'|_k
   \eea
\noi by Sobolev inequalities, by (\ref{3.8e}), by (\ref{3.12e}) with 
$k = (3 - \varepsilon )/2$
and by (\ref{4.3e}) (\ref{4.4e}). \par

If $k > 3/2$, so that $\delta = 3/2$ and $r = \infty$,
\bea
\label{4.12e}
\cdots &\leq& C \parallel \omega^{k+1} B_1(w,w) \parallel_2\ \left \{ 
t^{-\beta (k -
1/2)} \parallel \omega^{k} w'  \parallel_{2} + t^{-\beta} \parallel 
w'  \parallel_{\infty}
\right \}\nn \\
&\leq& C \ t^{-\beta} \ I_k \left ( \parallel \omega^{k}w 
\parallel_2\ \parallel
w\parallel_{\infty} \right ) \left (  \parallel \omega^{k}w' 
\parallel_2 + \parallel
w'\parallel_{\infty} \right ) \nn \\
&\leq& C \
t^{-\beta_1} I_{m_1} \left ( |w|_k^2 \right ) |w'|_k
   \eea
\noi by (\ref{3.8e}) (\ref{3.10e}), Lemma 3.2, Sobolev inequalities 
and (\ref{4.3e})
(\ref{4.4e}). Substituting (\ref{4.7e}) (\ref{4.8e}) (\ref{4.10e}) 
(\ref{4.11e}) (\ref{4.12e})
into (\ref{4.6e}) yields (\ref{4.1e}).\par

We now turn to the estimate of $s'$, namely to the proof of 
(\ref{4.2e}). For $0 \leq m \leq
\ell$, we estimate
\bea
\label{4.13e}
\left | \partial_t \parallel \omega^{m+1} s'\parallel_2 \right | 
&\leq&t^{-2} \left \{
\parallel [\omega^{m+1} , s] \cdot \nabla s' \parallel_2 + \parallel 
(\nabla \cdot s)
\omega^{m+1}\  s'\parallel_2 \right \} \nn \\
&&+ t^{-1} \parallel \omega^{m+2} B_L(w, w) \parallel_2  \ .
\eea
\noi The first bracket in the RHS of (\ref{4.13e}) is estimated by Lemma 3.2 as
\bea
\label{4.14e}
\{\cdot \} &\leq& C \left ( \parallel \nabla s\parallel_{\infty} \ 
\parallel \omega^{m+1} s'
\parallel_2 + \parallel \omega^{m+1} s\parallel_2 \ \parallel \nabla 
s'\parallel_{\infty}\right )
\nn \\  &\leq& C  |s|_{\ell}^{\dot{}} \ |s'|_{\ell}^{\dot{}}
\eea
\noi by Sobolev inequalities. \par

The contribution of $B_L$ for $m = \ell$ is estimated by (\ref{3.9e}) 
(\ref{3.10e})
(\ref{3.12e}) and Lemma 3.2 as
\bea
\label{4.15e}
\parallel \omega^{\ell + 2} B_L(w,w)\parallel_2  &\leq& C \ t^{\beta 
(\ell + 5/2 - 2k)}\ I_{2k -
3/2} \left (  \parallel \omega^k w\parallel_2^2 \right ) \qquad \ 
\hbox{for $k < 3/2$} \ , \nn \\
&&C \ t^{\beta (\ell + 1 - k)}\ I_{k} \left (  \parallel \omega^k 
w\parallel_2 \ \parallel
w\parallel_{\infty} \right ) \qquad \hbox{for $k > 3/2$} \ , \nn \\
&&C \ t^{\beta (\ell - 1/2 + \varepsilon)}\ I_{3/2 - \varepsilon} 
\left (  \parallel \omega^{(3
- \varepsilon )/2}  w\parallel_2^2 \right ) \quad \hbox{for $k = 3/2$} \
, \nn \\
&&\leq C \ t^{\beta_2} \ I_{m_1} \left (|w|_k^2 \right )
   \eea
\noi in all cases. The contribution of $B_L$ for $m = 0$ is estimated 
similarly as
\bea
\label{4.16e}
\parallel \nabla^2 B_L(w,w)\parallel_2 &\leq& C \ t^{\beta (5/2 - 
k)}\ I_{k- 3/2} \ \left (
\parallel \omega^k w\parallel_2  \ \parallel w\parallel_2\right ) 
\quad \ \hbox{for $k < 3/2$} \
, \nn \\  &&C \ t^{\beta }\ I_0 \left (  \parallel 
w\parallel_{\infty}  \ \parallel w\parallel_2
\right ) \qquad \qquad \qquad \quad \hbox{for $k > 3/2$} \ , \nn \\
&&C \ t^{\beta (1 + \varepsilon)}\ I_{- \varepsilon} \left ( 
\parallel \omega^{3/2 -
\varepsilon}  w\parallel_2 \  \parallel w\parallel_2 \right ) \qquad 
\hbox{for $k = 3/2$} \ ,
\nn \\
&&\leq C \ t^{\beta '_2} \ I_{m_1-k} \left (\parallel w \parallel_2 \ 
|w|_k \right )
   \eea
\noi in all cases, with
\beq
\label{4.17e}
\beta '_2 = \beta \left ( 1 + [3/2 - k]_+ \right ) \leq \beta_2
\eeq
\noi since $\ell \geq k$.

Collecting (\ref{4.14e}) (\ref{4.15e}) (\ref{4.16e}) yields (\ref{4.2e}).\par
\hfill $\sq$\par

We next estimate the difference of two solutions of the linearized 
system (\ref{2.21e})
corresponding to two different choices of $(w, s)$. We estimate that 
difference at a lower
level of regularity than the solutions themselves. \\

\noi {\bf Lemma 4.2.} {\it Let $1 < k \leq \ell$, $\ell > 3/2$ and 
$\beta > 0$. Let $I \subset
[1, \infty )$ be an interval and let $t_1 \in \bar{I}$. Let $B_0$ be 
sufficiently regular, for
instance $B_0 \in {\cal C}(I, H_3^k)$. Let $(w_i, s_i), (w'_i, s'_i) 
\in {\cal C} (I,
X^{k,\ell})$ with $w_i \in L^{\infty}(I, H^k)$, $i = 1,2$, and let 
$(w'_i, s'_i)$ be solutions
of the system (\ref{2.21e}) associated with $(w_i, s_i)$. Define 
$(w_{\pm}, s_{\pm}) = (1/2)
(w_1 \pm w_2, s_1 \pm s_2)$ and $(w'_{\pm} , s'_{\pm}) = 1/2 (w'_1 
\pm w'_2, s'_1 \pm s'_2)$.
Then the following estimates hold for all $t \in I$~:
\beq    \label{4.18e}
\left | \partial_t\parallel w'_-\parallel_2  \right | \leq C\ t^{-2}
|s_-|_{\ell_0}^{\dot{}}\ |w'_+|_k + C t^{-1- \beta_1} \ I_{m_1 - 
k}^{t_1} \left ( |w_+|_k
\parallel w_-\parallel_2 \right ) |w'_+|_k
\eeq
\bea
\label{4.19e}
\left | \partial_t | s'_-|_{\ell_0}  \right | &\leq& C\ t^{-2} \left 
( |s_+|_{\ell}^{\dot{}}\
|s'_-|_{\ell_0}^{\dot{}}  + |s_-|_{\ell_0}^{\dot{}}\ 
|s'_+|_{\ell}^{\dot{}}  \right ) \nn \\
&+& C \ t^{-1 + \beta_2} \ I_{m_1-k}^{t_1} \left ( |w_+|_k \ 
\parallel w_-\parallel_2 \right )
\eea
\noi where $\beta_1$, $m_1$ and $\beta_2$ are defined by (\ref{4.3e}) 
(\ref{4.4e}) (\ref{4.5e})
and where
\beq
\label{4.20e}
[3/2 - k]_+ \leq \ell_0 \leq \ell - k \ .
\eeq
}
\vskip 3 truemm

\noi {\bf Proof.} We again omit the superscript $t_1$ in the proof. 
Taking the difference of
the system (\ref{2.21e}) for $(w'_i, s'_i)$, we obtain the following 
system for $(w'_-, s'_-)$~:
\bea
\label{4.21e}
\left \{ \begin{array}{l} \partial_t w'_- = i (2t^2)^{-1} \Delta w'_- 
+ t^{-2} (Q(s_+, w'_-) +
Q(s_-, w'_+)) + it^{-1} B_0 w'_-  \\
\\
\qquad + i t^{-1} \left \{ \left ( B_S(w_+, w_+) + B_S(w_-,w_-) 
\right ) w'_- + 2B_S(w_+,
w_-)w'_+ \right \} \\
\\
\partial_t s'_- = t^{-2} \left ( s_+ \cdot \nabla s '_- + s_- \cdot 
\nabla s'_+ \right ) -
2t^{-1}\nabla B_L(w_+, w_-) \ .
  \end{array}
\right .
\eea
\noi We first estimate $w'_-$. From (\ref{4.21e}) we obtain
\beq
\label{4.22e}
\left | \partial_t \parallel w'_- \parallel_2 \right | \leq t^{-2} 
\parallel Q( s_-,w'_+)
\parallel_2 + 2t^{-1}\parallel B_S(w_+, w_-) w'_+ \parallel_2
\eeq
\noi where only those terms appear that do not preserve the 
$L^2$-norm. We estimate the first
norm in the RHS by H\"older and Sobolev inequalities as follows~: \par

If $k < 3/2$,
\begin{eqnarray*}
\parallel Q(s_-, w'_+) \parallel_2 \ &\leq& C \left ( \parallel 
s_-\parallel_{3/(k-1)} +
\parallel \nabla \cdot s_- \parallel_{3/k} \right ) \parallel 
\omega^k w'_+\parallel_2 \\
&\leq& C \parallel \omega^{3/2- k} \nabla s_- \parallel \ \parallel 
\omega^k w'\parallel_2 \ .
\end{eqnarray*}

If $k = 3/2$,
$$\parallel Q(s_-, w'_+) \parallel_2 \ \leq C  \parallel 
\omega^{\varepsilon} \nabla s_-
\parallel_{2} \ \parallel \omega^{3/2 - \varepsilon}\ w'_+ \parallel_2 \ .$$

If $k > 3/2$,
$$\parallel Q(s_-, w'_+) \parallel_2 \ \leq C  \parallel  \nabla s_-
\parallel_{2} \left ( \parallel \nabla w'_+\parallel_3 + \parallel 
w'_+\parallel_{\infty}
\right ) \ ,$$ \noi and in all cases
\beq
\label{4.23e}
\parallel Q(s_-, w'_+) \parallel_2 \ \leq C |s_-|_{\ell_0}^{\dot{}} \ |w'_+|_k
\eeq
\noi provided $\ell_0 \geq [3/2 - k]_+$. \par

We estimate the second norm in the RHS of (\ref{4.22e}) by 
(\ref{3.8e}) (\ref{3.10e}), by Lemma
3.2 and by the H\"older and Sobolev inequalities as follows~: \par

If $k < 3/2$,
\begin{eqnarray*}
&&\parallel B_S(w_+, w_-)w'_+ \parallel_2 \ \leq C \parallel 
\omega^{3/2-k} \ B_S(w_+,
w_-)\parallel_2 \ \parallel \omega^k w'_+ \parallel_2 \\
&&\leq C\ t^{-2\beta (k-1)} \parallel \omega^{k-1/2} \ B_1(w_+, 
w_-)\parallel_2 \ \parallel
\omega^k w'_+ \parallel_2 \\
&&\leq C\ t^{-2\beta (k-1)}\ I_{k-3/2} \left ( \parallel \omega^k w_+ 
\parallel_2 \ \parallel
w_-\parallel_2 \right ) \parallel \omega^k w'_+ \parallel_2 \ .
  \end{eqnarray*}

If $k = 3/2$,
$$\parallel B_S(w_+, w_-)w'_+ \parallel_2 \ \leq C\ t^{-\beta (1 - 2 
\varepsilon )} \
I_{-\varepsilon} \left ( \parallel \omega^{3/2 - \varepsilon} w_+ 
\parallel_2 \ \parallel
w_-\parallel_2 \right ) \parallel \omega^{3/2 - \varepsilon}  w'_+ 
\parallel_2 \ .$$

If $k > 3/2$,
\begin{eqnarray*}
&&\parallel B_S(w_+, w_-)w'_+ \parallel_2 \ \leq t^{-\beta} \parallel 
\nabla B_S
(w_+, w_-)\parallel_2\ \parallel w'_+ \parallel_{\infty} \\
&&\leq C\ t^{-\beta} \ I_{0} \left ( \parallel w_+ \parallel_{\infty} 
\ \parallel w_-\parallel_2
\right ) \parallel w'_+ \parallel_{\infty}
\end{eqnarray*}
\noi and in all cases
\beq
\label{4.24e}
\parallel B_S(w_+, w_-)w'_+ \parallel_2 \ \leq C\ t^{-\beta_1} 
I_{m_1-k} \left ( | w_+
|_k \parallel w_-\parallel_2 \right ) |w'_+|_k \ .
\eeq
\noi with $\beta_1$ and $m_1$ defined by (\ref{4.3e}) (\ref{4.4e}). 
Substituting (\ref{4.23e})
(\ref{4.24e}) into (\ref{4.22e}) yields (\ref{4.18e}). \par

We now turn to the estimate of $s'_-$, namely to the proof of 
(\ref{4.19e}). From (\ref{4.21e})
we estimate for $m \geq 0$
\bea
\label{4.25e}
&&\partial_t \parallel \omega^{m+1} s'_-\parallel_2 \ \leq t^{-2} 
\Big \{ \parallel
[\omega^{m+1}, s_+]\cdot \nabla s'_-\parallel_2 + \parallel (\nabla 
\cdot s_+) \omega^{m+1} s'_-
\parallel_2 \nn \\
&& + \parallel \omega^{m+1}(s_-\cdot \nabla s'_+)\parallel_2 \Big \} 
+ 2t^{-1}\parallel
\omega^{m+2} B_L(w_+, w_-) \parallel_2 \ .
\eea

If $m = 0$ (a case which has to be self-estimating if $\ell_0 = 0$, 
which is allowed if $k
> 3/2)$, we estimate the bracket in the RHS of (\ref{4.25e}) directly as
\bea
\label{4.26e}
&&\{ (m = 0)\} \leq \ \parallel \nabla s_+\parallel_{\infty}\ 
\parallel \nabla s'_- \parallel_2
+ \parallel \nabla s_- \parallel_2\left ( \parallel \nabla 
s'_+\parallel_{\infty} +  \parallel
\nabla^2 s'_+ \parallel_3 \right ) \nn \\
&&\leq C \left ( |s_+|_{\ell}^{\dot{}} \parallel \nabla s'_- \parallel_2 +
\parallel \nabla s_- \parallel_2 \ |s'_+|_{\ell}^{\dot{}} \right )
\eea
\noi since $\ell > 3/2$. \par

If $m > 0$, we estimate that bracket by Lemma 3.2 and Sobolev inequalities as
\bea
\label{4.27e}
&&\{ \cdot \} \leq C \Big \{ \parallel \nabla s_+\parallel_{\infty} \ 
\parallel \omega^{m+1} s'_-
\parallel_2 + \parallel \omega^{m+1}  s_+ \parallel_{3 / \delta} \ 
\parallel \nabla
s'_-\parallel_{r}\nn \\
&& +  \parallel \omega^{m+1} s_- \parallel_2 \  \parallel \nabla s'_+ 
\parallel_{\infty} +
\parallel s_- \parallel_{r'} \ \parallel \omega^{m+2}  s'_+ 
\parallel_{3 / \delta '} \Big \}
\eea
\noi where $0 < \delta = \delta (r) \leq 3/2$, $0 < \delta ' = \delta 
(r') \leq 3/2$. The first
and third term in the RHS of (\ref{4.27e}) are readily controlled by 
the corresponding terms in
(\ref{4.19e}) for $0 < m \leq \ell_0$ and $\ell > 3/2$. The remaining 
two terms are similarly
controlled through Sobolev inequalities provided
\begin{eqnarray*}
&&0 < \delta \leq [\ell_0 \wedge 3/2] \quad , \quad m + 3/2 - \delta 
\leq \ell \ . \\
&&1 \leq \delta ' \leq [( \ell_0 + 1) \wedge 3/2] \quad , \quad m + 
5/2 - \delta ' \leq \ell \ .
\end{eqnarray*}
\noi Those conditions are easily seen to be compatible in $\delta$ 
and $\delta '$ for all $m$,
$0 < m \leq \ell_0$, provided $\ell \geq [(\ell_0 + 1) \vee 3/2]$, 
which follows from $\ell >
3/2$ and $\ell \geq \ell_0 + k$. \par

We finally estimate the contribution of $B_L(w_+, w_-)$ by
$$\parallel \omega^{m+2} B_L(w_+, w_-) \parallel_2 \ \leq C\ t^{\beta 
m}  \parallel \nabla^2 B_L
(w_+, w_-) \parallel_2$$
\noi by (\ref{3.9e}) and we estimate the last norm in exactly the 
same way as in (\ref{4.16e}),
thereby obtaining
\beq
\label{4.28e}
\parallel \omega^{m+2} B_L(w_+, w_-) \parallel_2 \ \leq C\ t^{\beta 
'_2 + \beta m}  \
I_{m_1-k} \left ( |w_+|_k \parallel  w_- \parallel_2 \right ) \ .
\eeq
\noi Collecting (\ref{4.25e}) (\ref{4.26e}), (\ref{4.27e}) and the 
discussion that follows, and
(\ref{4.28e}) and noting that $\beta '_2 + \beta m \leq \beta_2$ for 
$m \leq \ell_0 \leq \ell -
k$, we obtain (\ref{4.19e}).\par \nobreak
\hfill $\sq$\par

With the estimates of Lemma 4.1 and 4.2 available, it is an easy 
matter to solve the Cauchy
problem globally in time for the linearized system (\ref{2.21e}). \\

\noi {\bf Proposition 4.1.} {\it Let $1 < k \leq \ell$, $\ell > 3/2$ 
and $\beta > 0$. Let $I
\subset [1, \infty )$ be an interval and let $t_1 \in \bar{I}$. Let 
$B_0$ satisfy the estimates
(\ref{3.17e}) for $0 \leq m \leq k$. Let $(w, s) \in {\cal C} (I, 
X^{k,\ell})$ with $w \in
L^{\infty}(I, H^k)$. Let $t_0 \in I$ and let $(w'_0, s'_0) \in 
X^{k,\ell}$. Then the system
(\ref{2.21e}) has a unique solution $(w', s') \in {\cal C} (I, 
X^{k,\ell})$ with $(w', s')(t_0)
= (w'_0, s'_0)$. That solution satisfies the estimates (\ref{4.1e}) 
(\ref{4.2e}) for all $t
\in I$. Two such solutions $(w'_i, s'_i)$ associated with $(w_i, 
s_i)$, $i = 1,2$, satisfy
the estimates (\ref{4.18e}) (\ref{4.19e}) for all $t \in I$.} \\

\noi {\bf Proof.} The proof proceeds in the same way as that of 
Proposition 4.1 of \cite{6r},
through a parabolic regularization and a limiting procedure, with the 
simplification that the
system (\ref{2.21e}) is linear. We define $U_1(t) = U(1/t)$, 
$\widetilde{w}'(t) = U_1(t)w'(t)$.
We first consider the case $t \geq t_0$. The system (\ref{2.21e}) 
with a parabolic
regularization added is rewritten in terms of the variables 
$(\widetilde{w}', s')$ as
\bea
\label{4.29e}
\left \{ \begin{array}{ll} \partial_t \widetilde{w}' &= \eta \Delta 
\widetilde{w}' +
t^{-2} U_1 Q(s, U_1^* \widetilde{w}') + i t^{-1} U_1 (B_0 + B_S(w,w)) 
U_1^* \widetilde{w}'
\\
\\
&\equiv \eta \Delta \widetilde{w}' + F(\widetilde{w}') \\
\\
\partial_t s ' &= \eta \Delta s' + t^{-2} s \cdot \nabla s' - t^{-1} 
\nabla B_L (w, w) \equiv
\eta \Delta s' + G(s')
  \end{array}
\right .
\eea
\noi where the parametric dependence of $F$, $G$ on $(w, s)$ has been 
omitted. The Cauchy
problem for the system (\ref{4.29e}) can be recast in the integral form
\beq
\label{4.30e}
{\widetilde{w}' \choose s'} (t) = V_{\eta}(t-t_0) {\widetilde{w}'_0 
\choose s'_0} +
\int_{t_0}^t dt' \ V_{\eta} (t-t') {F(\widetilde{w}') \choose G(s')} (t') \eeq
\noi where $V_{\eta} (t) = \exp (\eta t \Delta)$. The operator 
$V_{\eta}(t)$ is a contraction
in $X^{k,\ell}$ and satisfies the bound
$$\parallel \nabla V_{\eta}(t) ; {\cal L} (X^{k,\ell}) \parallel \ 
\leq C(\eta t)^{-1/2} \ .$$
\noi From those facts and from estimates on $F$, $G$ similar to and 
mostly contained in those
of Lemma 4.1, it follows by a contraction argument that the system 
(\ref{4.30e}) has a unique
solution $(\widetilde{w}'_{\eta}, s'_{\eta}) \in {\cal C}([t_0, t_0 + 
T], X^{k,\ell})$ for
some $T > 0$ depending only on $|w'_0|_k$, $|s'_0|_{\ell}^{\dot{}}$ 
and $\eta$. That solution
satisfies the estimates (\ref{4.1e}) and (\ref{4.2e}) and can 
therefore be extended to $I_+ =
I \cap \{ t : t \geq t_0\}$ by a standard globalisation argument 
using Gronwall's inequality.
\par

We next take the limit $\eta \to 0$. Let $\eta_1$, $\eta_2 > 0$ and 
let $(w'_i, s'_i) =
(w'_{\eta_i}, s'_{\eta_i})$, $i = 1,2$ be the corresponding 
solutions. Let $(w'_-, s'_-) =
(1/2)(w'_1 - w'_2, s'_1 - s'_2)$. By estimates similar to, but 
simpler than those of Lemma 4.2,
since in particular $(w_-, s_-) = 0$, we obtain
$$\left \{ \begin{array}{l} \partial_t \parallel w'_-\parallel_2^2 \ \leq
|\eta_1 - \eta_2 | \left ( \parallel \nabla w'_1 \parallel_2^2 + \parallel
\nabla w'_2 \parallel_2^2 \right ) \\
\\
\partial_t \parallel \nabla s'_-\parallel_2^2 \ \leq
|\eta_1 - \eta_2 | \left ( \parallel \nabla^2 s'_1 \parallel_2^2 + \parallel
\nabla^2 s'_2 \parallel_2^2\right ) + C\ t^{-2} \parallel \nabla s_+ 
\parallel_{\infty} \parallel
\nabla s'_-\parallel_2^2 \ .\end{array}\right .$$
\noi Those estimates imply that $(w'_{\eta}, s'_{\eta})$ converges in 
$X^{0,0}$ uniformly in
time in the compact subintervals of $I_+$, to a solution of the 
original system. It follows
then by a standard compactness argument using the estimates 
(\ref{4.1e}) (\ref{4.2e}) that the
limit belongs to ${\cal C}(I_+, X^{k,\ell})$. This completes the 
proof for $t \geq t_0$. The
case $t \leq t_0$ is treated similarly. \par \hfill
$\sq$\par

We now turn to the Cauchy problem for the auxiliary system 
(\ref{2.20e}). Because of the
difficulties described in Section 2, the problem of existence of 
solutions is scattered
with pitfalls, as the discussion below will show. On the other hand, 
the uniqueness problem of
suitably bounded solutions is a rather easy matter and we consider 
that problem first. The
proof relies entirely on Lemma 4.2 and therefore does not require any 
a priori estimate on
$B_0$. The snag of course is that it is difficult to prove the 
existence of solutions with the
required boundedness properties. \\

\noi {\bf Proposition 4.2.} {\it Let $1 < k \leq \ell$, $\ell > 3/2$ 
and $\beta > 0$. Let $I
\subset [1, \infty )$ be an interval and let $t_1 \in \bar{I}$. Let 
$B_0$ be sufficiently
regular, for instance $B_0 \in {\cal C}(I, H_3^k)$. \par

(1) Let $t_0 = t_1 < \infty$ and let $(w_0, s_0) \in X^{k,\ell}$. 
Then the system (\ref{2.20e})
has at most one solution $(w, s) \in {\cal C}(I, X^{k,\ell})$ with $w 
\in L^{\infty}(I, H^k)$
and $(w, s)(t_0) = (w_0, s_0)$. \par

Let now $\beta_2 < 1$, where $\beta_2$ is defined by (\ref{4.5e}), 
let $\ell_0$ satisfy
(\ref{4.20e}). Let $(w_i, s_i)$, $i = 1,2$ be two solutions of the 
system (\ref{2.20e}) in
$I$ such that $(w_i, t^{\eta - 1} s_i) \in ({\cal C} \cap L^{\infty}) 
(I, X^{k,\ell})$ for
some $\eta > 0$ and let
\beq
\label{4.31e}
\parallel w_i, L^{\infty}(H^k)\parallel \ \leq a \qquad , \qquad 
\parallel t^{\eta - 1} \nabla
s_i ; L^{\infty}(H^{\ell})\parallel \ \leq b \ .
\eeq

(2) Let $t_0 \in I$, $t_0 \leq t_1$, $t_0 < \infty$ and assume that 
$(w_1, s_1)(t_0) = (w_2,
s_2)(t_0)$. Then there exists $c = c(a,b)$ such that if
\beq
\label{4.32e}
\left ( t_0^{-(1 - \beta_2)} \vee t_0^{-\beta_1} \right ) \left ( 1 -
(t_0/t_1)^{\alpha}\right ) \leq c(a,b) \eeq
\noi where $\alpha = [k \wedge 3/2]- 1$, then $(w_1, s_1) = (w_2, 
s_2)$. In particular there
exists $T_0 = T_0 (a,b)$ such that if $t_0 \geq T_0$, then $(w_1, 
s_1) = (w_2, s_2)$. \par

(3) Let $t_1 = \infty$. Assume that $\parallel w_1 - w_2\parallel_2 
t^{\beta_2}$ and $|s_1 -
s_2|_{\ell_0}^{\dot{}}$ tend to zero when $t \to \infty$. Then $(w_1, 
s_1) = (w_2, s_2)$.}
\\

\noi {\bf Proof.} If $(w_i, s_i)$, $i = 1,2$ are two solutions of the 
system (\ref{2.20e}) in
${\cal C}(I, X^{k,\ell})$, then they satisfy the estimates 
(\ref{4.18e}) (\ref{4.19e}) with
$(w'_i, s'_i) = (w_i, s_i)$, which we denote (\ref{4.18e}=) 
(\ref{4.19e}=) and refrain from
rewriting for brevity. The proof consists in exploiting those 
estimates to prove that $(w_1,
s_1) = (w_2, s_2)$. We define $y = \parallel w_-\parallel_2$ and $z = 
|s_-|_{\ell_0}^{\dot{}}$.
\par

\noi \underbar{Part (1)}. With $t_0 = t_1 < \infty$, the estimates 
(\ref{4.18e}=)
(\ref{4.19e}=) take the general form
\bea
\label{4.33e}
&&|\partial_t y| \leq f_1(t) z + g_1(t) \int_{t_0}^t dt' \ h_1(t') \ y(t') \\
&&|\partial_t z| \leq f_2(t) z + g_2(t) \int_{t_0}^t dt' \ h_2(t') \ y(t')
\label{4.34e}
\eea
\noi for suitable continuous nonnegative functions $f_1$, $g_1$, 
$h_1$, $f_2$, $g_2$, $h_2$
(actually $h_1 = h_2$, but that is irrelevant). Furthermore $y(t_0) = 
z(t_0) = 0$. We shall
reduce the system (\ref{4.33e}) (\ref{4.34e}) to a standard form 
where Gronwall's inequality is
applicable. We restrict our attention to the case $t \geq t_0$ for 
definiteness. The case $t
\leq t_0$ can be treated similary. Defining $\widetilde{z}$ by
$$z(t) = E(t) \ \widetilde{z}(t) = \exp \left \{ \int_{t_0}^t dt'\ 
f_2(t') \right \}
\widetilde{z}(t) \ ,$$
\noi we reduce the system (\ref{4.33e}) (\ref{4.34e}) for $(y, z)$ to 
a similar system for $(y,
\widetilde{z})$, where $f_2$, $g_2$ and $f_1$ are replaced by 0, 
$E^{-1} g_2$ and $Ef_1$. We can
therefore assume that $f_2 = 0$. Then
$$z(t) \leq \int_{t_0}^t dt'' \ g_2(t'') \int_{t_0}^{t''} dt' \ 
h_2(t') \ y(t') \leq G_2(t)
\int_{t_0}^t dt' \ h_2(t') \ y(t')   $$
\noi where
$$G_2(t) = \int_{t_0}^t dt'' \ g_2(t'')$$
\noi so that
\begin{eqnarray*}
\partial_t y &\leq& f_1(t) G_2(t) \int_{t_0}^t dt' \ h_2(t') y(t') + 
g_1(t) \int_{t_0}^{t} dt'
\ h_1(t') \ y(t') \\
&\leq& ( f_1\ G_2 + g_1) \int_{t_0}^t dt' (h_1 \vee h_2)(t') \ y(t')
\end{eqnarray*}
\noi which is of the same form as (\ref{4.33e}) with $f_1 = 0$. 
Integrating the latter yields
$$
y \leq \int_{t_0}^t dt'' \ g_1(t'') \int_{t_0}^{t''} dt' \ h_1(t') \ y(t')
\leq G_1(t) \int_{t_0}^t dt' \ h_1(t') \ y(t')
$$
\noi where
$$G_1(t) = \int_{t_0}^t dt' \ g_1(t') \ , $$
\noi which together with $y(t_0) = 0$ implies $y(t) = 0$ for all $t$ 
by an easy variant of
Gronwall's inequality. Substituting that result into (\ref{4.34e}) 
(with $f_2 = 0$) yields $z =
0$ and therefore $(w_1, s_1) = (w_2, s_2)$. \par

We now turn to the proof of Parts (2) and (3). Introducing the 
assumption and notation
(\ref{4.31e}), changing the variable from $\nu$ to $t' = \nu t$ in 
the definition of
$I_m^{t_1}$, and omitting an absolute overall constant, we can 
rewrite (\ref{4.18e}=)
(\ref{4.19e}=) in the form
\bea
\label{4.35e}
  &&|\partial_t y| \leq t^{-2} \ az + t^{-1-\beta_1 + \alpha} \ a^2 
\int_{t}^{t_1} dt' \
t'^{-1-\alpha} \ y(t') \\
&&|\partial_t z| \leq t^{-1- \eta} \ bz + t^{-1+\beta_2 + \alpha} \ a 
\int_{t}^{t_1} dt' \
t'^{-1-\alpha} \ y(t')
\label{4.36e}
\eea
\noi where $\alpha = [k \wedge 3/2] - 1 > 0$, and the goal is to 
prove that (\ref{4.35e})
(\ref{4.36e}) with suitable initial conditions imply $y = z = 0$. \par

\noi \underbar{Part (2)}. Let $Y = \parallel y ; L^{\infty}([t_0, 
t_1])\parallel$. Then
\bea
\label{4.37e}
t^{\alpha} \int_t^{t_1} dt'\ t'^{-1 - \alpha}\ y(t') &\leq& Y\ 
\alpha^{-1}\left (1 -
(t/t_1)^{\alpha} \right ) \nn \\
&\leq& Y \ \alpha^{-1} \left ( 1 - t_0/t_1)^{\alpha} \right ) \equiv 
\bar{Y} \ .
\eea
\noi Substituting (\ref{4.37e}) into (\ref{4.36e}) and integrating 
with $z(t_0) = 0$ yields
\beq
\label{4.38e}
z \leq \exp \left ( b\ \eta^{-1}\ t_0^{-\eta } \right ) \ a\ \bar{Y}\ 
\beta_2^{-1} \ t^{\beta_2}
\ . \eeq
Substituting (\ref{4.37e}) (\ref{4.38e}) into (\ref{4.35e}), 
integrating with $y(t_0) = 0$ and
taking the Supremum over $t$ in $[t_0, t_1]$ yields
$$Y \leq \left \{ \exp \left ( b\ \eta^{-1}\ t_0^{-\eta } \right ) (1 
- \beta_2)^{-1}\ t_0^{-(1
- \beta_2)} + \beta_1^{-1} \ t_0^{-\beta_1} \right \} a^2 \bar{Y}$$
\noi which implies $Y = 0$ and therefore $y = z = 0$ provided
\beq
\label{4.39e}
a^2 \left \{ \exp \left ( b\ \eta^{-1}\ t_0^{-\eta } \right ) (1 - 
\beta_2)^{-1}\ t_0^{-(1
- \beta_2)} + \beta_1^{-1} \ t_0^{-\beta_1} \right \} \alpha^{-1} \left ( 1 -
(t_0/t_1)^{\alpha} \right ) < 1 \ ,  \eeq
\noi a condition which follows from (\ref{4.32e}) for suitable $c(a, b)$.\par

\noi \underbar{Part (3)}. We now take $t_1 = \infty$. The term $bz$ 
in (\ref{4.36e}) can be
exponentiated as in the proof of Part (2). Since in addition the 
statement does not involve
conditions on $a$ and $b$, we can and shall assume without loss of 
generality that $b = 0$ and
$a = 1$. Let
 
\beq
\label{4.40e}
\varepsilon (t) = \mathrel{\mathop {\rm Sup}_{t' \geq t}} 
t'^{\beta_2} y(t') \ .
\eeq
\noi Then $\varepsilon (t)$ is nonincreasing in $t$ and tends to zero 
as $t \to \infty$.
Furthermore for any $t_0 \in I$
\beq
\label{4.41e}
\int_{t_0}^{\infty} dt' \ t'^{-1 - \alpha} \ y(t') \leq \varepsilon 
(t_0) (\alpha +
\beta_2)^{-1} \ t_0^{-(\alpha + \beta_2)} \ . \eeq
\noi Let now $t_0 \in I$ ($t_0$ will eventually tend to $\infty$), 
$y_0 = y(t_0)$ and $z_0 =
z(t_0)$. We estimate $y$ and $z$ for $t \leq t_0$ by integrating 
(\ref{4.35e}) (\ref{4.36e})
(with $t_1 = \infty$, $a = 1$ and $b = 0$) between $t$ and $t_0$. 
Integrating (\ref{4.36e})
yields
\bea
\label{4.42e}
z(t) &\leq& z_0 + (\alpha + \beta_2)^{-2} \varepsilon (t_0) + 
\int_t^{t_0} dt'' \ t''^{-1 +
\beta_2 + \alpha} \int_{t''}^{t_0} dt' \ t'^{-1-\alpha}\ y(t') \nn \\
&\leq& \cdots + \int_t^{t_0} dt' \ t'^{-1 - \alpha} \ y(t') 
\int_t^{t'} dt' \ t''^{-1+
\beta_2 +\alpha}\nn \\
&\leq& z_0 + (\alpha + \beta_2)^{-2} \ \varepsilon (t_0) + (\alpha + 
\beta_2)^{-1}\ Y(t)
\eea
\noi where we have used (\ref{4.41e}) and where
\beq
\label{4.43e}
Y(t) = \int_t^{t_0} dt'\ t'^{-1 + \beta_2}\ y(t') \ .
\eeq
\noi Substituting (\ref{4.42e}) into (\ref{4.35e}), integrating and 
using the fact that $Y(t)$
is decreasing in $t$, we obtain
\beq
\label{4.44e}
y(t) \leq y_0 + t^{-1}\left ( z_0 + (\alpha + \beta_2 )^{-2} 
\varepsilon (t_0)\right ) + (\alpha
+ \beta_2)^{-1}\ t^{-1}Y(t) + y_1(t) \eeq
\noi where
$$y_1(t) = \int_t^{t_0} dt'' \ t''^{-1-\beta_1 + \alpha} 
\int_{t''}^{\infty} dt' \ t'^{-1 -
\alpha} \ y(t') \ .$$
\noi Substituting (\ref{4.44e}) into (\ref{4.43e}) yields
\bea
\label{4.45e}
Y(t) &\leq& y_0 \ \beta_2^{-1} \ t_0^{\beta_2} + \left ( z_0 + 
(\alpha + \beta_2)^{-2}
\varepsilon (t_0) \right ) (1 - \beta_2)^{-1} \ t^{-(1 - \beta_2)} \nn \\
&+& (\alpha + \beta_2)^{-1} \int_t^{t_0} dt'\ t'^{-2 + \beta_2} \ 
Y(t') + Y_1(t)
\eea
\noi where
\bea
&&Y_1(t) = \int_t^{t_0} dt''' \ t'''^{-1 + \beta_2} \ y_1(t''') \nn \\
&&= \int_t^{t_0} dt'' \ t''^{-1 - \beta_1 + \alpha} \int_t^{t''} 
dt'''\ t'''^{-1 + \beta_2}
\int_{t''}^{\infty} dt' \ t'^{- 1 - \alpha}\ y(t')  \nn \\
&&\leq \beta_2^{-1} \int_t^{t_0} dt''\ t''^{-1 - \beta_1 + \alpha + \beta_2}
\int_{t''}^{\infty} dt'\ t'^{-1-\alpha}\ y(t') \nn \\
&&\leq  \beta_2^{-1}\ t^{-\beta_1} (\alpha + \beta_2)^{-2} \ 
\varepsilon (t_0) + \beta_2^{-1}
\int_t^{t_0} dt''\ t''^{-1 - \beta_1} \int_{t''}^{t_0} dt'\
t'^{-1+ \beta_2}\ y(t') \nn \\
&&\leq  \beta_2^{-1}\ t^{-\beta_1} (\alpha + \beta_2)^{-2} \ 
\varepsilon (t_0) + \beta_2^{-1}
\int_t^{t_0} dt''\ t''^{-1 - \beta_1}\ Y(t'') \ .
\label{4.46e}
\eea
\noi Substituting (\ref{4.46e}) into (\ref{4.45e}) yields the 
following inequality for
$Y(t)$~:
\beq
\label{4.47e}
Y(t) \leq f(t) + \int_t^{t_0} dt'\ g(t') \ Y(t')
\eeq
\noi where
$$f(t) = \beta_2^{-1} \varepsilon (t_0) + z_0 (1 - \beta_2)^{-1} \ 
t^{-(1-\beta_2)} + (\alpha +
\beta_2)^{-2} \varepsilon (t_0) \left ( (1 - 
\beta_2)^{-1}t^{-(1-\beta_2)} + \beta_2^{-1}
t^{-\beta_1} \right )$$
$$g(t) = (\alpha + \beta_2)^{-1}\ t^{-2+\beta_2} + \beta_2^{-1} \ 
t^{-1-\beta_1} \ .$$
\noi Note that $f$ and $g$ are decreasing in $t$ and that $g$ is 
integrable at infinity. Let
\begin{eqnarray*}
&&G(t) = \int_t^{\infty} g(t') \ dt' \ , \\
&&\bar{Y}(t) = \int_t^{t_0} dt'\ g(t') \ Y(t')\ .
\end{eqnarray*}
\noi Then (\ref{4.47e}) can be rewritten as
$$\partial_t \bar{Y} = - g Y \geq - gf - g\bar{Y}$$
\noi which is readily integrated with $\bar{Y}(t_0) = 0$ to yield
$$\bar{Y}(t) \leq \int_t^{t_0} dt'\ g(t') \ f(t') \exp \left \{ 
\int_t^{t'} dt''\ g(t'') \right
\} \leq f(t)\ G(t) \exp (G(t))$$
\noi and therefore by (\ref{4.47e}) again
$$Y(t) \leq f(t) \left \{ 1 + G(t) \exp (G(t))\right \} \ .$$
\noi Now $G$ is independent of $t_0$ while $f$ tends to zero when 
$t_0 \to \infty$ for fixed
$t$ under the assumptions made. Letting $t_0 \to \infty$ then shows that
$$Y(t) = \int_t^{t_0} dt' \ t'^{-1 + \beta_2} \ y(t') \longrightarrow 
0 \quad \hbox{when}\  t_0
\to \infty$$
\noi which implies that $y = 0$, from which it follows easily that $z 
= 0$, and therefore
$(w_1, s_1) = (w_2, s_2)$.\par \nobreak\hfill
$\sq$\par

\noi {\bf Remark 4.1.} The necessity of some condition of the type 
(\ref{4.32e}) in Part (2) is
easily understood on the simpler example \\
\beq
\label{4.48e}
|\partial_t y | \leq a^2 \int_t^{t_1} dt'\ y(t')
\eeq
\noi with $t_0 = 0$ and $y(0) = 0$. Defining $Y = \parallel y ; 
L^{\infty}([0, t_1])\parallel$
we obtain
$$|\partial_t y| \leq a^2 (t_1 - t) Y$$
\noi and therefore by integration
$$Y \leq a^2 Y \mathrel{\mathop {\rm Sup}_{0 \leq t \leq 
t_1}}\int_0^t dt'(t_1 - t') = (1/2) a^2\
t_1^2\ Y$$
\noi which implies $Y = 0$ if $a t_1 < \sqrt{2}$. However if $at_1 = 
\pi /2$, (\ref{4.48e})
admits the nonvanishing solution $y = \sin at$. \\

We next prove another property which follows easily from estimates 
similar to those of Lemma
4.2, namely the fact that for suitably bounded solutions $(w, s)$ of 
the auxiliary system,
$w(t)$ tends to a limit $w_+$ when $t \to \infty$.\\

\noi {\bf Proposition 4.3.} {\it Let $k >1$, $\ell_0 \leq [3/2-k]_+$ 
and $\beta > 0$. Let $T
\geq 1$, $t_1 = \infty$ and $I = [T,\infty )$. Let $B_0$ satisfy the 
estimate (\ref{3.17e})
for $m = 0$. Let $(w, s) \in {\cal C} (I, X^{k,\ell_0})$ with $(w, 
t^{\eta - 1} s) \in
L^{\infty}(I, X^{k, \ell_0})$ for some $\eta > 0$ and let $(w, s)$ 
satisfy the first equation
of the system (\ref{2.20e}). Then there exists $w_+ \in H^k$ such 
that $w(t)$ tends to $w_+$
weakly in $H^k$ and strongly in $H^{k'}$ for $0 \leq k' < k$ when $t 
\to \infty$. Furthermore
the following estimate holds for all $t \in I$~:
\beq
\label{4.49e}
\parallel \widetilde{w}(t) - w_+ \parallel_2 \ \leq C\ t^{-\alpha_1}
\eeq
\noi where $\widetilde{w}(t) = U(1/t)w$, and
$$\alpha_1 = \eta \wedge [1/2 \wedge k/3] \wedge (\beta_1 + \beta k)$$
\noi with $\beta_1$ defined by (\ref{4.3e}).}\\

\noi {\bf Proof.} Let $t_0 \in I$, $\widetilde{w}_0 = \widetilde{w}(t_0)$ and
$$a = \parallel w; L^{\infty}(I, H^k)\parallel \quad , \quad b = 
\parallel t^{\eta -1} \nabla s;
L^{\infty}(I,H^{\ell_0})\parallel\ .$$
\noi The first equation of the system (\ref{2.20e}) can be rewritten as
$$\partial_t \left ( \widetilde{w}(t) - \widetilde{w}_0 \right ) = 
t^{-2}\ U(1/t) \ Q(s,w) +
it^{-1}\ U(1/t) \left ( B_0 + B_S(w,w)\right )w$$
\noi where we have omitted the superscript $\infty$ in $B_S$, and therefore
\beq
\label{4.50e}
\partial_t \parallel \widetilde{w}(t) - \widetilde{w}_0 \parallel_2\ 
\leq  t^{-2}\parallel
Q(s,w)\parallel_2 + t^{-1} \parallel B_0w\parallel_2 + t^{-1} 
\parallel B_S(w,w)w\parallel_2 \ .
\eeq
\noi By exactly the same estimate as in (\ref{4.23e}), we obtain
\beq
\label{4.51e}
\parallel Q(s,w)\parallel_2 \ \leq C |s|_{\ell_0}^{\dot{}}\ |w|_k 
\leq C\ ab\ t^{1-\eta}\ .
\eeq
\noi We next estimate
\beq
\label{4.52e}
\parallel B_0w\parallel_2 \ \leq \parallel B_0\parallel_{3/\delta} \ 
\parallel w\parallel_r \
\leq \ C\ ab_0\ t^{-[1/2 \wedge k/3]}  \eeq
\noi with $\delta = \delta (r) = [k \wedge 3/2]$, by (\ref{3.17e}) 
and Sobolev inequalities,
and
\bea
\label{4.53e}
\parallel B_S(w, w)w\parallel_2 &\leq& \parallel B_S(w, 
w)\parallel_{3/\delta} \ \parallel
w\parallel_r \nn \\
&\leq& C \parallel \omega^{[3/2 - k]_+}\ B_S(w,w)\parallel_{2} \ |w|_k
\eea
\noi with the same $\delta$. The last norm of $B_S$ in (\ref{4.53e}) 
is estimated exactly as in
the proof of Lemma 4.1 (see (\ref{4.10e}) (\ref{4.11e}) (\ref{4.12e})) as
\beq
\label{4.54e}
\parallel \omega^{[3/2 - k]_+}\ B_S(w,w)\parallel_{2} \ \leq C\ 
t^{-\beta_1 - \beta k}\
I_{m_1}^{\infty} (|w|_k^2) \leq C\ a^2 \ t^{-\beta_1 - \beta k} \ . \eeq
\noi Substituting (\ref{4.51e}) (\ref{4.52e}) (\ref{4.53e}) 
(\ref{4.54e}) into (\ref{4.50e})
and integrating between $t_0$ and $t$ yields
$$\parallel \widetilde{w}(t) - \widetilde{w}(t_0) \parallel_2\ \leq 
C\left \{ (t \wedge
t_0)^{-\eta} ab + (t \wedge t_0)^{-[1/2 \wedge k/3]} ab_0 + (t \wedge 
t_0)^{-\beta_1 - \beta
k} a^3 \right \}$$
\noi from which it follows that $\widetilde{w}(t)$ and therefore also $w(t)$
has a strong limit $w_+$ in $L^2$ when $t\to \infty$ and that 
(\ref{4.49e}) holds. Since in
addition $w(t)$ is bounded in $H^k$, it follows by a standard 
compactness argument that $w_+ \in
H^k$ and that $w(t)$ tends to $w_+$ in the other topologies stated in 
the Proposition. \par
\hfill $\sq$\par

We now turn to the problem of existence of solutions of the auxiliary 
system (\ref{2.20e}),
with the aim of proving that that system defines an asymptotic 
dynamics for large times and
preferably up to infinity in time. Here however, we encounter the 
difficulties described in
Section 2 and arising from the different propagation properties of 
the Schr\"odinger and wave
equations. First of all for $t_1 = t_0$, even if $B_0 = 0$, the 
estimates of Lemma 4.1 are
insufficient to prevent blow up of the solutions in a finite time 
after $t_0$, independently of
the size of $t_0$ and of the initial data for $(w, s)$ at $t_0$. In 
fact, if in the estimates
(\ref{4.1e}=) (\ref{4.2e}=) we set $b_0 = 0$, omit the second 
inequality and take $s = 0$ in the
first one, we obtain the following stronger estimate for $y = |w'|_k = |w|_k$
\beq
\label{4.55e}
\partial_t\ y \leq C\  t^{-1-\beta_1}\ y \int_{t_0/t}^1 d\nu \ \nu^{-1-m}\ y(\nu t)^p
\eeq
\noi where $m = m_1 + 1/2 > \beta_1$ for $\beta \l 
eq 1$, and $p = 2$. We shall prove that
(\ref{4.55e}) does not prevent finite time blow up by showing that 
equality in (\ref{4.55e})
implies such a blow up. Taking $y^p$ instead of $y$ as the unknown 
function and rescaling $t$
and $y$, we can take $p = 1$, $t_0 = 1$ and $C = 1$ without loss of 
generality. We are
therefore led to consider the equation
\beq
\label{4.56e}
\partial_t\ y = t^{-1-\beta_1}\ y \int_{1/t}^1 d\nu \ \nu^{-1-m}\ y(\nu t)
\eeq
\noi or equivalently
\beq
\label{4.57e}
\partial_t\ y = t^{-1-\beta_1+ m }\ y \int_{1}^t dt' \  t'^{-1-m}\ y(t') \ .
\eeq

\noi {\bf Warning 4.1.} {\it Let $0 < \beta_1 < m$. Then the solution 
of the equation
(\ref{4.57e}) with initial data $y(1) = y_0 > 0$ blows up in a finite time.} \\

The proof will be given in Appendix A. \par

The previous result encourages us to take $t_1 > t_0$ and actually 
the situation improves in
that case and in particular we shall prove the existence of solutions 
defined in $[t_0, t_1]$ if
$B_0 = 0$ in Proposition 4.4 below. Of course for $t_1 < \infty$, by 
the previous argument, we
shall be unable to exclude finite time blow up after $t_1$. On the 
other hand, if $B_0 \not=
0$, we cannot exclude finite time blow up between $t_0$ and $t_1$ if 
$t_1$ is sufficiently
large. Actually, we shall show that equality in a stronger version of 
(\ref{4.1e}) implies
such a blow up. We again drop the inequality (\ref{4.2e}) and take $s = 0$ in
(\ref{4.1e}). Omitting in addition the first term with $b_0$, we are left with
\beq
\label{4.58e}
\partial_t\ y = C\left\{ y^{1- 1/k} + t^{-1-\beta_1}\ y 
\int_{1}^{t_1/t} \nu^{-1-m} \  y(\nu
t)^2 \right \} .
\eeq
\noi Since the solution of (\ref{4.58e}) is increasing in time for $t 
\geq t_0$, blow up for
(\ref{4.58e}) is implied by blow up for the equation
\beq
\label{4.59e}
\partial_t\ y = C\left\{ y^{1- 1/k} + t^{-1-\beta_1}\ y^3 m^{-1} 
\left ( 1 - (t/t_1)^m \right )
\right \} \ .
\eeq
\noi Now if blow up occurs for $t \leq T^*$ for the equation
\beq
\label{4.60e}
\partial_t\ y = C\left ( y^{1- 1/k} + t^{-1-\beta_1}\ y^3 m^{-1} (1 - 
2^{-m}) \right )
\eeq
\noi then a fortiori blow up will occur for the equation 
(\ref{4.59e}) if $t_1 \geq 2T^*$. It
is therefore sufficient to prove blow up for (\ref{4.60e}), which 
after rescaling can be
rewritten as
\beq
\label{4.61e}
\partial_t\ y = k\ y^{1- 1/k} + t^{-1-\beta_1}\ y^3\ .
\eeq

\noi {\bf Warning 4.2.} {\it Let $2k > \beta_1$. Then the solution of 
the equation
(\ref{4.61e}) with initial data $y(t_0) > 0$ at time $t_0 \geq 1$ 
blows up in a finite time.}
\\

The proof will be given in Appendix A. The condition $2k > \beta_1$ 
is always satisfied in the
present situation. \par

We now prove the main result of this section, namely the existence of 
solutions of the
auxiliary system (\ref{2.20e}) defined up to $t_1$, possibly with 
$t_1 = \infty$, for $B_0 = 0$
and for initial data given at sufficiently large $t_0 < t_1$. \\

\noi {\bf Proposition 4.4.} {\it Let $B_0 = 0$. Let $1 < k \leq 
\ell$, $\ell >3/2$
and $0 < \beta < 1$. Let $\beta_2 < 1$, where $\beta_2$ is defined by 
(\ref{4.5e}). Let
$(w_0, \widetilde{s}_0) \in X^{k,\ell}$ and let $y_0 = |w_0|_k$ and 
$\widetilde{z}_0 =
|\widetilde{s}_0|_{\ell}^{\dot{}}$. Then there exists $T_0 < \infty$ 
depending on $(y_0,
\widetilde{z}_0)$ such that for all $t_0 \geq T_0$, there exists $T < 
t_0$, depending on
$(y_0, \widetilde{z}_0)$ and on $t_0$, such that for all $t_1$, $t_0 
\leq t_1 \leq \infty$,
the auxiliary system (\ref{2.20e}) with initial data $(w, s)(t_0) = 
(w_0, t_0^{\beta_2}
\widetilde{s}_0)$ has a unique solution $(w,s)$ in the interval $I = 
[T, t_1)$ such that
$(w, t^{-\beta_2}s) \in ({\cal C} \cap L^{\infty})(I, X^{k,\ell})$. 
One can take
\bea
\label{4.62e}
&&T_0 = C \left \{ \left ( \widetilde{z}_0 + y_0^2 \right )^{1/(1 - 
\beta_2)} \vee
y_0^{2/\beta_1} \right \} \ , \\
&&T = t_0^{\beta_2} \ T_0^{1-\beta_2} \ ,
\label{4.63e}
\eea
\noi and the solution $(w, s)$ is estimated for all $t \in I$ by
\bea
\label{4.64e}
&&|w(t)|_k \leq 2y_0 \ , \\
&&|s(t)|_{\ell}^{\dot{}} \leq \left ( 2 \widetilde{z}_0 + C\ y_0^2 
\right ) (t_0
\vee t)^{\beta_2} \ .  \label{4.65e}
\eea
}\vskip 3 truemm

\noi {\bf Proof.} The proof consists in exploiting the estimates of 
Lemmas 4.1 and 4.2 in order
to show that the map $\Gamma : (w, s) \to (w', s')$, where $(w', s')$ 
is defined from $(w, s)$
by Proposition 4.1, is a contraction of a suitable subset of ${\cal 
C}(I, X^{k, \ell})$ for a
suitably time rescaled norm of $L^{\infty}(I, X^{0,\ell_0})$. We 
first consider the interval
$I = [t_0, t_1)$ and we define the set
$${\cal R} = \left \{ (w, s) \in {\cal C}(I, X^{k,\ell}) : \ 
\parallel w; L^{\infty}(I,
H^k)\parallel \ \leq Y, \parallel t^{-\beta_2} \nabla 
s;L^{\infty}(I,H^{\ell})\parallel \ \leq Z
\right \} \  ,$$
\noi for $Y > 0$, $Z > 0$. Let $(w, s) \in {\cal R}$ and $(w', s') = \Gamma (w,
s)$. Let $y = |w(t)|_k$, $y' = |w'(t)|_k$, $z = 
|s(t)|_{\ell}^{\dot{}}$ and $z' =
|s'(t)|_{\ell}^{\dot{}}$. From Lemma 4.1, namely (\ref{4.1e}) 
(\ref{4.2e}) with $b_0 = 0$ and
with an overall constant omitted, we obtain
\beq
\label{4.66e}
\left \{ \begin{array}{l} \partial_t \ y' \leq t^{-2+ \beta_2} \ Zy' +
t^{-1-\beta_1} \ Y^2 y' \\
\\
\partial_t \ z' \leq t^{-2+ \beta_2} \ Zz' +
t^{-1+\beta_2} \ Y^2 \ .
\end{array}\right .
\eeq
\noi Integrating from $t_0$ to $t$ with $(y', z')(t_0) = (y,z)(t_0) = 
(y_0, z_0)$ where $z_0 =
\widetilde{z}_0 t_0^{\beta_2}$, we estimate
\beq
\label{4.67e}
Y' = \parallel y' ; L^{\infty}(I)\parallel \quad , \quad Z' = 
\parallel z'; L^{\infty}(I)
\parallel
  \eeq
\noi by
\beq
\label{4.68e}
\left \{ \begin{array}{l} Y' \leq y_0 \exp \left \{ (1 - 
\beta_2)^{-1} \ t_0^{-1 + \beta_2}\ Z
+ \beta_1^{-1} \ t_0^{-\beta_1}\ Y^2 \right \} \\ \\
Z' \leq \left ( \widetilde{z}_0 + \beta_2^{-1}\ Y^2 \right ) \exp 
\left \{ (1 - \beta_2)^{-1}\
t_0^{-1 + \beta_2} \ Z \right \} \ . \end{array}\right .
\eeq
\noi We now impose
\beq
\label{4.69e}
\left \{ \begin{array}{l} (1 - \beta_2) \ t_0^{1-\beta_2} \geq 2(\ell 
n 2)^{-1} Z \\
\\
\beta_1 \ t_0^{\beta_1} \geq 2(\ell n 2)^{-1}\ Y^2
\end{array}\right .
\eeq
\noi and choose
\beq
\label{4.70e}
Y = 2y_0 \quad , \quad Z = \sqrt{2} \left ( \widetilde{z}_0 + 4 
\beta_2^{-1} \ y_0^2 \right )
\eeq
\noi thereby ensuring that $Y' \leq Y$, $Z' \leq Z$, so that the set 
${\cal R}$ is mapped into
itself by $\Gamma$. The conditions (\ref{4.69e}) can be rewritten as
\beq
\label{4.71e}
\left \{ \begin{array}{l} (1 - \beta_2) \ t_0^{1-\beta_2} \geq 
2\sqrt{2} (\ell n 2)^{-1} \left
( \widetilde{z}_0 + 4 \beta_2^{-1} \ y_0^2 \right ) \\ \\
\beta_1 \ t_0^{\beta_1} \geq 8(\ell n 2)^{-1}\ y_0^2
\end{array}\right .
\eeq
\noi and hold for all $t_0 \geq T_0$ for $T_0$ satisfying 
(\ref{4.62e}) with suitable $C$. \par

We next show that the map $\Gamma$ is a contraction on ${\cal R}$. We 
use the notation of Lemma
4.2 and in addition
\beq
\label{4.72e}
\left \{ \begin{array}{l} y_- = \parallel w_-(t)\parallel_2 \quad , \quad z_- =
|s_-(t)|_{\ell_0}^{\dot{}} \ , \\ \\
Y_- = \parallel y_- ; L^{\infty}(I) \parallel \quad , \quad Z_- = 
\parallel t^{-\beta_2}\ z_- ;
L^{\infty}(I) \parallel
\end{array}\right .
\eeq
\noi and a similar notation for primed quantities. From Lemma 4.2, in 
particular (\ref{4.18e})
(\ref{4.19e}), and again with an overall constant omitted, we obtain
\beq
\label{4.73e}
\left \{ \begin{array}{l} \partial_t \ y'_- \leq t^{-2} \ Y\ z_- +
t^{-1-\beta_1} \ Y^2 Y_- \\
\\
\partial_t \ z'_- \leq t^{-2+ \beta_2} \ Z(z_- + z'_-) +
t^{-1+\beta_2} \ Y\ Y_-
\end{array}\right .
\eeq
\noi and by integration with $(y'_-, z'_-)(t_0) = 0$,
\beq
\label{4.74e}
\hskip - 0.5 truecm \left \{ \begin{array}{l} Y'_- \leq (1 - 
\beta_2)^{-1}\ t_0^{-1+ \beta_2}
YZ_- + \beta_1^{-1} \ t_0^{-\beta_1} \ Y^2 \ Y_- \\
\\
Z'_- \leq \exp \left ( (1 - \beta_2)^{-1}\ t_0^{-1+ \beta_2} Z \right 
) \left \{ (1-
\beta_2)^{-1} \ t_0^{-1 + \beta_2} \ ZZ_-\ + \beta_2^{-1} \ YY_- \right \} \ .
  \end{array}\right .
\eeq
\noi The second inequality in (\ref{4.74e}) reduces to
\beq
\label{4.75e}
Z'_- \leq \sqrt{2} \left \{ (1 - \beta_2)^{-1} \ t_0^{-1 + \beta_2} \ 
ZZ_- + \beta_2^{-1} \ Y
Y_- \right \}
\eeq
\noi under the first condition in (\ref{4.69e}) imposed previously. \par

We now ensure that the map $\Gamma$ is a contraction for the norms 
defined by (\ref{4.72e}) in
the form
\beq
\label{4.76e}
\left \{ \begin{array}{l} Y'_- \leq ( c^{-1}\ Z_- + Y_- )/4 \\
\\
Z'_- \leq (Z_- + c\ Y_-)/4
  \end{array}\right .
\eeq
\noi which implies
\beq
\label{4.77e}
Z'_- + c\ Y'_- \leq (Z_- + c \ Y_-)/2
\eeq
\noi by taking $c = 8 \beta_2^{-1} Y$ and imposing the conditions
$$\left \{ \begin{array}{l} (1 - \beta_2)t_0^{1-\beta_2} \geq 8 Z 
\quad , \quad \beta_1
t_0^{\beta_1} \geq 4Y^2 \\ \\
(1 - \beta_2)\ t_0^{1-\beta_2} \geq 4c\ Y = 32 \beta_2^{-1} \ Y^2
  \end{array}\right .$$
\noi which follow again from (\ref{4.62e}) for all $t_0 \geq T_0$. \par

We have proved that $\Gamma$ maps ${\cal R}$ into itself and is a 
contraction for the norms
(\ref{4.72e}). By a standard compactness argument, ${\cal R}$ is 
easily shown to be closed for
the latter norms. Therefore $\Gamma$ has a unique fixed point in 
${\cal R}$, which completes
the proof for $t \geq t_0$. \par

We now turn to the case $t \leq t_0$, namely we consider the interval 
$I = [T, t_0]$ for some
$T < t_0$. The proof proceeds in exactly the same way, with however 
slightly different norms.
In addition, one has to take into account the following fact~: the 
various integrals
$I_m^{t_1}$ that occur in (\ref{4.1e}) (\ref{4.2e}) and (\ref{4.18e}) 
(\ref{4.19e}) involve $w$
and $w_1$, $w_2$ up to time $t_1$. In the subinterval $[t_0, t_1]$, 
one takes $w = w_1 = w_2 =$
the solution constructed at the previous step (in particular $w_- = 
0$ for $t \geq t_0$, so
that actually no contribution from the interval $[t_0, t_1]$ occurs 
in (\ref{4.18e})
(\ref{4.19e})). In (\ref{4.1e}) (\ref{4.2e}) the contribution of the 
interval $[t_0, t_1]$ is
taken into account by using the fact that all the integrals over time 
in the relevant
$I_m^{t_1}$ are convergent at infinity and that we shall eventually 
use the same ansatz
$|w(t)|_k \leq Y = 2y_0$ both for $t \leq t_0$ and $t \geq t_0$. With 
this in mind, we complete
the proof by simply giving the computational details. We consider the set
$${\cal R}_< = \left \{ (w, s) \in {\cal C}(I, X^{k,\ell}) ; 
\parallel w; L^{\infty}(I,
H^k)\parallel \leq Y, \parallel \nabla s;L^{\infty}(I,H^{\ell})\parallel \leq Z
\right \} \  .$$
\noi For $(w, s) \in {\cal R}_<$, $(w',s') = \Gamma (w, s)$ and $y$, 
$y'$, $z$ and $z'$
defined as previously, we estimate by Lemma 4.1, again with an 
overall constant omitted
\beq
\label{4.78e}
\left \{ \begin{array}{l} |\partial_t \ y'| \leq t^{-2} \ Zy' +
t^{-1-\beta_1} \ Y^2 y' \\
|\partial_t \ z'| \leq t^{-2} \ Zz' + t^{-1+\beta_2} \ Y^2
\end{array}\right .
\eeq
\noi and by integration from $t$ to $t_0$ with $(y', z')(t_0) = 
(y,z)(t_0) = (y_0, z_0)$, we
obtain $y'(t) \leq Y'$, $z'(t) \leq Z'$ with
\beq
\label{4.79e}
\left \{ \begin{array}{l} Y' \leq y_0 \exp \left ( t^{-1}Z + 
\beta_1^{-1} \ t^{-\beta_1} \ Y^2
\right ) \\ \\
Z' \leq \left ( z_0 + \beta_2^{-1} \ Y^2\ t_0^{\beta_2} \right ) \exp 
(t^{-1} Z) \ .
  \end{array}\right .
\eeq
\noi We now impose
\beq
\label{4.80e}
\left \{ \begin{array}{l} t \geq 2(\ell n \ 2)^{-1}\ Z  \\
\\
\beta_1 \ t^{\beta_1}  \geq 2(\ell n \ 2)^{-1}\ Y^2
\end{array}\right .
\eeq
\noi and choose
\beq
\label{4.81e}
Y = 2y_0 \quad , \quad Z = \sqrt{2} \left ( z_0 + 4 \beta_2^{-1} \ 
y_0^2 \ t_0^{\beta_2} \right
) \eeq
\noi thereby ensuring that $Y' \leq Y$, $Z' \leq Z$ so that the set 
${\cal R}_<$ is mapped
into itself by $\Gamma$. The conditions (\ref{4.80e}) can be rewritten as
\beq
\label{4.82e}
\left \{ \begin{array}{l} \beta_1\ t^{\beta_1}  \geq 8(\ell n \ 
2)^{-1}\ y_0^2  \\
\\
t  \geq 2\sqrt{2} (\ell n \ 2)^{-1}\left ( \widetilde{z}_0 + 4 
\beta_2^{-1} \ y_0^2 \right )
t_0^{\beta_2}
  \end{array}\right .
\eeq
\noi and hold for all $t \geq T$ with $T$ defined by (\ref{4.63e}) 
and $T_0$ satisfying
(\ref{4.62e}) for suitable $C$. \par

We next prove that $\Gamma$ is a contraction on ${\cal R}_<$ for the 
norm in $L^{\infty}(I,
X^{0,\ell_0})$. With the notation of Lemma 4.2 and in addition   

\beq
\label{4.83e}
\left \{ \begin{array}{l}  y_- = \ \parallel w_-(t)\parallel_2 \quad 
, \quad z_- =
|s_-(t)|_{\ell_0}^{\dot{}}\\
\\
Y_- = \ \parallel y_- ; L^{\infty}(I)\parallel \quad , \quad Z_- = 
\parallel z_- ; L^{\infty} (I)
\parallel \end{array}\right .
\eeq
 
\noi and a similar notation for primed quantities, we obtain from 
(\ref{4.18e}) (\ref{4.19e})
\beq
\label{4.84e}
\left \{ \begin{array}{l} |\partial_t \ y'_-| \leq t^{-2} \ Yz_- + 
t^{-1-\beta_1} \ Y^2 \ Y_- \\
\\
|\partial_t \ z'_-| \leq t^{-2} \ Z(z_- + z'_-) + t^{-1+\beta_2} \ YY_- \ .
\end{array}\right .
\eeq
\noi By integration between $t$ and $t_0$, we deduce therefrom
\beq
\label{4.85e}
\left \{ \begin{array}{l} Y'_- \leq t^{-1}YZ_- + \beta_1^{-1} \ 
t^{-\beta_1}\ Y^2Y_-  \\
\\
Z'_- \leq \left ( t^{-1} ZZ_- + \beta_2^{-1} \ t_0^{\beta_2} \ YY_- \ 
\right ) \exp
(t^{-1} Z)
  \end{array}\right .
\eeq
\noi thereby ensuring the contraction in the form (\ref{4.76e}) which 
implies (\ref{4.77e}) by
taking $c = 8 \beta_2^{-1} Y t_0^{\beta_2}$ and imposing
$$t \geq 8 Z \quad , \quad \beta_1 t^{\beta_1} \geq 4 Y^2 \quad , 
\quad t \geq 4c \ Y = 32
\beta_2^{-1} \ t_0^{\beta_2}\ Y^2$$
\noi which hold for all $t \geq T$ with the choice (\ref{4.81e}) 
under the conditions
(\ref{4.62e}) (\ref{4.63e}). With the previous estimates available, 
the proof proceeds as in
the case $t \geq t_0$.\par \hfill
$\sq$\par

\mysection{Asymptotics and wave operators for the linear system}
\hspace*{\parindent} In this section we study the asymptotic 
properties of solutions of the
linear equation (\ref{2.22e}) in the form (\ref{2.24e}) at the level 
of regularity of $H^k$ with
$k \geq 1$ for $w$. In particular we solve the Cauchy problem at 
infinity, thereby
constructing the wave operators in $H^k$. For the linear equation 
(\ref{2.24e}), the wave
operators in $L^2$ can be easily constructed by a variant of Cook's 
method and the construction
of the wave operators in $H^k$ reduces to a regularity problem for 
the $L^2$ wave operators
thereby obtained. As a preliminary to that study, we shall first 
solve the Cauchy problem for the
equation (\ref{2.24e}) with finite initial time. We emphasize the 
fact that in this section we
do not strive after any kind of optimality in the treatment of the 
linear equation, since we
are mainly interested in a form of that treatment that can be 
incorporated in that of the fully
interacting system. \\
 
\noi {\bf Proposition 5.1.} {\it Let $I = [1 , \infty )$, let $k \geq 
1$ and let $B_0$ satisfy
the estimates (\ref{3.17e}) for $0 \leq m \leq k$. Let $t_0 \in I$ 
and $w_0 \in H^k$. Then the
equation (\ref{2.24e}) has a unique solution $w \in {\cal C}(I, L^2)$ 
with $w(t_0) = w_0$.
Furthermore $w \in {\cal C}(I, H^k) \cap L^{\infty}(I, L^2)$ and $w$ 
satisfies the conservation
law
$$\parallel w(t)\parallel_2\ = {\rm const.}$$
\noi and the estimate
\beq
\label{5.1e}
|w(t)|_k \leq \left ( 1 + C_k |t- t_0|(t \vee t_0)^{\bar{k}-1} \right ) |w_0|_k
\eeq
\noi where $\bar{k} = k$ for integer $k$ and $\bar{k} = k + 
\varepsilon$ with $\varepsilon > 0$
for noninteger $k$.} \\

\noi {\bf Proof.} It follows easily from standard arguments and from 
Lemma 3.2 that $w$ exists
and satisfies the properties stated except possibly the estimate 
(\ref{5.1e}), and we
concentrate on the proof of the latter, assuming without loss of 
generality that $|w_0|_k =
1$. We first prove (\ref{5.1e}) by induction  for integer $k \geq 1$. 
Let $0 \leq j \leq k -1$
and $y_j = \parallel \omega^jw\parallel_2$. From (\ref{2.24e}) and 
from the Leibnitz formula
and Sobolev inequalities, we obtain
$$\left | \partial_t\ y_{j+1} \right | \leq C\ t^{-1} \Big \{ 
\parallel \nabla B_0
\parallel_{\infty} \ \parallel \omega^j w\parallel_2 + \sum_{|\alpha 
| = j+1} \parallel
\partial^{\alpha} B_0 \parallel_{\infty} \ \parallel w \parallel_2 \Big \}$$
\noi and therefore by (\ref{3.17e})
\beq
\label{5.2e}
\left | \partial_t\ y_{j+1} \right | \leq C \ b_0 \left ( y_j + t^j 
\right ) \ .
\eeq
\noi Substituting the induction assumption for $y_j$ and integrating 
(\ref{5.2e}) between $t_0$
and $t$, we obtain
\begin{eqnarray*}
y_{j+1} &\leq& 1 + C \ b_0 \left ( 1 + (C_j + 1) (t \vee t_0)^j 
\right ) |t-t_0| \\
&\leq& 1 + C_{j+1} |t-t_0| (t \vee t_0)^j
\end{eqnarray*}
\noi with $C_{j+1} = Cb_0(C_j + 2)$. This completes the proof for 
integer $k$. \par

Let now $k = k_0 + \theta$ with integer $k_0 \geq 1$ and $0 < \theta 
< 1$. We estimate
\bea
\label{5.3e}
&&\left | \partial_t \parallel \omega^kw\parallel_2 \right | \leq 
t^{-1} \parallel [\omega^k ,
B_0 ] w \parallel_2 \nn \\
&&\leq C t^{-1} \left \{ \parallel \nabla B_0 \parallel_{\infty} \ 
\parallel \omega^{k-1}
w\parallel_2 + \parallel \omega^k B_0 \parallel_{3/\delta} \ 
\parallel w \parallel_r \right \}
\eea
\noi by Lemma 3.2, with $0 < \delta = \delta (r) \leq 1$,
$$\cdots \leq C \ b_0 \left ( \parallel \omega^{k-1} w\parallel_2 + 
t^{k-1-\delta /3} \parallel
\omega^{\delta} w\parallel_2 \right )$$
\noi by (\ref{3.17e}) and Sobolev inequalities. We next interpolate and obtain
$$\cdots \leq C \ b_0 \left ( \parallel \omega^{k_0-1} 
w\parallel_2^{1 - \theta} \ \parallel
\omega^{k_0} w\parallel_2^{\theta}+ t^{k-1-\delta /3} \parallel 
w\parallel_2^{1 - \delta} \
\parallel \nabla w \parallel_2^{\delta} \right ) \ .$$
\noi We finally substitute the estimate (\ref{5.1e}) for the integer 
values $k_0 - 1$, $k_0$,
and 1 and integrate between $t_0$ and $t$, thereby obtaining
$$\parallel \omega^kw \parallel_2 \ \leq 1 + C\ b_0 |t-t_0| (t \vee 
t_0)^{k-1 + 2 \delta /3}$$
\noi which yields (\ref{5.1e}) with $\varepsilon = 2 \delta /3$. \par
\hfill
$\sq$\par

The fact that a direct $H^k$ estimate of the solution does not 
prevent its $H^k$ norm to
increase as a power of $t$ is a warning of the fact that the 
construction of the wave operators
at that level of regularity is not trivial. The same fact appeared 
already in Section 4 above
in Warning 4.2 and compelled us to assume $B_0 = 0$ in Proposition 4.4. \par

We next construct the $L^2$-wave operators for (\ref{2.24e}).\\

\noi {\bf Proposition 5.2.} {\it Let $I = [1 , \infty )$ and let $B_0$ satisfy
the estimates (\ref{3.17e}) for $m = 0$. \par

(1) Let $W \in {\cal C}(I,L^2)$ with $U(1/t)W \in {\cal C}^1(I, 
L^2)$, satisfying
\beq
\label{5.4e}
\parallel R(W) \parallel_2 \ \leq \ c_0 \ t^{-1 - \lambda_0}
\eeq
\noi for some $\lambda_0 > 0$ and for all $t \in I$. Then there 
exists a unique solution $w \in
{\cal C}(I, L^2)$ of the equation (\ref{2.24e}), such that $w(t) - 
W(t)$ tends to zero
strongly in $L^2$ when $t \to \infty$. Furthermore, for all $t \in I$,
\beq
\label{5.5e}
\parallel w(t) - W(t) \parallel_2 \ \leq \ c_0 \ \lambda_0^{-1} \ 
t^{- \lambda_0} \ .
\eeq
\indent The solution $w$ is the norm limit in $L^{\infty}(I,L^2)$ as 
$t_0 \to \infty$ of the
solution $w_{t_0}$ of the equation (\ref{2.24e}) with initial 
condition $w_{t_0}(t_0) = W(t_0)$
obtained in Proposition 5.1, and the following estimate holds for all 
$t \in I$~:
\beq
\label{5.6e}
\parallel w_{t_0}(t)  - w(t) \parallel_2 \ \leq \ c_0 \ 
\lambda_0^{-1} \ t_0^{- \lambda_0} \ .
\eeq
\indent (2) Let in addition $W \in L^{\infty}(I, H^k)$ for some $k$, 
$0 < k < 3/2$. Then there
exists $w_+ \in H^k$ such that $W(t)$ tends to $w_+$ strongly in 
$L^2$ and weakly in $H^k$ when
$t \to \infty$, and the following estimate holds for all $t \in I$~:
\beq
\label{5.7e}
\parallel W(t) - w_+ \parallel_2 \ \leq \ C \left ( t^{- \lambda_0} + 
t^{- k/3} \right ) \ .
\eeq
\noi Conversely let $w_+ \in H^k$ for some $k$, $0 < k < 3/2$, and 
let $W_1 = U^*(1/t)w_+$.
Then $W_1$ satisfies the assumptions of Part (1) with $\lambda_0 = k/3$. \par

Let $W$, $w_+$ and $W_1$ be related as above. Then the solutions of 
the equation (\ref{2.24e})
constructed in Part (1) from $W$ and $W_1$ coincide. \par

(3) Let $w_+ \in L^2$. Then the equation (\ref{2.24e}) has a unique 
solution $w \in {\cal C}(I,
L^2)$ such that $w(t)$ tends to $w_+$ strongly in $L^2$ when $t \to 
\infty$.} \\

\noi {\bf Proof.} \underbar{Part (1)}. Following the sketch of 
Section 2, we look for $w$ in
the form $w = W + q$, so that $q$ satisfies the equation
\beq
\label{5.8e}
\partial_t \ q = i(2t^2)^{-1} \Delta q + i \ t^{-1} \ B_0 \ q - R(W)
\eeq
\noi and therefore the a priori estimate
\beq
\label{5.9e}
|\partial_t \parallel q \parallel_2 | \leq \ \parallel R(W) 
\parallel_2 \ \leq c_0 \
\lambda_0^{-1} \  t^{-1 - \lambda_0} \ .
  \eeq
\noi Define $w_{t_0}$ as in Part (1) and let $w_{t_0} = W + q_{t_0}$ 
so that $q_{t_0}(t_0) = 0$.
Integrating (\ref{5.9e}) between $t_0$ and $t$ yields
\beq
\label{5.10e}
\parallel q_{t_0}(t) \parallel_2  \ \leq c_0 \ \lambda_0^{-1}
| t^{- \lambda_0} - t_0^{- \lambda_0}|
  \eeq
\noi and therefore, by $L^2$ norm conservation for the difference of 
two solutions
\beq
\label{5.11e}
\parallel q_{t_0}(t) - q_{t_1} (t)\parallel_2 | \ = \ \parallel 
q_{t_0}(t_1) \parallel_2
\ \leq c_0 \ \lambda_0^{-1} | t_1^{- \lambda_0} - t_0^{- \lambda_0}|
  \eeq
\noi for any $t_0$ and $t_1$, $1 \leq t_0$, $t_1 < \infty$. This 
proves convergence of $q_{t_0}$
and therefore of $w_{t_0}$ in norm in $L^{\infty}(I, L^2)$. Let $w$ 
be the limit of $w_{t_0}$.
Taking the limit $t_0 \to \infty$ in (\ref{5.10e}) yields 
(\ref{5.5e}), while taking the
limit $t_1 \to \infty$ in (\ref{5.11e}) yields (\ref{5.6e}). Clearly 
$w$ satisfies the
equation (\ref{2.24e}). \par

\noi \underbar{Part (2)}. $W$ satisfies the equation
\beq
\label{5.12e}
\partial_t\ U(1/t)\ W = i\ t^{-1}\ U(1/t) \ B_0\ W + U(1/t)\ R(W) \ .
\eeq
\noi From (\ref{3.17e}) and Sobolev inequalities, we obtain
\beq
\label{5.13e}
\parallel B_0 W \parallel_2 \ \leq \ \parallel B_0\parallel_{3/k} \ 
\parallel W\parallel_r \
\leq C\ ab_0\ t^{-k/3} \eeq
\noi where $a = \parallel W;L^{\infty}(I, H^k)\parallel$ and $k = 
\delta (r)$. Integrating
(\ref{5.12e}) between $t_1$ and $t_2$ and using (\ref{5.4e}) and 
(\ref{5.13e}), we obtain
\beq
\label{5.14e}
\parallel U(1/t_1) W(t_1) - U(1/t_2) W(t_2)\parallel_2 \ \leq C \left 
\{ |t_1^{-k/3} -
t_2^{-k/3} |+|t_1^{-\lambda_0} - t_2^{-\lambda_0} | \right \} \eeq
\noi for any $t_1$ and $t_2$, $1 \leq t_1$, $t_2 < \infty$. Therefore 
$U(1/t) W(t)$ and therefore
also $W(t)$ has a strong limit $w_+$ in $L^2$, and
\beq
\label{5.15e}
\parallel U(1/t) \ W(t) - w_+ \parallel_2 \ \leq C \left ( t^{-k/3} + 
t^{-\lambda_0} \right ) \
, \eeq
\noi from which (\ref{5.7e}) follows. Furthermore by a standard 
compactness argument, $w_+ \in
H^k$ with $|w_+|_k \leq a$ and $W(t)$ tends to $w_+$ weakly in $H^k$. \par

Let now $w_+ \in H^k$ and $W_1 = U^*(1/t) w_+$. Then
\beq
\label{5.16e}
R(W_1) = - i t^{-1} \ B_0 \ U^*(1/t) \ w_+
\eeq
\noi so that
\beq
\label{5.17e}
\parallel R(W_1) \parallel_2 \ \leq C\ b_0|w_+|_k \ t^{-k/3}
\eeq
\noi by (\ref{5.13e}). The last statement follows from the fact that 
$W$ and $W_1$ have the
same limit $w_+$ in $L^2$ and from $L^2$ norm conservation for the 
equation (\ref{2.24e}). \par

\noi \underbar{Part (3)} follows from Parts (1) and (2) by a standard 
density argument.
\par\nobreak
\hfill $\sq$\par

We next prove that the solutions with asymptotic properties in $L^2$ 
obtained in Proposition
5.2 exhibit similar asymptotic properties in $H^k$ under suitable 
additional assumptions. \\

\noi {\bf Proposition 5.3.} {\it Let $I = [1 , \infty )$, let $k \geq 
1$ and let $B_0$ satisfy
the estimates (\ref{3.17e}) for $0 \leq m \leq k$. Let $\lambda > 0$ 
and $\lambda_0 >
\lambda + k$ and let $U(1/t) W \in {\cal C}^1(I,H^k)$ satisfy the 
estimates (\ref{5.4e}) and
\beq
\label{5.18e}
\parallel \omega^k \ R(W) \parallel_2 \ \leq c_1 \ t^{-1-\lambda}
\eeq
\noi for all $t \in I$.\par

(1) Let $w$ be the solution of the equation (\ref{2.24e}) obtained in 
Proposition 5.2 part
(1). Then $w \in {\cal C}(I, H^k)$ and $w$ satisfies the estimates 
(\ref{5.5e}) and
\beq
\label{5.19e}
\parallel \omega^k (w(t) - W(t))\parallel_2 \ \leq C\ t^{-\lambda}
\eeq
\noi for all $t \in I$. \par

(2) Let $w_{t_0}$ be the solution of the equation (\ref{2.24e}) 
defined in Proposition 5.2 part
(1). When $t_0 \to \infty$, $w_{t_0}$ converges to $w$ strongly in 
$L^{\infty} ([1, \bar{T}],
H^{k'})$ for $0 \leq k' < k$ and in the weak $*$ sense in 
$L^{\infty}([1, \bar{T}], H^k)$ for
any $\bar{T} < \infty$.} \\

\noi {\bf Proof.} Part (1) will be proved together with the limiting 
properties stated in
Part (2). We know from Proposition 5.1 that $w_{t_0} \in {\cal C}(I, 
H^k)$. The main point of
the proof consists in estimating $q_{t_0} = w_{t_0} - W$ in $H^k$ 
uniformly in $t_0$ for $t
\leq t_0$. We know already from (\ref{5.10e}) that
\beq
\label{5.20e}
\parallel q_{t_0} (t) \parallel_2 \ \leq Y_0 \ t^{-\lambda_0}
\eeq
\noi for $t \leq t_0$, with $Y_0 = c_0\lambda_0^{-1}$. We next 
estimate $y \equiv \parallel
\omega^k q_{t_0}\parallel_2$. From (\ref{5.8e}) we obtain
\beq
\label{5.21e}
\left | \partial_t \parallel \omega^k q_{t_0}\parallel_2 \right | 
\leq t^{-1} \parallel
[\omega^k, B_0]q_{t_0}\parallel_2 \ + \ \parallel \omega^k 
R(W)\parallel_2 \eeq
\noi so that by Lemma 3.2, in the same way as in (\ref{4.7e}),
\bea
\label{5.22e}
\left | \partial_t \parallel \omega^k q_{t_0}\parallel_2 \right | 
&\leq& C t^{-1} \left (
\parallel \nabla B_0 \parallel_{\infty} \ \parallel \omega^{k-1} 
q_{t_0}\parallel_2 \ + \
\parallel \omega^k B_0\parallel_{3/\delta} \ \parallel 
q_{t_0}\parallel_r \right ) \nn \\
&&+\parallel \omega^k\ R(W)\parallel_2  \eea
\noi with $0 < \delta = \delta (r) < k \wedge 3/2$. Using the 
estimates (\ref{3.17e}), Sobolev
inequalities and interpolation with the help of (\ref{5.20e}) and the 
assumption (\ref{5.18e}),
we obtain
\bea
\label{5.23e}
|\partial_ty| &\leq& C\ b_0 \Big \{ Y_0^{1-\delta /k}\ y^{\delta /k}\ 
t^{k-1-\delta /3 -
\lambda_0 (1 - \delta /k)} \nn \\
&&+ Y_0^{1/k} \ y^{1- 1/k}\ t^{-\lambda_0/k} \Big \} + c_1 \
t^{-1-\lambda} \ .\eea
\noi We now define $Y \equiv \parallel t^{\lambda}y;L^{\infty}([1, 
t_0])\parallel$, substitute
that definition into the RHS of (\ref{5.23e}), integrate the latter 
between $t$ and $t_0$, and
obtain
\beq
\label{5.24e}
y \leq Cb_0 \left ( Y_0^{1- \delta/k} \ Y^{\delta /k}(\lambda + 
\nu_1)^{-1} \ t^{-\lambda
-\nu_1} + Y_0^{1/k} \ Y^{1-1/k}(\lambda + \nu_2)^{-1} \ t^{-\lambda - 
\nu_2} \right ) +
c_1\ \lambda^{-1} \ t^{-\lambda} \eeq
\noi where
\bea
\label{5.25e}
&&\nu_1 = (\lambda_0 - \lambda ) (1 - \delta /k) - k + \delta /3 \\
&&\nu_2 = (\lambda_0 - \lambda )/k - 1
\label{5.26e}
\eea
\noi provided $\lambda + \nu_1 > 0$ and $\lambda + \nu_2 > 0$. We 
impose in addition $\nu_1
\geq 0$, $\nu_2 \geq 0$, multiply (\ref{5.24e}) by $t^{\lambda}$, 
take the Supremum over $t$
and obtain
\beq
\label{5.27e}
Y \leq \lambda^{-1} \ C \ b_0 \left ( Y_0^{1-\delta /k} \ Y^{\delta 
/k} + Y_0^{1/k} \ Y^{1-1/k}
\right ) + \lambda^{-1} \ c_1 \eeq
\noi which is uniform in $t_0$. The condition $\nu_1 \geq 0$ reduces to
\beq
\label{5.28e}
\lambda_0 \geq \lambda + k + 2\delta k/3(k - \delta )
\eeq
\noi and can be satisfied for $\lambda_0 > \lambda  + k$ by taking 
$\delta$ sufficiently
small. It implies $\nu_2 > 0$. Changing the notation to $x = 
YY_0^{-1}$, $b = \lambda^{-1}
Cb_0$ and $c = \lambda^{-1} c_1 Y_0^{-1}$, we rewrite (\ref{5.27e}) as
\beq
\label{5.29e}
x \leq b \left ( x^{\delta /k} + x^{1 - 1/k} \right ) + c \ .
\eeq
\noi Assuming $\delta \leq k-1$ without loss of generality and using
$$x^{\theta} \leq \varepsilon x + \varepsilon^{-\theta/(1 - \theta )}$$
\noi for $x > 0$, $\varepsilon > 0$ and $0 < \theta < 1$, we obtain 
from (\ref{5.29e})
$$x \leq 2b \left ( \varepsilon x + \varepsilon^{1-k} \right ) + c$$
\noi and for $\varepsilon = (4b)^{-1}$
$$x \leq (4b)^k + 2c$$
\noi or equivalently
\beq
\label{5.30e}
Y \leq \left ( 4C\lambda^{-1} b_0 \right )^k Y_0 + 2 \lambda^{-1}\ c_1 \ ,
\eeq
\noi which completes the proof of the estimate of $q_{t_0}$ in $H^k$ 
uniformly in $t_0$ for $t
\leq t_0$. \par

Let now $\bar{T} < \infty$ and $J = [1, \bar{T}]$. We know from 
(\ref{5.30e}) that $w_{t_0}$
is estimated in $L^{\infty} (J, H^k)$ uniformly in $t_0$ for $t_0 
\geq \bar{T}$ and that
$w_{t_0}$ converges to $w$ in norm in $L^{\infty}(J, L^2)$ by 
Proposition 5.2 part (1). It
follows therefrom by a standard compactness argument that $w \in 
(L^{\infty} \cap {\cal C}_w)
(J, H^k)$, that $w$ also satisfies the estimate (\ref{5.30e}) and 
that $w_{t_0}$ converges to
$w$ in the topologies described in Part (2). Strong continuity of $w$ 
in $H^k$ follows from
Proposition 5.1. \par \nobreak \hfill
$\sq$\par

In order to complete the construction of the wave operators in $H^k$, 
we now have to construct
model functions $W$ satisfying the assumptions (\ref{5.4e}) and 
(\ref{5.18e}). In view of
Proposition 5.2 part (2), we restrict our attention to $W$ of the form
\beq
\label{5.31e}
W(t) = U^*(1/t) w_+
\eeq
\noi for some fixed $w_+ \in H^{k_+}$ and we take $k_+ > 3/2$. With 
that choice, we obtain
\beq
\label{5.32e}
R(W) = R (U^*(1/t)W) = - i \ t^{-1} \ B_0 \ U^*(1/t) w_+ \ .
\eeq
\noi However, with no further assumptions on $w_+$, we are restricted 
to $\lambda_0 \leq 1/2$
and consequently to $k < 1/2$. In fact, from (\ref{3.17e}) we obtain
\beq
\label{5.33e}
\parallel R(W)\parallel_2 \ \leq t^{-1} \parallel B_0 \parallel_2 \ 
\parallel U^*(1/t) w_+
\parallel_{\infty} \ \leq C b_0 \ t^{-3/2} |w_+|_{k_+} \ . \eeq
\noi Furthermore, from Lemma 3.2 and (\ref{3.17e}) we obtain for $k < 1/2$

$$\parallel \omega^k R(W) \parallel_2 \ \leq C\ t^{-1}\left \{ 
\parallel \omega^k
B_0\parallel_2 \ \parallel U^*(1/t)w_+ \parallel_{\infty} \ + \ 
\parallel B_0\parallel_r \
\parallel \omega^k U^*(1/t) w_+ \parallel_{3/\delta} \right \}$$
\beq
\label{5.34e}
\hskip - 6.5 truecm\leq Cb_0\ t^{-3/2+k} \ |w_+|_{k_+}
\eeq
\noi with $0 < \delta = \delta (r) \leq 3k < 3/2$. Together with an 
extension of Proposition
5.2 to $k \leq 1/2$, which we have not performed, the estimates 
(\ref{5.33e}) (\ref{5.34e})
would allow us to complete the construction of the wave operators for 
$0 < k < 1/2$, with
$\lambda_0 = 1/2$ and $0 < \lambda < 1/2 - k$. \par

In order to cover higher values of $k$, and in particular for $k > 
1$, as will be needed for
the nonlinear system (1.1) (1.2), we shall need additional conditions on $w_+$
and $B_0$. We first exhibit a set of local sufficient conditions in 
the form of joint decay
estimates for $w_+$, and $B_0$, where the nonlocal operator 
$U^*(1/t)$ no longer appears. \\

\noi {\bf Lemma 5.1.} {\it Let $\lambda_0 > 0$ and let $\bar{m}$ be a 
nonnegative integer. Let
$B_0$ satisfy the estimates (\ref{3.17e}) for $0 \leq m \leq 
\bar{m}$. Let $w_+ \in H^{k_+}$
with $k_+ \geq 2 \lambda_0 \vee \bar{m}$ and let $a_+ = |w_+|_{k_+}$. 
Assume that $B_0$ and
$w_+$ satisfy the estimates
\beq
\label{5.35e}
\parallel \left ( \partial^{\alpha_1} B_0 \right ) \left ( 
\partial^{\alpha_2}w_+\right )
\parallel_2 \ \leq b_1 \ t^{-\lambda_0 + |\alpha_1| + |\alpha_2|/2} \eeq
\noi for all multi-indices $\alpha_1$, $\alpha_2$ with $0 \leq 
|\alpha_1| \leq \bar{m}$ and $0
\leq |\alpha_2| < 2 \lambda_0$, and for all $t \geq 1$. Then the 
following estimates hold for
all $m$, $0 \leq m \leq \bar{m}$, and for all $t \geq 1$}
\bea
\label{5.36e}
\parallel \omega^m R(U^*(1/t) w_+ ) \parallel_2 &=& t^{-1} \parallel 
\omega^m \left ( B_0 \
U^*(1/t) w_+\right ) \parallel_2 \nn \\
&\leq& C \left ( b_1 + b_0 a_+\right ) t^{-1-\lambda_0+m} \ .
\eea
\vskip 3 truemm

\noi {\bf Proof.} By interpolation, it is sufficient to prove 
(\ref{5.36e}) for integer $m$. Let
$\alpha$ be a multi-index with $|\alpha | = m$. We estimate
\beq
\label{5.37e}
\parallel \partial^{\alpha} \left ( B_0\ U^*(1/t) w_+ \right 
)\parallel_2 \ \leq C \sum_{\alpha_1
+ \alpha_3 = \alpha}  \parallel \left ( \partial^{\alpha_1} B_0\right 
)  U^*(1/t) \
\partial^{\alpha_3} w_+ \parallel_2 \ . \eeq
\noi If $|\alpha_3| < 2 \lambda_0$, we expand $U^*(1/t)$ through the relation
$$\left | e^{ix} - \sum_{j \leq p} (j!)^{-1} (ix)^j \right | \leq 
2(p!)^{-1} \ |x|^{p+ \theta}$$
\noi with $p + \theta = \lambda_0 - |\alpha_3|/2$ and $0 < \theta 
\leq 1$, so that
\bea
\label{5.38e}
&&\parallel \left ( \partial^{\alpha_1} B_0\right ) U^*(1/t) 
\partial^{\alpha_3} w_+ \parallel_2
\ \leq C \sum_{j < \lambda_0 - |\alpha_3|/2} t^{-j}  \parallel \left 
( \partial^{\alpha_1}
B_0\right )  \Delta^j\ \partial^{\alpha_3} w_+ \parallel_2 \nn \\
&&+ C  \parallel \partial^{\alpha_1} B_0 \parallel_{\infty} \ 
t^{-(\lambda_0 - |\alpha_3|/2)}
\parallel \omega^{2\lambda_0} w_+ \parallel_2 \nn \\
&&\leq C(b_1 + b_0 a_+) t^{-\lambda_0 + |\alpha_1| + |\alpha_3|/2} = 
C(b_1 + b_0 a_+)
t^{-\lambda_0 + m - |\alpha_3|/2}
\eea
\noi by (\ref{5.35e}) and (\ref{3.17e}), which proves (\ref{5.36e}) 
in this case. \par

If $|\alpha_3| \geq 2 \lambda_0$, the last norm in (\ref{5.37e}) is 
estimated by the use of
(\ref{3.17e}) as
\bea
\label{5.39e}
&&\parallel \left ( \partial^{\alpha_1} B_0\right ) U^*(1/t) 
\partial^{\alpha_3} w_+ \parallel_2
\ \leq C \parallel \partial^{\alpha_1} B_0\parallel_{\infty} \ 
\parallel  \partial^{\alpha_3}
  w_+ \parallel_2 \nn \\ &&\leq C  \ b_0 \ a_+ \ t^{|\alpha_1|} \leq 
Cb_0 \ a_+ \ t^{m-2\lambda_0}
\eea
\noi since $|\alpha_1| = m - |\alpha_3| \leq m - 2 \lambda_0$, which 
completes the proof of
(\ref{5.36e}). \par \nobreak \hfill
$\sq$\par

We shall apply Lemma 5.1 with $\bar{m} = \{ k\}$, the smallest 
integer $\geq k$. Then
(\ref{5.36e}) with $m = 0$ and $m = k$ reduces to (\ref{5.4e}) and 
(\ref{5.18e}) with $\lambda =
\lambda_0 - k$ respectively. For $\lambda_0 > k \geq 1$, one can take 
$k_+ = 2 \lambda_0$. \par

We now give sufficient conditions that ensure the assumption 
(\ref{5.35e}). We first remark
that (\ref{5.35e}) is trivially satisfied under suitable support 
properties of $w_+$ and
of the initial data $(A_+, \dot{A}_+)$ of the scalar field $A_0$ at 
time $t = 0$ (see
(\ref{2.3e})). In fact, assume that
\beq
\label{5.40e}
{\rm Supp}\ (A_+, \dot{A}_+) \subset \{ x : |x| \leq R \} \ .
\eeq
\noi Then, by the Huyghens principle
\beq
\label{5.41e}
{\rm Supp}\ A_0 \subset \{ (x,t) : | |x| - t|\leq R \}
\eeq
\noi so that
\beq
\label{5.42e}
{\rm Supp}\ B_0 \subset \{ (x,t) : | |x| - 1|\leq R/t \} \ .
\eeq
\noi If on the other hand
\beq
\label{5.43e}
{\rm Supp}\ w_+ \subset \{ x : | |x| - 1|\geq \eta \}
\eeq
\noi for some $\eta$, $0 < \eta < 1$, then $(\partial^{\alpha_1}B_0) 
\partial^{\alpha_2}w_+ =
0$ for $t \geq R/\eta$ for any multi-indices $\alpha_1$ and 
$\alpha_2$, which ensures
(\ref{5.35e}) in a trivial way. \par

We shall now give more general conditions that ensure (\ref{5.35e}), 
keeping the support
condition (\ref{5.43e}) on $w_+$, and relaxing the support condition 
(\ref{5.40e}) on $(A_+,
\dot{A}_+)$ to space decay conditions. \\

\noi {\bf Lemma 5.2.} {\it Let $\lambda_0 > 0$, $k_+ \geq 2 
\lambda_0$ and let $k_+ > 3/2$. Let
$w_+ \in H^{k_+}$. Let $\alpha_1$ be a multi-index.\par

(1) Let $B_0$ satisfy (\ref{3.17e}) for $m = |\alpha_1|$ and in addition
\beq
\label{5.44e}
\parallel \chi_0\ \partial^{\alpha_1} B_0 \parallel_2 \ \leq b_2 \ 
t^{-\lambda_0 + |\alpha_1|}
\eeq
\noi where $\chi_0$ is the characteristic function of the support of 
$w_+$. Then (\ref{5.35e})
holds for any multi-index $\alpha_2$, $0 \leq |\alpha_2| \leq 
2\lambda_0$. The constant in
(\ref{5.35e}) can be taken as $b_1 = C(b_0 \vee b_2) a_+$. \par

(2) Let $w_+$ satisfy the support property (\ref{5.43e}) and let 
$(A_+, \dot{A}_+)$ satisfy the
following conditions for all $R \geq R_0$ for some $R_0 >0$~:
\bea
\label{5.45e}
\parallel \chi(|x| \geq R)  \partial^{\alpha_1} A_+ \parallel_2 
&\leq& C \ R^{-\lambda_0
+1/2} \ , \\
\label{5.46e}
\parallel \chi(|x| \geq R)  \dot{A}_+ \parallel_{6/5}  &\leq& C \ 
R^{-\lambda_0 +1/2} \quad
\hbox{if} \ \alpha_1 = 0 \ , \\
  \parallel \chi(|x| \geq R) \partial^{\alpha '_1}  \dot{A}_+ 
\parallel_{2}  &\leq& C \
R^{-\lambda_0 +1/2} \quad \hbox{if} \ \alpha_1 \not= 0 \ ,
\label{5.47e}
\eea
\noi where $\alpha '_1$ is a multi-index satisfying $\alpha '_1 \leq 
\alpha_1$, $|\alpha '_1| =
|\alpha_1 | - 1$, and where $\chi (|x| \geq R)$ is the characteristic 
function of $\{x: |x|
\geq R\}$. Then (\ref{5.44e}) holds.} \\

\noi {\bf Proof.} \underbar{Part (1)}. We estimate by the H\"older 
inequality and interpolation
between (\ref{3.17e}) and (\ref{5.44e})
\bea
\label{5.48e}
&&\parallel \left ( \partial^{\alpha_1} B_0\right ) \left ( 
\partial^{\alpha_2} w_+\right )
\parallel_2 \ \leq \ \parallel \chi_0 \left ( \partial^{\alpha_1} 
B_0\right ) \parallel_r \ \parallel
\partial^{\alpha_2} w_+ \parallel_{3/\delta} \nn \\
&&\leq (b_0 \vee b_2) t^{-\lambda_0 + |\alpha_1| + 2 \lambda_0 \delta 
/3} \parallel
\partial^{\alpha_2} w_+ \parallel_{3/\delta}\nn \\
&&= (b_0 \vee b_2) t^{-\lambda_0 + |\alpha_1| + |\alpha_2|/2} \parallel
\partial^{\alpha_2} w_+ \parallel_{3/\delta}
\eea
\noi where $\delta = \delta (r) = 3|\alpha_2|/4\lambda_0$, so that $0 
\leq \delta \leq 3/2$.
The last norm in (\ref{5.48e}) is estimated by $|w_+|_{k_+}$ through 
Sobolev inequalities
since $|\alpha_2| + 3/2 - \delta \equiv 3/2 + |\alpha_2|(1 - 
3/4\lambda_0)$ ranges from 3/2
to $2\lambda_0$ when $|\alpha_2|$ ranges from 0 to $2\lambda_0$. \par

\noi \underbar{Part (2)}. Using the support properties of $w_+$ and 
returning to the variable
$A_0$, we see that (\ref{5.44e}) is implied by
\beq
\label{5.49e}
\parallel \chi (| |x| - t | \geq \eta t ) \partial^{\alpha_1} A_0 (t) 
\parallel_2 \ \leq b_2\
t^{-\lambda_0 + 1/2} \ . \eeq
\noi Let now $R > 0$, let $\chi_1 \in {\cal C}^{\infty} 
({I\hskip-1truemm R}^3)$, $0 \leq
\chi_1 \leq 1$, $\chi_1 (x) = 0$ for $|x| \leq 1$, $\chi_1(x) = 1$ 
for $|x| \geq 2$ and let
$\chi_R(x) = \chi (x/R)$. Let $\widetilde{A}_R$ be the solution of 
the wave equation $\sq
\widetilde{A}_R = 0$ with initial data
$$\left ( \widetilde{A}_R, \partial_t \widetilde{A}_R\right )(0) = 
\left ( \chi_R
\ \partial^{\alpha_1} A_+, \chi_R \ \partial^{\alpha_1} \dot{A}_+ 
\right ) \quad \hbox{at}\ t = 0
\ .$$
\noi By the Huyghens principle, $\widetilde{A}_R(t) = 
\partial^{\alpha_1} A_0(t)$ for $||x| -
t| \geq 2R$ so that
\beq
\label{5.50e}
\parallel \chi (|| x| - t| \geq 2R ) \partial^{\alpha_1} A_0(t) 
\parallel_2 \ = \
\parallel \chi (||x| - t| \geq 2R) \widetilde{A}_R(t) \parallel_2 \ 
\leq \ \parallel
\widetilde{A}_R(t) \parallel_2 \ .
  \eeq
\noi It follows now from (\ref{3.16e}) that
\beq
\label{5.51e}
\parallel \widetilde{A}_R(t) \parallel_2 \ \leq \ \parallel \chi_R
\ \partial^{\alpha_1}A_+\parallel_2 \ + \ \parallel \omega^{-1} \ 
\chi_R \ \partial^{\alpha_1}
\dot{A}_+ \parallel_2 \ . \eeq
\noi If $\alpha_1 = 0$, we estimate the last norm in (\ref{5.51e}) by
\beq
\label{5.52e}
\parallel \omega^{-1} \ \chi_R \ \dot{A}_+ \parallel_2 \ \leq C 
\parallel \chi_R \ \dot{A}_+
\parallel_{6/5} \ . \eeq
\noi If $\alpha_1 \not= 0$, we rewrite $\partial^{\alpha_1} = 
\partial_j \partial^{\alpha '_1}$
and estimate
\bea
\label{5.53e}
\parallel \omega^{-1} \ \chi_R \ \partial^{\alpha_1} \dot{A}_+ 
\parallel_2 &\leq&  \parallel
\omega^{-1} \ \partial_j\  \chi_R \ \partial^{\alpha '_1} \dot{A}_+ 
\parallel_2 \ + \parallel
\omega^{-1} (\partial_j \ \chi_R) \partial^{\alpha '_1} \dot{A}_+ 
\parallel_2 \nn \\
&\leq& \parallel \chi_R \ \partial^{\alpha '_1} \dot{A}_+ \parallel_2\ + \ C
\parallel \partial_j\  \chi_1\parallel_3 \ \parallel \chi (|x| \geq 
R) \ \partial^{\alpha '_1}
\dot{A}_+ \parallel_2 \nn \\
\eea
\noi by the Sobolev and H\"older inequalities. Collecting 
(\ref{5.50e})-(\ref{5.53e}) and using
the assumption (\ref{5.45e})-(\ref{5.47e}), we obtain
$$\parallel \chi (| | x| - t| \geq 2R) \partial^{\alpha_1} A_0 (t) 
\parallel_2 \ \leq C \
R^{-\lambda_0 + 1/2}$$
\noi from which (\ref{5.49e}) follows by taking $2R = \eta t$.\par 
\nobreak \hfill
$\sq$\par

Collecting the previous results essentially yields the wave operators 
in $H^k$ for the equation
(\ref{2.24e}) in the form of Proposition 5.4 below. In that 
proposition, we have kept the
assumptions on $B_0$ in the implicit form of the estimates 
(\ref{3.17e}) and (\ref{5.35e}). If
so desired, those assumptions can be replaced by sufficient 
conditions on $(w_+, A_+,
\dot{A}_+)$ by the use of Lemmas 3.5 and 5.2. \\

\noi {\bf Proposition 5.4.} {\it Let $k \geq 1$, $k_+ > 2k$, let 
$\lambda_0$ and $\lambda$
satisfy $\lambda > 0$, $k + \lambda < \lambda_0 \leq k_+/2$. Let $w_+ 
\in H^k$, let $B_0$
satisfy the estimates (\ref{3.17e}) for $0 \leq m \leq k$, and let 
$(w_+, B_0)$ satisfy the
estimates (\ref{5.35e}) for all multi-indices $\alpha_1$, $\alpha_2$ 
with $0 \leq |\alpha_1| \leq
\bar{k}$ and $0 \leq |\alpha_2| < 2\lambda_0$, where $\bar{k}$ is the 
smallest integer $\geq
k$. Then the equation (\ref{2.24e}) has a unique solution $w \in 
{\cal C}([1, \infty ), L^2)$
such that
$$\parallel w(t) - w_+ \parallel_2 \ \to 0 \quad \hbox{when} \ t \to 
\infty \ .$$
\noi Furthermore $w \in {\cal C}([1, \infty ), H^k)$ and $w$ satisfies the estimates
\bea
\label{5.54e}
&&\parallel w(t) - U^*(1/t) w_+ \parallel_2 \ \leq C \ t^{-\lambda_0}
\\
\label{5.55e}
&&\parallel \omega^k (w(t) - U^*(1/t) w_+) \parallel_2 \ \leq C \
t^{-\lambda} \eea
\noi for all $t \geq 1$.\\}

\noi {\bf Proof.} The results follow from Propositions 5.2 and 5.3 
and from Lemma 5.1.\par \nobreak \hfill
$\sq$\par

The existence of the wave operators for $u$ in the usual sense at the 
corresponding level of
regularity is an easy consequence of Proposition 5.4. We refrain from 
giving a formal statement
at this stage. The same question will be considered in Section 8 in 
the case of the interacting
system (1.1) (1.2).

\mysection{Cauchy problem at infinity for the auxiliary system}
\hspace*{\parindent} In this section we begin the construction of the 
wave operators for the
auxiliary system (\ref{2.20e}) by solving the Cauchy problem at 
infinity for that system in
the difference form (\ref{2.30e}), for large or infinite initial 
time, and for a given choice
of $(W, S)$ satisfying a number of a priori estimates. The 
construction of $(W,S)$ satisfying
those estimates is deferred to the next section. In the same spirit 
as in Section 4, we
solve the system (\ref{2.30e}) in two steps. We first solve the 
linearized version of that
system (\ref{2.31e}), thereby defining a map $\Gamma : (q, \sigma ) 
\to (q', \sigma ')$. We
then show that this map is a contraction in suitable norms on a 
suitable set. \par

The basic tool of this section again consists of a priori estimates 
for suitably regular
solutions of the linearized system (\ref{2.31e}). In order to handle 
efficiently a
non-vanishing $B_0$, those estimates have to be much more elaborate 
than those of Section
4. \par

We first estimate  a single solution of the linearized system 
(\ref{2.31e}) at the level
of regularity where we shall eventually solve the auxiliary system 
(\ref{2.30e}). \\

\noi {\bf Lemma 6.1.} {\it Let $1 < k \leq \ell$, $\ell > 3/2$ and 
$\beta >0$. Let $I
\subset [1, \infty)$ be an interval and let $t_1 \in \bar{I}$. Let 
$B_0$ satisfy the
estimates (\ref{3.17e}) for $0 \leq l \leq k$. Let $(U(1/t))W,S) \in 
{\cal C} (I_+, X^{k+1,
\ell + 1}) \cap {\cal C}^1(I_+, X^{k,\ell})$ with $W \in 
L^{\infty}(I_+, H^{k+1})$ and let
\beq
\label{6.1e}
a = \parallel W; L^{\infty}(I_+, H^{k+1})\parallel \ .
\eeq
\noi Let $(q, \sigma )$, $(q', \sigma ') \in {\cal C} (I, 
X^{k,\ell})$ with $q \in L^{\infty}(I,
H^k) \cap L^2(I,L^2)$ if $t_1 = \infty$, and let $(q', \sigma ')$ be 
a solution of the system
(\ref{2.31e}) in $I$. Then the following estimates hold for all $t \in I$~:
\bea
\label{6.2e}
\left | \partial_t \parallel q'\parallel_2 \right | &\leq& C \Big \{ 
t^{-2} a \parallel \nabla
\sigma \parallel_2 + t^{-1-\beta} \ a^2 \ I_0\left ( \parallel q 
\parallel_2\right ) \nn \\
&&+ t^{-1}\ a \ I_{-1}\left ( \parallel q\parallel_2\ \parallel 
q\parallel_3 \right ) \Big \} +
\ \parallel R_1 (W, S) \parallel_2 \ , \eea
\bea
\label{6.3e}
&&\left | \partial_t \parallel \omega^k q'\parallel_2 \right | \leq C 
\Big \{ b_0 \left (
\parallel \omega^{k-1}q'\parallel_2 \ + t^{k-1-\delta /3} \ \parallel 
q'\parallel_r \right )\nn
\\ &&+ t^{-2}\ a \left ( \parallel \omega^k \nabla \sigma \parallel_2 
\ + \ \parallel \sigma
\parallel_{\infty} \ + \ \parallel \nabla \sigma \parallel_3 \right ) \nn \\
&& t^{-2} \left ( \parallel \nabla s \parallel_{\infty} \ \parallel 
\omega^k q'\parallel_2 \ +
\ \parallel \omega^{[k\vee 3/2]}\nabla s \parallel_2 \ \parallel 
\omega^{[k\wedge
3/2]}q'\parallel_2 \right . \nn \\
&&\left . + \chi (k > 3/2) \parallel \omega^k \nabla s \parallel_2 \ \parallel
q'\parallel_{\infty} \right ) \nn \\
&&+ t^{-1} \ a^2 \left ( I_{k-1} \left ( \parallel \omega^{k-1} q 
\parallel_2 \ + \ \parallel q
\parallel_2 \right ) + I_0 \left ( \parallel q \parallel_2 \right ) + 
\ \parallel
\omega^{k-1}q'\parallel_2 \ + \ \parallel q'\parallel_2 \right ) \nn \\
&&+ t^{-1} \ a \left ( I_{k-1} \left ( \parallel \omega^{k} q \parallel_2  \
\parallel q \parallel_3 \right ) + I_0 \left ( \parallel \nabla q 
\parallel_2\ \parallel q
\parallel_3 \right )  \right . \nn \\
&&+ I_{1/2} \left (  \parallel \omega^{1/2}q\parallel_2 \right )
\parallel \omega^k q'\parallel_2\ + \left . I_{k-1/2} \left ( 
\parallel \omega^{k-1/2} q
\parallel_2 \ + \ \parallel q \parallel_2 \right ) \ \parallel \nabla 
q' \parallel_2 \right )
\nn \\ && + t^{-1} \left ( I_{1/2} \left ( \parallel \nabla q 
\parallel_2^2 \right )
\parallel \omega^k q' \parallel_2 \ + I_{k-1/2} \left (  \parallel 
\omega^{k}q\parallel_2
\  \parallel \nabla q\parallel_2 \right ) \parallel \nabla 
q'\parallel_2 \right ) \Big \} \nn \\
&&+\ \parallel \omega^k R_1(W,S) \parallel_2
  \eea
\noi where $s = S + \sigma$ and $0 < \delta = \delta (r) \leq [k \wedge 3/2]$,
\bea
\label{6.4e}
&&\left | \partial_t \parallel \omega^m \nabla \sigma '\parallel_2 
\right | \leq C \ t^{-2} \Big
\{ \parallel \nabla s\parallel_{\infty} \ \parallel \omega^m \nabla 
\sigma ' \parallel_2 \
+ \ \parallel \omega^m \nabla s \parallel_2 \ \parallel \nabla \sigma 
' \parallel_{\infty}\nn \\
&&\parallel \omega^m \nabla \sigma\parallel_{2} \ \parallel \nabla S
\parallel_{\infty} \ + \ \parallel \sigma \parallel_{\infty} \ 
\parallel \omega^m \nabla^2
S\parallel_{2} \Big \}\nn \\
&&+ C \Big \{ t^{-1+ \beta (m+1)} a\ I_0\left ( \parallel 
q\parallel_2 \right ) + t^{-1 + \beta
(m + 5/2)} \ I_{-3/2} \left ( \parallel q \parallel_2^2 \right ) \Big \} \nn \\
&&+ \ \parallel \omega^m \nabla R_2(W,S)\parallel_2
\eea
\noi for $0 \leq m \leq \ell$,
$$\left | \partial_t \parallel \nabla \sigma '\parallel_2 \right | 
\leq C \ t^{-2} \Big
\{ \parallel \nabla s\parallel_{\infty} \ \parallel \nabla \sigma ' 
\parallel_2 \
+ \ \parallel \nabla \sigma \parallel_2 \left ( \parallel \nabla S 
\parallel_{\infty}\ +
\parallel \omega^{3/2} \nabla S \parallel_2 \right ) \Big \} $$
$$+ C \Big \{ t^{-1+\beta}\ a\ I_0 \left ( \parallel q 
\parallel_2\right ) + t^{-1+5\beta /2}\
I_{-3/2} \left ( \parallel q \parallel_2^2 \right ) \Big \} + \ 
\parallel \nabla
R_2(W,S)\parallel_2 \eqno(6.4)_0$$
\noi where the time parameter is $t_1$ in all the estimating 
functions $I_m$, and the
superscript $t_1$ is omitted for brevity.}\\

\noi {\bf Remark 6.1.} All the norms of $(q, \sigma )$ and $(q', 
\sigma ')$ that appear in
(\ref{6.2e})-(\ref{6.4e}) are controlled by the norms in $X^{k,\ell}$ 
through Sobolev
inequalities. Furthermore all the integrals $I_m$ are convergent if 
$t_1 = \infty$. This
follows from boundedness of $q$ in $H^k$ in all cases where $m > - 
1/2$, namely in all cases
but two. The exceptions are
$$I_{-3/2} \left ( \parallel q\parallel_2^2 \right ) = 
\int_1^{\infty} d\nu \parallel q(\nu
t)\parallel_2^2$$
\noi in (\ref{6.4e}) and
$$I_{-1} \left ( \parallel q\parallel_2 \ \parallel q\parallel_3\right ) \leq C
\int_1^{\infty} d\nu \ \nu^{-1/2} \  \parallel q\parallel_2^{3/2} \ 
\parallel \nabla
q\parallel_2^{1/2}$$
$$\leq C\parallel \nabla q ; L^{\infty}([t, \infty ), 
L^2)\parallel^{1/2} \ \parallel q; L^2
((t, \infty ), L^2)\parallel^{3/2}$$
  \noi in (\ref{6.2e}), both of which are controlled under the
additional assumption that $q \in L^2(I,L^2)$. \par

Finally it is easy to see by estimates similar to, but simpler than, 
those of Lemma 4.1 that
all the norms of the remainders $R_1$ and $R_2$ that occur in 
(\ref{6.2e})-(\ref{6.4e}) are
finite under the assumptions made on $(W,S)$. \\

\noi {\bf Proof of Lemma 6.1.} In all the proof, the time superscript 
in $B_S$, $B_L$ and in
the various $I_m$ is omitted, except in dubious cases. That time 
superscript is in general
$t_1$, except in $B_S(W,W)$ where it is $\infty$. \par

\noi \underbar{Proof of (6.2)}. From (\ref{2.31e}), we estimate
\beq
\label{6.5e}
\left | \partial_t \parallel q'\parallel_2 \right | \leq t^{-2} 
\parallel Q(\sigma, W)
\parallel_{2} \ + t^{-1} \parallel \left ( B_S(q,q) + 2B_S(W,q) 
\right ) W\parallel_2 \ + \
\parallel R_1(W,S)\parallel_2 \ .
\eeq
\noi We next estimate by (\ref{3.8e}) (\ref{3.10e}) and Sobolev inequalities
\bea
\label{6.6e}
\hskip - 1 truecm &&\parallel Q(\sigma, W) \parallel_{2} \leq C 
\parallel \nabla \sigma
\parallel_{2} \left ( \parallel \nabla W \parallel_{3} + \parallel 
W\parallel_{\infty} \right )
\leq C a  \parallel \nabla \sigma \parallel_{2}  , \\
\label{6.7e}
\hskip - 1 truecm &&\parallel B_S(W,q) W\parallel_{2}  \leq C 
t^{-\beta}  \parallel W
\parallel_{\infty} \ I_0 \left ( \parallel W \parallel_{\infty } 
\parallel q\parallel_{2}
\right ) \leq C  a^2  t^{-\beta}  I_0\left (\parallel q \parallel_{2} 
\right )  , \\
\hskip - 1 truecm &&\parallel B_S(q,q) W\parallel_{2} \leq C 
\parallel W \parallel_{\infty}
I_{-1} \left ( \parallel q \parallel_{2 } \parallel q\parallel_{3} 
\right ) \leq
C a I_{-1}\left (\parallel q \parallel_{2} \ \parallel q 
\parallel_{3}\right )  .
\label{6.8e}
\eea
\noi Substituting (\ref{6.6e}) (\ref{6.7e}) (\ref{6.8e}) into 
(\ref{6.5e}) yields
(\ref{6.2e}).\par

\noi \underbar{Proof of (6.3)}. From (\ref{2.31e}), we estimate
\bea
\label{6.9e}
&&\left | \partial_t \parallel \omega^k q'\parallel_2 \right | \leq t^{-1}
\parallel [\omega^{k}, B_0]q'\parallel_2 \ + t^{-2} \left (  \parallel
[\omega^k,s]\cdot \nabla q'\parallel_2 \right .\nn \\
&&\left . + \parallel (\nabla \cdot s)\omega^k q' \parallel_2 \ + \ 
\parallel \omega^k \left (
(\nabla \cdot s)q'\right ) \parallel_{2} \ + \ \parallel \omega^k 
Q(\sigma, W) \parallel_2
\right ) \nn \\ && t^{-1} \left ( \parallel [\omega^k ,B_S(w,w)]q' 
\parallel_{2} \ + \ \parallel
\omega^k \left ( B_S(q,q) + 2B_S(q,W)\right ) W \parallel_2 \right )  \nn \\
&&+\ \parallel \omega^k \ R_1(W,S) \parallel_2
  \eea
\noi and we estimate the various terms in the RHS successively. \par

The contribution of $B_0$ is estimated by Lemma 3.2 and by 
(\ref{3.17e}) exactly as in Section
4 (see (\ref{4.7e})) and yields
\beq
\label{6.10e}
\parallel [\omega^k, B_0]q'\parallel_2\ \leq \ C\ b_0 \left ( t 
\parallel \omega^{k-1}
q'\parallel_2 \ + t^{k-\delta /3} \ \parallel q'\parallel_r \right ) \ .\eeq
\indent The contribution of $Q(s, q')$ is estimated by (\ref{4.8e}) as
\bea
\label{6.11e}
&&\parallel [\omega^k, s]\cdot \nabla q'\parallel_2\ + \ \parallel
(\nabla \cdot s) \omega^k q'\parallel_2 \ + \ \parallel \omega^k 
\left ( (\nabla \cdot
s)q'\right ) \parallel_2 \nn \\
&&\leq C \Big \{ \parallel \nabla s\parallel_{\infty} \ \parallel
\omega^k q'\parallel_2 \ + \ \parallel \omega^{[k\vee 3/2]} \nabla s 
\parallel_2 \ \parallel
\omega^{[k\wedge 3/2]} q' \parallel_2\nn \\
&&+ \chi (k > 3/2) \parallel \omega^k \nabla s \parallel_2 \ 
\parallel q'\parallel_{\infty} \Big
\} \ .
\eea
\indent The contribution of $Q(\sigma , W)$ is estimated by Lemma 3.2 
and Sobolev inequalities
as \bea
\label{6.12e}
&&\parallel \omega^kQ(\sigma , W)\parallel_2\ \leq C \Big \{ 
\parallel\sigma \parallel_{\infty}
\ \parallel \omega^k \nabla W\parallel_2 \ + \  \parallel \omega^k 
\sigma \parallel_6\ \parallel
\nabla W\parallel_3\nn \\
&&+ \ \parallel \omega^k\nabla \sigma \parallel_{2} \ \parallel
W\parallel_{\infty} \ + \ \parallel \nabla \sigma \parallel_3 \ \parallel
\omega^kW \parallel_6 \Big \}\nn \\
&&\leq C\ a \Big \{ \parallel \omega^k \nabla \sigma \parallel_2 \ + 
\ \parallel \sigma
\parallel_{\infty} \ + \ \parallel \nabla \sigma \parallel_3 \Big \} \ .  \eea
\noi The contribution of $B_S$ with $w = W + q$ yields a number of 
terms which we order by
increasing number of $q$ or $q'$ occurring therein. We first expand
$$B_S^{t_1, \infty}(w, w) = B_S^{\infty}(W,W) + 2B_S^{t_1}(W,q) + 
B_S^{t_1}(q,q) \ .$$
\noi By Lemma 3.2 and Sobolev inequalities, we estimate
\bea
\label{6.13e}
&&\parallel [\omega^k, B_S(W,W)]q'\parallel_2\ \leq C \Big \{ 
\parallel\nabla B_S(W,W)
\parallel_{\infty}  \ \parallel \omega^{k-1} q' \parallel_2 \nn \\
&&+ \ \parallel \omega^{k+3/2-\delta} B_1(W,W) \parallel_{2} \ \parallel
q'\parallel_r \Big \}
\eea
\noi where we take $0 < \delta = \delta (r) = (k-1)\wedge 1/2$ so that
$$\parallel q'\parallel_r \ \leq C \left ( \parallel \omega^{k-1} 
q'\parallel_2 \ + \ \parallel
q'\parallel_2 \right ) \ .$$
\noi Furthermore
$$\parallel \nabla B_S(W,W) \parallel_{\infty} \ \leq C \parallel \nabla^2
B_1(W,W)\parallel_2^{1-\theta} \ \parallel \omega^{k+3/2} \ 
B_1(W,W)\parallel_2^{\theta}$$
\noi with $\theta = 1/(2k-1)$, and by (\ref{3.10e}) and Lemma 3.2
$$\parallel \omega^{m+1} B_1(W,W) \parallel_{2} \ \leq C \ I_m \left 
( \parallel \omega^m
W\parallel_2\ \parallel W\parallel_{\infty}\right ) \leq C\ a^2$$
\noi which we use with $m = 1$, $k + 1/2$ and $k + 1/2 - \delta$. 
Substituting those estimates
into (\ref{6.13e}) yields
\beq
\label{6.14e}
\parallel [\omega^k, B_S(W,W)]q'\parallel_2\ \leq C \ a^2 \left ( 
\parallel\omega^{k-1}
q'\parallel_{2}  \ + \ \parallel  q' \parallel_2 \right ) \  .
\eeq
\indent In a similar way, we estimate
\bea
\label{6.15e}
&&\parallel [\omega^k, B_S(W,q)]q'\parallel_2\ \leq C \Big \{ 
\parallel\omega^{3/2} B_1(W,q)
\parallel_{2}  \ \parallel \omega^{k} q' \parallel_2 \nn \\
&&+ \ \parallel \omega^{k+1/2} B_1(W,q) \parallel_{2} \ \parallel
\nabla q'\parallel_2 \Big \} \ .
\eea
\noi By Lemma 3.2 again and by (\ref{3.10e})
\begin{eqnarray*}
&&\parallel \omega^{m+1}B_1(W, q)\parallel_2 \ \leq I_m \left ( 
\parallel \omega^m Wq
\parallel_2 \right ) \ ,\\
&&\parallel \omega^{m}W q\parallel_2 \ \leq C \left ( \parallel W 
\parallel_{\infty} \  \parallel
\omega^mq \parallel_{2} \ + \ \parallel \omega^{m+3/2-\delta} W 
\parallel_{2} \ \parallel
q\parallel_r \right )
\end{eqnarray*}
\noi with $0 < \delta = \delta (r) = m \wedge 1/2$, so that for $m 
\leq k - 1/2$
\beq
\label{6.16e}
\parallel \omega^{m+1}B_1(W,q) \parallel_2 \ \leq C \ a \ I_m \left ( 
\parallel \omega^m q
\parallel_{2} \  + \ \parallel q\parallel_{2} \right ) \ ,
\eeq
\noi where $\parallel q\parallel_2$ can be omitted for $m \leq 1/2$. 
Substituting (\ref{6.16e})
with $m = 1/2$ and $m = k - 1/2$ into (\ref{6.15e}) yields
\bea
\label{6.17e}
&&\parallel [\omega^{k}, B_S(W,q)]q' \parallel_2 \ \leq C \ a \Big \{ 
I_{1/2} \left ( \parallel
\omega^{1/2} q \parallel_{2} \right )  \parallel \omega^k 
q'\parallel_{2}  \nn \\
&&+ I_{k-1/2} \left ( \parallel \omega^{k-1/2} q \parallel_2 \ + \ 
\parallel q \parallel_2
\right ) \parallel \nabla q'\parallel_2 \Big \} \ .
  \eea
\indent By Lemma 3.2 and Sobolev inequalities again, we next estimate
\bea
\label{6.18e}
&&\parallel [\omega^{k}, B_S(q,q)]q' \parallel_2 \ \leq C  \Big \{ \parallel
\omega^{3/2} B_1(q,q) \parallel_{2} \   \parallel \omega^k 
q'\parallel_{2} \nn \\
&&+ \parallel \omega^{k+1/2} B_1(q,q) \parallel_2 \ \parallel \nabla 
q' \parallel_2
\Big \}
  \eea
\noi followed by (see also (\ref{3.10e}))
\begin{eqnarray*}
&&\parallel \omega^{3/2} B_1(q, q)\parallel_2 \ \leq C  I_{1/2} \left 
(  \parallel
\nabla q \parallel_2^2 \right ) \\
&&\parallel \omega^{k+1/2} B_1(q, q)\parallel_2 \ \leq C I_{k-1/2} 
\left ( \parallel
\omega^k q \parallel_{2} \  \parallel \nabla q \parallel_{2} \right ) 
\end{eqnarray*}
\noi so that
\bea
\label{6.19e}
&&\parallel [\omega^{k}, B_S(q,q)]q' \parallel_2 \ \leq C  \Big \{ 
I_{1/2} \left ( \parallel
\nabla q \parallel_{2}^2 \right )   \parallel \omega^k q'\parallel_{2} \nn \\
&&+ I_{k-1/2} \left ( \parallel \omega^{k} q\parallel_2 \ \parallel \nabla q
\parallel_2\right ) \parallel \nabla q'\parallel_2 \Big \} \ .
  \eea
\indent We now turn to the second contribution of $B_S$ to 
(\ref{6.9e}). By Lemma 3.2 and
Sobolev inequalities again
\bea
\label{6.20e}
&&\parallel \omega^{k}(B_S(W,q)W) \parallel_2 \ \leq C  \Big \{ \parallel
\omega^{k} B_1(W,q) \parallel_{2} \   \parallel W\parallel_{\infty} \nn \\
&&+ \ \parallel \nabla B_1(W,q) \parallel_2 \ \parallel 
\omega^{k+1/2} W \parallel_2
\Big \}
  \eea
\noi and by (\ref{6.16e}) with $m = k - 1$ and $m = 0$
\beq
\label{6.21e}
\parallel \omega^{k}(B_S(W,q)W) \parallel_2 \ \leq C  \ a^2\left ( 
I_{k-1} \left (  \parallel
\omega^{k-1} q \parallel_{2} \ + \   \parallel q\parallel_{2} \right 
) + I_0 \left (  \parallel
q\parallel_2 \right ) \right ) \ .
\eeq
\noi Similarly
\bea
\label{6.22e}
&&\parallel \omega^{k}(B_S(q,q)W) \parallel_2 \ \leq C \Big \{ \parallel
\omega^{k} B_1(q,q) \parallel_{2} \  \parallel W\parallel_{\infty}  \nn \\
&&+ \ \parallel \nabla B_1(q,q) \parallel_2 \ \parallel 
\omega^{k+1/2} W \parallel_2
  \Big \}
  \eea
\noi followed by (see the proof of (\ref{6.19e}))
\begin{eqnarray*}
&&\parallel \omega^{k} B_1(q, q)\parallel_2 \ \leq C \ I_{k-1} \left 
(  \parallel
\omega^k q\parallel_2 \  \parallel q\parallel_3 \right )  \\
&&\parallel \nabla B_1(q, q)\parallel_2 \ \leq C \ I_{0} \left (  \parallel
\nabla q\parallel_2 \  \parallel q\parallel_3 \right )
\end{eqnarray*}
\noi yields
\beq
\label{6.23e}
\parallel \omega^{k}(B_S(q,q)W) \parallel_2 \ \leq C  \ a\left ( 
I_{k-1} \left (  \parallel
\omega^{k} q \parallel_{2} \ \parallel q\parallel_{3} \right ) + I_0 
\left (  \parallel
\nabla q\parallel_2 \ \parallel q\parallel_{3} \right ) \right ) \ .
\eeq
\noi Substituting (\ref{6.10e}) (\ref{6.11e}) (\ref{6.12e}) 
(\ref{6.14e}) (\ref{6.17e})
(\ref{6.19e}) (\ref{6.21e}) and (\ref{6.23e}) into (\ref{6.9e}) and 
reordering the
contributions of $B_S$ by increasing powers of $(q,q')$ yields 
(\ref{6.3e}). \par

\noi \underbar{Proof of (6.4)}. From (\ref{2.31e}), we estimate
\bea
\label{6.24e}
&&\left | \partial_t \parallel \omega^{m+1} \sigma '\parallel_2 
\right | \leq t^{-2}
\left ( \parallel [\omega^{m+1}, s]\cdot \nabla \sigma '\parallel_2 \ 
+ \parallel
(\nabla \cdot s)\omega^{m+1} \sigma '\parallel_2 \right .\nn \\
&&\left . + \parallel \omega^{m+1} (\sigma \cdot \nabla S) 
\parallel_2 \right ) + t^{-1}
\parallel \omega^{m+2} \left ( B_L(q,q) + 2B_L(W,q)\right ) \parallel_{2}\nn \\
&&+\ \parallel \omega^{m+1} \
R_2(W,S) \parallel_2 \ .
  \eea
\noi We next estimate by Lemma 3.2 again
\beq
\label{6.25e}
\parallel [\omega^{m+1} , s]\cdot \nabla \sigma '\parallel_2 \ \leq C \left (
\parallel \nabla s\parallel_{\infty} \  \parallel \omega^{m+1} \sigma 
'\parallel_2 \ + \
\parallel \omega^{m+1} s\parallel_2 \ \parallel \nabla \sigma 
'\parallel_{\infty}\right )
\eeq
\beq
\label{6.26e}
\parallel \omega^{m+1} (\sigma \cdot \nabla S)\parallel_2 \ \leq C \left (
\parallel \omega^{m+1} \sigma \parallel_{2} \  \parallel \nabla 
S\parallel_{\infty} \
+ \ \parallel \sigma\parallel_{\infty} \ \parallel \omega^{m+1} \nabla
S\parallel_{2}\right ) \eeq
\noi while by (\ref{3.9e}) and Lemma 3.4
\bea
\label{6.27e}
\parallel \omega^{m+2} B_L(W,q) \parallel_2 &\leq& t^{\beta (m+1)} \ 
I_0 \left (
\parallel W \parallel_{\infty} \ \parallel q\parallel_{2} \right ) \nn \\
&\leq& t^{\beta (m+1)} \ a\ I_0 \left ( \parallel q\parallel_{2} \right )\ ,
\eea
\beq
\label{6.28e}
\parallel \omega^{m+2} B_L(q,q) \parallel_2 \leq C\  t^{\beta 
(m+5/2)} \ I_{-3/2} \left (
\parallel q \parallel_2^2\right ) \ .
\eeq
\noi Substituting (\ref{6.25e}) (\ref{6.26e}) (\ref{6.27e}) 
(\ref{6.28e}) into (\ref{6.24e})
yields (\ref{6.4e}). For $m = 0$, the term $\parallel \omega^m\nabla 
s\parallel_2 \ \parallel
\nabla \sigma '\parallel_{\infty}$ can be omitted, and the term 
$\sigma \cdot \nabla S$ can be
estimated in a sllightly different way, thereby leading to 
(6.4)$_0$.\par \nobreak \hfill
$\sq$\par

We next estimate the difference of two solutions of the linearized 
system (\ref{2.31e})
corresponding to two different choices of $(q, \sigma )$, but to the 
same choice of $(W, S)$.
As in Section 4, we estimate that difference at a lower level of 
regularity than the solutions
themselves. \\

\noi {\bf Lemma 6.2.} {\it Let $1 < k \leq \ell$, $\ell > 3/2$ and 
$\beta > 0$. Let $I \subset
[1 , \infty )$ be an interval and let $t_1 \in \bar{I}$. Let $B_0$ be 
sufficiently regular, for
instance $B_0 \in {\cal C}(I, H_3^k)$. Let $(W,S)$ satisfy the 
assumptions of Lemma 6.1. Let
$(q_i, \sigma_i)$, $(q'_i , \sigma '_i) \in {\cal C}(I, X^{k,\ell})$ 
with $q_i \in
L^{\infty}(I, H^k) \cap L^2(I, L^2)$, $i = 1,2$, if $t_1 = \infty$, 
and let $(q'_i, \sigma
'_i)$ be solutions of the system (\ref{2.31e}) associated with $(q_i, 
\sigma _i)$ and $(W,S)$.
Define $(q_{\pm}, \sigma_{\pm}) = (1/2) (q_1 \pm q_2, \sigma_1 \pm 
\sigma_2)$ and $(q'_{\pm},
\sigma '_{\pm}) = (1/2) (q'_1 \pm q'_2, \sigma '_1 \pm \sigma '_2)$. 
Then the following
estimates hold for all $t \in I$.
\bea
&&\left | \partial_t \parallel q'_- \parallel_2 \right | \leq C \Big 
\{ t^{-2} \ a \parallel
\nabla \sigma_-\parallel_2 + t^{-2} \left ( \parallel \omega^{[3/2 - 
k]_+} \nabla \sigma_-
\parallel_2 \ \parallel \omega^{[k\wedge 3/2]}q'_+ \parallel_2\right .\nn \\
&&\left . + \chi (k > 3/2) \parallel \nabla \sigma_-\parallel_2 \ \parallel
q'_+\parallel_{\infty} \right ) + t^{-1-\beta} \ a^2\ I_0\left ( 
\parallel q_-\parallel_2 \right
) \nn \\
&&+ t^{-1} \ a\left ( \parallel q'_+\parallel_3\ I_0 \left ( 
\parallel q_-\parallel_2
\right ) + I_{-1} \left ( \parallel q_+ \parallel_3 \ \parallel q_-\parallel_2
\right ) \right ) \nn \\
&&+ t^{-1} \parallel q'_+\parallel_6\ I_{-1/2} \left ( \parallel q_+\parallel_6
\ \parallel q_- \parallel_2 \right ) \Big \} \ ,
\label{6.29e}
\eea
\bea
&&\left | \partial_t \parallel \omega^m \nabla \sigma '_- \parallel_2 
\right | \leq C \
t^{-2} \Big \{ \parallel \nabla s_+\parallel_{\infty} \ \parallel
\omega^{m} \nabla \sigma '_- \parallel_2 \ + \ \parallel \nabla s'_+
\parallel_{\infty} \ \parallel \omega^m \nabla \sigma_-\parallel_2 \nn \\
&&+ \parallel \omega^{m-m'+3/2} \nabla s_+\parallel_2 \
  \parallel \omega^{m'} \nabla \sigma '_-\parallel_2 \ + \ \parallel 
\omega^{m+5/2 - \delta}
\nabla s'_+ \parallel_2\ \parallel \sigma_-\parallel_r  \Big \} \nn \\
&& + C \Big \{ t^{-1+ \beta (m+1)} \ a\ I_0 \left ( \parallel
q_- \parallel_2 \right ) + t^{-1+ \beta (m + 5/2)} I_{-3/2} \left ( 
\parallel q_+\parallel_2
\ \parallel q_-\parallel_2\right ) \Big \} \ ,
\label{6.30e}
\eea
$$\left | \partial_t \parallel \nabla \sigma '_- \parallel_2 \right | \leq C \
t^{-2} \Big \{ \parallel \nabla s_+\parallel_{\infty} \ \parallel
\nabla \sigma '_- \parallel_2 $$
$$+ \ \parallel \nabla \sigma_-
\parallel_{2} \left (  \parallel \nabla s'_+\parallel_{\infty} + 
\parallel \omega^{3/2} \nabla
s'_+ \parallel_2 \right ) \Big \}$$
$$ + C \Big \{ t^{-1+ \beta } \ a\ I_0 \left ( \parallel q_-
\parallel_2 \right ) + t^{-1+ 5\beta /2} I_{-3/2} \left (  \parallel 
q_+\parallel_2  \ \parallel
q_-\parallel_2\right ) \Big \} \eqno(6.30)_0$$ \noi where $s_+ = S + 
\sigma_+$, $s'_+ = S +
\sigma '_+$,  \bea
\label{6.31e}
&&0 \leq m \leq \ell_0 \ , \ m' = m \wedge 1/2\ , \ \delta = \delta 
(r) = [(m+1) \wedge 3/2] \
,\nn \\
&&\qquad [3/2 - k]_+ \leq \ell_0 \leq \ell - 1 \ ,
\eea
\noi and the superscript $t_1$ is again omitted in the estimating 
functions $I_m$.}\\

\noi {\bf Remark 6.2.} Under the assumptions made, all the norms in 
the first part of the RHS
of (\ref{6.30e}) are controlled by $|s_+|_{\ell}^{\dot{}}$, 
$|s'_+|_{\ell}^{\dot{}}$,
$|\sigma_-|_{\ell_0}^{\dot{}}$ and $|\sigma '_-|_{\ell_0}^{\dot{}}$. \\

\noi {\bf Proof.} Taking the difference of the system (\ref{2.31e}) 
for $(q'_i , \sigma '_i)$,
or equivalently and more simply rewriting (\ref{4.21e}) with $(w_-, 
s_-) = (q_-, \sigma_-)$,
$(w'_-, s'_-) = (q'_-, \sigma '_-)$, and accounting for the 
replacement of $B_S^{t_1}(w,w)$ by
$B_S^{t_1 , \infty}(w, w)$, we obtain
\beq
\label{6.32e}
\left \{ \begin{array}{l} \partial_t q'_- = i (2t^2)^{-1} \Delta q'_- 
+ t^{-2}\Big \{
Q(s_+, q'_-) + Q(\sigma_-, w'_+)\Big \} + it^{-1} B_0 q'_-\\
\\
+ i t^{-1} \Big \{ \left ( B_S^{t_1, \infty} (w_+, w_+) + 
B_S^{t_1}(q_-, q_-)\right )
q'_- + 2 B_S^{t_1} (w_+, q_-) w'_+ \Big \}\\
\\
\partial_t \sigma '_- = t^{-2} \left ( s_+ \cdot \nabla \sigma '_- + 
\sigma_- \cdot
\nabla s'_+ \right ) - 2 t^{-1} \nabla B_L^{t_1} (w_+, q_-) \ 
.\end{array} \right .
\eeq
\indent We first estimate $q'_-$. From (\ref{6.32e}) we obtain
\beq
\label{6.33e}
\left | \partial_t \parallel q'_- \parallel_2 \right | \leq t^{-2}
\parallel Q(\sigma_-, w'_+) \parallel_{2} \ + 2 t^{-1}  \parallel
B_S(w_+, q_-) w'_+ \parallel_2 \ .
\eeq
\noi We expand (\ref{6.33e}) by using
$$(w_+, s_+) = (W, S) + (q_+, \sigma_+) \ , \ (w'_+ , s'_+) = (W, S) 
+ (q'_+ , \sigma '_+)$$
\noi and we estimate the various terms successively. \par

 From (\ref{6.6e}) we obtain
\beq
\label{6.34e}
\parallel Q(\sigma_-, W) \parallel_2 \ \leq C\ a \parallel \nabla 
\sigma_-\parallel_2 \ .
\eeq
\noi By the same estimates as in the proof of (\ref{4.23e}), we next obtain
\bea
\label{6.35e}
\parallel Q(\sigma_-, q'_+) \parallel_2 &\leq& C\Big \{   \parallel 
\omega^{[3/2 -
k]_+} \nabla \sigma_- \parallel_2 \ \parallel \omega^{[k \wedge 3/2]} 
q'_+ \parallel_2 \nn \\
&&+ \chi (k > 3/2) \parallel \nabla \sigma_-\parallel_2 \ \parallel 
q'_+\parallel_{\infty} \Big
\} \ . \eea
\noi We next estimate by (\ref{6.7e}) (\ref{6.8e})
\beq
\label{6.36e}
\parallel B_S(W,q_-)W \parallel_2 \ \leq C\ a^2\ t^{-\beta} \ I_0 
\left (  \parallel q_-
\parallel_2 \right ) \  ,
\eeq
\beq
\label{6.37e}
\parallel B_S(q_+,q_-)W \parallel_2 \ \leq C\ a\ I_{-1} \left (  \parallel q_+
\parallel_3\ \parallel q_-
\parallel_2 \right ) \  .
  \eeq
\noi The remaining terms are new. Using (\ref{3.10e}) and Sobolev 
inequalities, we obtain
successively \beq
\label{6.38e}
\parallel B_S(W,q_-)q'_+ \parallel_2 \ \leq C\ a\parallel q'_+ 
\parallel_3 \ I_0 \left
(  \parallel q_- \parallel_2 \right ) \  ,
\eeq
\bea
\label{6.39e}
\parallel B_S(q_+,q_-)q'_+ \parallel_2 &\leq& C\parallel B_S(q_+, 
q_-)\parallel_3
\ \parallel q'_+ \parallel_6 \nn \\
&\leq& C\ I_{-1/2} \left ( \parallel q_+ \parallel_6 \
\parallel q_- \parallel_2 \right ) \parallel q'_+ \parallel_6 \ .
\eea
\noi Substituting (\ref{6.34e})-(\ref{6.39e}) into (\ref{6.33e}) 
yields (\ref{6.29e}). \par

We next estimate $\sigma '_-$. From (\ref{6.32e}) we obtain
\bea
\label{6.40e}
&&\left | \partial_t \parallel \omega^{m+1} \sigma '_- \parallel_2 
\right | \leq t^{-2}
\left ( \parallel [\omega^{m+1}, s_+]\cdot \nabla \sigma '_- \parallel_2 \ + \
\parallel (\nabla \cdot s_+)\omega^{m+1} \sigma '_- \parallel_2 \right .\nn \\
&&+\left .  \parallel \omega^{m+1}(\sigma_- \cdot \nabla s'_+ 
)\parallel_2 \right ) + 2t^{-1}
\parallel \omega^{m+2} (B_L(W,q_-) + B_L(q_+, q_-) ) \parallel_2
\eea
\noi and we estimate the various terms successively. \par

 From Lemma 3.2 and Sobolev inequalities, we obtain
\bea
\label{6.41e}
&&\parallel [\omega^{m+1} , s_+] \cdot \nabla \sigma '_- \parallel_2 \ \leq C
\left ( \parallel \nabla s_+ \parallel_{\infty} \
\parallel \omega^{m+1} \sigma '_- \parallel_2 \right .\nn \\
&&+ \left . \parallel \omega^{m+1}s_+ \parallel_{3/m'} \ \parallel 
\nabla \sigma '_-
\parallel_{r'}\right ) \nn \\
&&\leq C \left ( \parallel \nabla s_+\parallel_{\infty} \ \parallel 
\omega^{m} \nabla \sigma '_-
\parallel_2 \ +\ \parallel \omega^{m-m'+3/2} \nabla s_+ \parallel_2 \ 
\parallel \omega^{m'}
\nabla \sigma '_-  \parallel_2 \right )\eea
\noi with $m' = \delta (r') = m \wedge 1/2$, and
\bea
\label{6.42e}
&&\parallel \omega^{m+1} (\sigma_- \nabla s'_+) \parallel_2 \ \leq C
\left ( \parallel \nabla s'_+ \parallel_{\infty} \
\parallel \omega^{m} \nabla \sigma_- \parallel_2 \ + \parallel 
\omega^{m+1}\nabla s'_+
\parallel_{3/\delta} \ \parallel \sigma_-  \parallel_{r}\right )
\nn \\
&&\leq C \left ( \parallel \nabla s'_+\parallel_{\infty} \ \parallel 
\omega^{m} \nabla
\sigma_-  \parallel_2 \ +\ \parallel \omega^{m+5/2-\delta} \nabla 
s'_+ \parallel_2 \ \parallel
\sigma_-  \parallel_r \right )
\eea
\noi with $\delta = \delta (r) = [(m+1)\wedge 3/2]$. \par

The contribution of $B_L$ to (\ref{6.40e}) is estimated exactly as in 
the proof of (\ref{6.4e})
(see (\ref{6.27e}) and (\ref{6.28e})) by
\bea
\label{6.43e}
&&\parallel \omega^{m+2} \left ( B_L(W,q_-) + B_L(q_+, q_-)\right ) 
\parallel_2 \ \leq C
\left ( a\ t^{\beta (m+1)}\ I_0 \left ( \parallel q_-\parallel_2 
\right )\right . \nn \\
&&\left . \qquad + t^{\beta (m+5/2)} \ I_{-3/2} \left ( \parallel q_+ 
\parallel_2 \ \parallel
q_- \parallel_2 \right ) \right ) \ .
\eea
\noi Substituting (\ref{6.41e}) (\ref{6.42e}) (\ref{6.43e}) into 
(\ref{6.40e}) yields
(\ref{6.30e}) and (6.30$_0$), where one term from (\ref{6.41e}) can 
be omitted. \par \nobreak \hfill
$\sq$\par

We now begin the study of the Cauchy problem for the auxiliary system 
in the difference form
(\ref{2.30e}) and for that purpose we first study that problem for 
the linearized version
of that system. For finite initial time $t_0$, that problem is solved 
by a minor
modification of Proposition 4.1. The following proposition is simply 
a compilation of that
result and of Lemmas 6.1 and 6.2. \\

\noi {\bf Proposition 6.1.} {\it Let $1 < k \leq \ell$, $\ell > 3/2$ 
and $\beta > 0$. Let $I
\subset [1, \infty )$ be an interval and let $t_1 \in \bar{I}$. let 
$B_0$ satisfy the estimates
(\ref{3.17e}) for $0 \leq m \leq k$. Let $(U(1/t)W,S) \in {\cal C} (I_+ , X^{k+1, \ell + 1})
\cap {\cal C}^1(I_+, X^{k,\ell})$ with $W \in 
L^{\infty}(I_+, H^{k+1})$ and define $a$ by
(\ref{6.1e}). Let $(q, \sigma) \in {\cal C}(I, X^{k,\ell})$ with $q 
\in L^{\infty}(I, H^k) \cap
L^2(I, L^2)$ if $t_1 = \infty$. Let $t_0 \in I$ and let $(q'_0, 
\sigma '_0) \in X^{k,\ell}$.
Then the system (\ref{2.31e}) has a unique solution $(q', \sigma ') 
\in {\cal C}(I,
X^{k,\ell})$ with $(q', \sigma ')(t_0) = (q'_0 , \sigma '_0)$. That 
solution satisfies the
estimates (\ref{6.2e}) (\ref{6.3e}) (\ref{6.4e}) of Lemma 6.1 for all 
$t \in I$. Two such
solutions $(q'_i, \sigma '_i)$ associated with $(q_i, \sigma _i)$, $i 
= 1,2$ and with the same
$(W, S)$ satisfy the estimates (\ref{6.29e}) (\ref{6.30e}) of Lemma 
6.2 for all $t \in I$}.\\

We shall be eventually interested in solving the Cauchy problem for 
the auxiliary system
(\ref{2.30e}) with infinite initial time $t_0$. As a preliminary, we 
need to solve the same
problem for the linearized system (\ref{2.31e}). This is done in the 
following proposition,
which of course requires much stronger assumptions on the asymptotic 
behaviour in time of
$(W,S)$ and $(q, \sigma )$. With the study of the nonlinear system in 
view, we already make the
assumptions that will be needed for that purpose, although they could 
be slightly weakened for
the linear problem. Since we want to take $t_0 = \infty$, we also 
take $t_1 = \infty$.\\

\noi {\bf Proposition 6.2.} {\it Let $1 < k \leq \ell$, $\ell > 3/2$. 
Let $\beta$, $\lambda_0$
and $\lambda$ satisfy
\beq
\label{6.44e}
0 < \beta < 1 \quad , \quad \lambda > 0 \quad , \quad \lambda_0 > 
\lambda + k \quad , \quad
\lambda_0 > \beta (\ell + 1) \ .
\eeq
\noi Let $t_1 = \infty$, let $1 \leq T < \infty$, and $I = [T, \infty 
)$. Let $B_0$ satisfy the
estimates (\ref{3.17e}) for $0 \leq m \leq k$. Let $(W, S)$ satisfy 
the assumptions of
Proposition 6.1 with
\bea
\label{6.45e}
&&\qquad |W|_{k+1} \leq a \ , \\
&&\parallel \omega^m \nabla S \parallel_2 \ \leq b\ t^{1-\eta + \beta 
(m - 3/2)}
\label{6.46e}
\eea
\noi for some $\eta > 0$ and for $0 \leq m \leq \ell + 1$,
\bea
  \label{6.47e}
&&\parallel R_1(W, S) \parallel_2 \ \leq c_0\ t^{-1-\lambda_0} \quad 
, \quad \parallel \omega^k
R_1(W,S) \parallel_2 \ \leq c_1\ t^{-1-\lambda}\ , \\
&&\parallel \omega^m \nabla R_2 (W,S) \parallel_2 \ \leq c_2\ 
t^{-1-\lambda_0 + \beta (m + 1)}
\quad \hbox{for} \ 0 \leq m \leq \ell \ , \label{6.48e}
\eea
\noi for all $t \in I$. Let $(q, \sigma )\in {\cal C}(I, X^{k,\ell})$ satisfy
\bea
\label{6.49e}
&&\parallel q\parallel_2 \ \leq Y_0\ t^{-\lambda_0} \quad , \quad 
\parallel \omega^k
q \parallel_2 \ \leq Y\ t^{-\lambda}\ , \\
&&\parallel \omega^m \nabla \sigma \parallel_2 \ \leq Z\ 
t^{-\lambda_0 + \beta (m + 1)} \quad
\hbox{for} \ 0 \leq m \leq \ell \ , \label{6.50e}
\eea
\noi for all $t \in I$. Then the system (\ref{2.31e}) has a (unique) 
solution $(q', \sigma ')
\in {\cal C} (I, X^{k,\ell})$ satisfying
\bea
\label{6.51e}
&&\parallel q' \parallel_2 \ \leq Y'_0\ t^{-\lambda_0} \quad , \quad 
\parallel \omega^k
q' \parallel_2 \ \leq Y'\ t^{-\lambda}\ ,  \\
&&\parallel \omega^m \nabla \sigma ' \parallel_2 \ \leq Z'\ 
t^{-\lambda_0 + \beta (m + 1)} \quad
\hbox{for} \ 0 \leq m \leq \ell \ , \label{6.52e}
\eea
\noi for some $Y'_0$, $Y'$, $Z'$ depending on $k$, $\ell$, $\beta$, 
$\lambda_0$, $\lambda$,
$a$, $b$, $c_0$, $c_1$, $c_2$, $Y_0$, $Y$, $Z$ and $T$, for all $t 
\in I$. That solution
satisfies the estimates (\ref{6.2e}) (\ref{6.3e}) (\ref{6.4e}) of 
Lemma 6.1 for all $t \in I$.
Two such solutions $(q'_i, \sigma '_i)$ associated with $(q_i, 
\sigma_i)$, $i = 1,2$, satisfy
the estimates (\ref{6.29e}) (\ref{6.30e}) of Lemma 6.2 for all $t \in 
I$. The solution $(q',
\sigma ')$ is actually unique in ${\cal C}(I, X^{k,\ell})$ under the 
condition that $(q',
\sigma ')$ tends to zero in $X^{0,0}$ norm when $t \to \infty$.}\\

\noi {\bf Proof.} The proof consists in showing that the solution 
$(q'_{t_0}, \sigma '_{t_0})$
of the linearized system (\ref{2.31e}) with $t_1 = \infty$ and with 
initial data $(q'_{t_0},
\sigma '_{t_0})(t_0) = 0$ for finite $t_0$, obtained from Proposition 
6.1, satisfies the
estimates (\ref{6.51e}) (\ref{6.52e}) uniformly in $t_0$ for $t \leq 
t_0$ (namely with $Y'_0$,
$Y'$ and $Z'$ independent of $t_0$), and that when $t_0 \to \infty$, 
that solution converges on
the compact subintervals of $I$ uniformly in the norms considered in 
Lemma 6.2. \par

We first derive the estimates (\ref{6.51e}) (\ref{6.52e}) for that 
solution, omitting the
subscript $t_0$ for brevity in that part of the proof. Let
\beq
\label{6.53e}
y'_0 = \parallel q'\parallel_2 \quad , \quad y' = \parallel \omega^k 
q'\parallel_2 \quad ,
\quad z'_m = \parallel \omega^m \nabla \sigma ' \parallel_2 \ . \eeq
\noi We first estimate $y'_0$. Substituting (\ref{6.47e}) 
(\ref{6.49e}) (\ref{6.50e}) into
(\ref{6.2e}), and omitting an overall constant, we obtain
\bea
\label{6.54e}
|\partial_t \ y'_0| &\leq& a\ Z\ t^{-2-\lambda_0+\beta} + a^2 \ Y_0 \ 
t^{-1-\beta -
\lambda_0}\nn \\
&&+ a \left ( Y_0^3 \bar{Y}\right )^{1/2} \ t^{-1 - 2 \lambda_0 + 
(\lambda_0 - \lambda )/2k} +
c_0 \ t^{-1-\lambda_0}
\eea
\noi where $\bar{Y} = Y \vee Y_0$. Integrating (\ref{6.54e}) from $t_0$
to $t$ with $y'_0(t_0) = 0$, using the fact that $\lambda_0 > 1$ and 
$\lambda_0 - (\lambda_0 -
\lambda )/2k > k - 1/2 + \lambda$, and defining
\beq
\label{6.55e}
Y'_0 = \ \parallel t^{\lambda_0}\ y'_0 ; L^{\infty}([T, t_0])\parallel \ ,
\eeq
\noi we obtain
\beq
\label{6.56e}
Y'_0  \leq c_0 + a\ Z\ T^{-(1- \beta )} + a^2 \ Y_0 \ T^{-\beta} + a 
\left ( Y_0^3
\bar{Y}\right )^{1/2} \ T^{-(k-1/2+ \lambda)}  \ . \eeq
\noi That estimate is manifestly uniform in $t_0$.\par

We next estimate $y'$, wasting part of the time decay in order to 
alleviate the computation.
In particular when estimating $s = S + \sigma$, we use the fact that 
the time decay of
$\sigma$ is better than that of $S$ by at least a power $1 - \eta$. 
Furthermore in the
contributions coming from $B_S$, we eliminate $Y_0$ and $\lambda_0$ 
by using $Y_0 \leq \bar{Y}
= Y \vee Y_0$ and $\lambda_0 > \lambda + k$. In particular we estimate
\bea
\label{6.57e}
\parallel \omega^m q\parallel_2 &\leq& \parallel q 
\parallel_2^{1-m/k} \ \parallel
\omega^kq\parallel_2^{m/k} \nn \\
&\leq& \bar{Y} \ t^{-\lambda_0(1-m/k)-\lambda m/k} \leq \bar{Y}\ 
t^{-\lambda + m -k}
\eea
\noi for $0 \leq m \leq k$, and similarly
\bea
\label{6.58e}
&&\parallel \omega^m q'\parallel_2 \ \leq y{'}^{m/k} \ \left ( Y'_0 \ 
t^{-\lambda_0}\right
)^{1-m/k} \nn \\
&&\leq t^{m-k} \ y'^{m/k}\left ( Y'_0 \ t^{-\lambda}\right 
)^{1-m/k}\leq t^{m-k} \left ( y' +
Y'_0 \ t^{-\lambda}\right ) \ .
  \eea
\noi Substituting (\ref{6.47e}) (\ref{6.49e}) (\ref{6.50e}) into 
(\ref{6.3e}), using
(\ref{6.57e}) (\ref{6.58e}) and omitting an overall constant, we obtain
\bea
\label{6.59e}
&&\left | \partial_t y' \right | \leq b_0 \Big \{ y{'}_0^{1/k} \ 
y{'}^{1-1/k} + y{'}_0^{1 -
\delta /k} \ y'^{\delta /k} \ t^{k-1- \delta /3} \Big \} \nn \\
&&+ a\ Z \ t^{-2-\lambda -k + \beta (k+1)} + t^{-1} \left ( 
bt^{-\eta} + Zt^{-1} \right ) \left
( y' + Y'_0\ t^{-\lambda}\right ) \nn \\
&&+ a^2 \ t^{-2} \left ( y' + Y'_0\ t^{-\lambda} + \bar{Y} \ 
t^{-\lambda} \right ) + a\
t^{-k-1/2-\lambda} \left ( y' + Y'_0\ t^{-\lambda} + \bar{Y} \ 
t^{-\lambda} \right ) \bar{Y}\nn
\\
&&+ t^{1-2k-2\lambda} \ \bar{Y}^2\left ( y' + Y'_0 \ 
t^{-\lambda}\right ) + c_1\ t^{-1-\lambda} \
. \eea
\noi The terms linear in $y'$ in the RHS of (\ref{6.59e}) can be 
eliminated by changing
variables from $y'$ to $y' \exp (-E(t))$, where
\bea
\label{6.60e}
E(t) = b \eta^{-1} \ t^{-\eta} + Z t^{-1} + a^2 t^{-1} &+& (k - 1/2 + 
\lambda )^{-1} a \bar{Y}
t^{-(k-1/2+ \lambda)} \nn \\
&+& (2(k-1 + \lambda ))^{-1} \ \bar{Y}^2 \ t^{-2(k-1 + \lambda )}
\eea
\noi so that it is sufficient to estimate $y'$ from (\ref{6.59e}) 
with those terms omitted and
to multiply the end result by $\exp (E(t))$. With those terms 
omitted, and with the help of
the estimate of $y'_0$, (\ref{6.59e}) can be rewritten as
\bea
\label{6.61e}
\left | \partial_t y' \right | &\leq& b_0 \Big \{ Y{'}_0^{1-\delta 
/k} \ y{'}^{\delta /k} \
t^{k-1-\delta /3-\lambda_0(1 - \delta /k)} + Y{'}_0^{1/k} \ 
y'^{1-1/k} \ t^{-\lambda_0/k}\Big
\}  \nn \\
&&+ t^{-1-\lambda} \ C_1(t)
\eea
\noi where
\bea
\label{6.62e}
&&C_1(t) = a\ Z\ t^{-(1-\beta )(k+1)} + \left ( b\ t^{-\eta} + 
Zt^{-1} \right ) Y'_0\nn \\
&&+ \left ( a^2 \ t^{-1} + a \bar{Y} \ t^{-(k-1/2+ \lambda)} \right ) 
\left ( Y'_0 + \bar{Y}
\right ) + \bar{Y}^2 \ Y'_0\ t^{-2(k-1+ \lambda )} + c_1 \ . \eea
\noi In particular $C_1(t)$ is decreasing in $t$. \par

The inequality (\ref{6.61e}) is essentially identical with 
(\ref{5.24e}), up to notational
change and replacement of $c_1$ by $C_1(t)$. Proceeding as in Section 
5, defining
\beq
\label{6.63e}
Y' = \ \parallel t^{\lambda} y' ; L^{\infty}([T, t_0])\parallel
\eeq
\noi and reintroducing the factor $\exp (E(t))$, we obtain (see (\ref{5.30e}))
\beq
\label{6.64e}
Y' \leq \exp (E(T)) \Big \{ \left ( 4\lambda^{-1} b_0 \right )^k \ 
Y'_0 + 2 \lambda^{-1} \
C_1(T) \Big \} \ ,
   \eeq
\noi an estimate which is again manifestly uniform in $t_0$. This 
completes the proof of
(\ref{6.51e}). \par

We next estimate $z'_m$ for $0 \leq m \leq \ell$. By interpolation, 
it suffices to estimate
$z'_0$ and $z'_{\ell}$. We define
\beq
\label{6.65e}
Z'_m = \ \parallel t^{\lambda_0 - \beta (m+1)} \ z'_m ; 
L^{\infty}([T, t_0])\parallel
\eeq
\noi and
\beq
\label{6.66e}
Z' = \mathrel{\mathop{\rm Sup}_{0 \leq m \leq \ell}} Z'_m = Z'_0 \vee 
Z'_{\ell} \ .
\eeq
\noi We first estimate $z'_0$. Substituting (\ref{6.48e}) 
(\ref{6.49e}) (\ref{6.50e}) into
(6.4)$_0$ and omitting an overall constant, we obtain
\bea
\label{6.67e}
&&\left | \partial_t \ z'_0 \right | \leq \left ( b \ t^{-1-\eta} + 
Zt^{-2} \right ) z'_0 +
bZt^{-1-\eta - \lambda_0 + \beta} \nn \\
&&+ a\ Y_0\ t^{-1-\lambda_0 + \beta} + Y_0^2 \ t^{-1-2\lambda_0 + 5 
\beta /2} + c_2\
t^{-1-\lambda_0 + \beta} \ . \eea
\noi Integrating (\ref{6.67e}) from $t$ to $t_0$, we obtain
\bea
\label{6.68e}
&&Z'_0 \leq \exp \left ( b \ \eta^{-1} \ T^{-\eta} + Z \ T^{-1} 
\right ) \Big \{ \left ( \eta +
\lambda_0 - \beta \right )^{-1} \ b \ Z\ t^{-\eta} \nn \\
&&\qquad \qquad + ( \lambda_0 - \beta )^{-1} \left ( a\ Y_0 + Y_0^2\ 
T^{-\beta} + c_2
\right ) \Big \}  \eea
\noi where we have used again the fact that $\lambda_0 > 5 \beta /2$. \par

We next estimate $z'_{\ell}$. Using the inequality
\bea
\label{6.69e}
\parallel \nabla \sigma ' \parallel_{\infty} &\leq& C \parallel 
\omega^{\ell} \nabla \sigma
'\parallel_2^{3/2\ell} \ \parallel \nabla \sigma '\parallel_2^{1 - 
3/2 \ell} \nn \\
&\leq& C\ t^{3\beta /2} \left ( t^{-\beta \ell} \parallel 
\omega^{\ell} \nabla \sigma '
\parallel_2 \ + \ \parallel \nabla \sigma ' \parallel_2 \right ) \ ,
\eea
\noi substituting (\ref{6.48e}) (\ref{6.49e}) (\ref{6.50e}) into 
(\ref{6.4e}) with $m = \ell$,
and omitting again an overall constant, we obtain
\bea
\label{6.70e}
&&\left | \partial_t \ z'_{\ell} \right | \leq \left ( b \ 
t^{-1-\eta} + Zt^{-2} \right )
\left ( z'_{\ell} + Z'_0\ t^{-\lambda_0 + \beta (\ell + 1)} \right ) \nn \\
&&+ b\ Z\ t^{-1-\eta - \lambda_0 + \beta (\ell + 1)} + a\ Y_0 \ 
t^{-1-\lambda_0 + \beta (\ell +
1)} \nn \\
&&+ Y_0^2\ t^{-1-2\lambda_0 + \beta (\ell + 5/2)} + c_2\ 
t^{-1-\lambda_0 + \beta (\ell + 1)} \ .
\eea
\noi Integrating (\ref{6.70e}) as before, we obtain
\bea
\label{6.71e}
&&Z'_{\ell} \leq \exp \left ( b \ \eta^{-1} \ T^{-\eta} + Z \ T^{-1} 
\right ) \Big \{ b (Z +
Z'_0) \eta^{-1}\ T^{-\eta} + Z Z'_0\ T^{-1} \nn \\
&&\qquad \qquad + \nu^{-1} \left ( a\ Y_0 + Y_0^2\ T^{-\beta} + c_2 
\right ) \Big \}
\eea
\noi where $\nu = \lambda_0 - \beta (\ell + 1) > 0$, which together 
with (\ref{6.68e})
completes the proof of (\ref{6.52e}). \par

We have proved that the solution $(q'_{t_0}, \sigma '_{t_0})$ of the 
system (\ref{2.31e}),
vanishing at $t_0$, satisfies the estimates (\ref{6.51e}) 
(\ref{6.52e}) for $t \in [T, t_0]$,
with $Y'_0$, $Y'$, $Z'$ satisfying (\ref{6.56e}) (\ref{6.64e}) 
(\ref{6.66e}) (\ref{6.68e})
(\ref{6.71e}), which are uniform in $t_0$. We now prove that 
$(q'_{t_0}, \sigma '_{t_0})$
tends to a limit when $t_0 \to \infty$. For that purpose, we first 
let $(q'_i, \sigma
'_i)$, $i = 1,2$, be two solutions of the system (\ref{2.31e}) 
corresponding to the same $(q,
\sigma )$ and defined in an interval $[T, t_0)$ for some $t_0 > T$. 
Let $(q'_-, \sigma '_-) =
(1/2)  (q'_1 - q'_2, \sigma '_1 - \sigma '_2)$. $L^2$ norm 
conservation for $q'_i$ and
therefore for $q'_-$ implies
\beq
\label{6.72e}
\parallel q'_-(t)\parallel_2 \ = \ \parallel q'_-(t_0)\parallel_2 
\qquad \hbox{for all} \ t \in
[T, t_0] \ . \eeq
\noi Furthermore, the simple case $q_- = 0$, $\sigma_- = 0$ of 
(6.30)$_0$ implies
\beq
\label{6.73e}
\partial_t \parallel \nabla \sigma '_- \parallel_2 \ \leq \ C \ 
t^{-1} \left ( b\ t^{-\eta} +
Z\ t^{-1} \right ) \parallel \nabla \sigma'_-\parallel_2 \eeq
\noi and therefore
\beq
\label{6.74e}
\parallel \nabla \sigma '_-(t) \parallel_2 \ \leq \exp \left ( C 
\left ( \eta^{-1}\ b\
t^{-\eta} + Z\ t^{-1} \right ) \right ) \parallel \nabla \sigma '_- 
(t_0)\parallel_2 \ .
\eeq
\noi Let now $T < t'_1 < t'_2 < \infty$ and let $(q'_i, \sigma '_i) = 
(q'_{t'_i}, \sigma
'_{t'_i})$, $i = 1,2$. Then
\beq
\label{6.75e}
\parallel q'_-(t) \parallel_2 \ = \ \parallel q'_-(t'_1)\parallel_2 \ 
= (1/2) \parallel
q'_2(t'_1)\parallel_2\ \leq Y'_0\ t_1{'}^{-\lambda_0} \quad \hbox{for 
all}\ t < t'_1
\eeq
\noi by (\ref{6.72e}) with $t_0 = t'_1$ and (\ref{6.51e}) for $q' = 
q'_2$ and $t = t'_1$
Similarly, by (\ref{6.52e}) with $\sigma ' = \sigma '_2$ and $t = t'_1$,
$$\parallel \nabla \sigma '_-(t'_1)\parallel\ = (1/2) \parallel 
\nabla \sigma '_2 (t'_1)\parallel
\ \leq Z'_0\ t{'}_1^{-\lambda_0 + \beta}$$
\noi so that by (\ref{6.74e}) with $t_0 = t'_1$,
\beq
\label{6.76e}
\parallel \nabla \sigma '_-(t)\parallel\ \leq \ \exp \left ( C \left 
( \eta^{-1}\ b\ t^{-\eta}
+ Z\ t^{-1} \right ) \right )  Z'_0\ t{'}_1^{-\lambda_0 + \beta} 
\quad \hbox{for all}\ t < t'_1 \
. \eeq
\indent From (\ref{6.75e}) (\ref{6.76e}), it follows that $(q'_{t_0}, 
\sigma '_{t_0})$
converges to a limit $(q', \sigma ') \in {\cal C}(I, X^{0,0})$ 
uniformly on the compact
subintervals of $I$. From the uniform estimates (\ref{6.51e}) 
(\ref{6.52e}) and from Lemma 6.1,
it then follows by a standard compactness argument that $(q', \sigma 
') \in {\cal C}(I,
X^{k,\ell})$ and that $(q', \sigma ')$ also satisfies the estimates 
(\ref{6.51e})
(\ref{6.52e}). Clearly $(q', \sigma ')$ satisfies the system 
(\ref{2.31e}). This completes
the existence part of the proof. \par

The uniqueness statement follows immediately from (\ref{6.72e}) 
(\ref{6.74e}) by letting $t_0
\to \infty$. \par \nobreak \hfill
$\sq$\par

We now turn to the main result of this section, namely the fact that 
for $T$ sufficiently large
(depending on $(W, S)$), the auxiliary system in difference form 
(\ref{2.30e}) has a solution
$(q, \sigma )$ defined for all $t \geq T$ and decaying at infinity in 
a suitable sense. In the
same spirit as for Proposition 4.4, this will be done by showing that 
the map $\Gamma : (q,
\sigma ) \to (q', \sigma ')$ defined by Proposition 6.2 is a 
contraction in suitable
circumstances. According to our intuition of scattering, another 
natural route towards the
same result would be to construct first the solution $(q_{t_0}, 
\sigma_{t_0})$ of the
auxiliary system (\ref{2.30e}) vanishing at $t_0$ and to take the 
limit of that solution as
$t_0 \to \infty$. That route can also be followed, but it is slightly 
more complicated than
the previous one. One of the complications comes from the fact that 
the system (\ref{2.30e})
depends on $t_1$. In view of Warning 4.2, for finite $t_0$, we expect 
difficulties if we
take $t_1 >t_0$. This prompts us to take $t_1 = t_0$. The comparison 
of two solutions
$(q_{t_0}, \sigma_{t_0})$ corresponding to different values of $t_0$ 
is then complicated
by the fact that they do not solve exactly the same system, so that 
Lemma 6.2 is not
directly applicable and additional terms occur in the comparison. On 
the other hand, for
$B_0 \not= 0$, the construction of the solution $(q_{t_0}, 
\sigma_{t_0})$ of (\ref{2.30e})
is expected to meet difficulties for $t \geq t_0$ because of Warning 
4.1. We shall
therefore undertake it for $t \leq t_0$ only, which is sufficient 
anyway to take the limit
$t_0 \to \infty$. That construction proceeds again by a contraction 
starting from the
solutions obtained for the linearized system. The corresponding proof 
for $t_0 < \infty$
is not significantly simpler than for $t_0 = \infty$, which is 
another reason why the
second method is more complicated than the first one, since in addition to that
construction, a limiting procedure is needed. \par

We now state the main result and formalize the previous heuristic 
discussion in the
following proposition. \\

\noi {\bf Proposition 6.3.} {\it Let $1 < k \leq \ell$ and $\ell > 
3/2$. Let $\beta$,
$\lambda_0$ and $\lambda$ satisfy (\ref{6.44e}) and in addition
$$1 + \lambda > \beta (5/2 - k) \ .$$
\noi Let $B_0$ satisfy the estimates (\ref{3.17e}) for $0 \leq m \leq 
k$. Let $(W,S)$ satisfy
the assumptions of Proposition 6.2 in $[1, \infty )$. Then there 
exists $T$, $1 \leq T <
\infty$, and positive constants $Y_0$, $Y$ and $Z$, depending on $k$, 
$\ell$, $\beta$,
$\lambda_0$, $\lambda$, $a$, $b$, $c_0$, $c_1$ and $c_2$, such that 
the following holds. \par

(1) For all $t_0$, $T \leq t_0 < \infty$, the system (\ref{2.30e}) 
with $t_1 = t_0$ has a
unique solution $(q, \sigma ) \in {\cal C} (I, X^{k,\ell})$ with $I = 
[T, t_0]$ and $(q,
\sigma)(t_0) = 0$. That solution satisfies the estimates 
(\ref{6.49e}) (\ref{6.50e}) for all
$t \in I$. \par

(2) The system (\ref{2.30e}) with $t_1 = \infty$ has a unique 
solution $(q, \sigma ) \in
{\cal C}(I, X^{k,\ell})$, where $I = [T, \infty )$ satisfying the 
estimates (\ref{6.49e})
(\ref{6.50e}) for all $t \in I$. \par

(3) Let $(q_{t_0}, \sigma_{t_0})$ be the solution defined in Part (1) 
for $t_0 < \infty$ and
let $(q, \sigma )$ be the solution defined in Part (2) for $t_0 = 
\infty$. When $t_0 \to
\infty$, $(q_{t_0}, \sigma_{t_0})$ converges to $(q, \sigma )$ 
strongly in $L^{\infty}(J,
X^{k',\ell '})$ for $0 \leq k' < k$, $0 \leq \ell ' < \ell$, and in 
the weak-$*$ sense in
$L^{\infty}(J, X^{k,\ell })$ for any interval $J = [T, \bar{T}]$ with 
$\bar{T} < \infty$.}\\

\noi {\bf Proof.} \underbar{Parts (1) and (2)}. We prove Parts (1) 
and (2) together, because
the proof is exactly the same for both. It consists in showing that 
the map $\Gamma : (q,
\sigma ) \to (q', \sigma ')$ defined by solving the linearized system 
(\ref{2.31e}) is a
contraction on a suitable subset of ${\cal C}(I, X^{k,\ell})$ in the 
lower norms used in Lemma
6.2. For $t_0 < \infty$, the map $\Gamma$ is defined by Proposition 
6.1, restricted to those
$(q, \sigma )$ satisfying $(q, \sigma )(t_0) = 0$, with the initial 
data $(q', \sigma ') (t_0)
= 0$. For $t_0 = \infty$, the map $\Gamma$ is defined by Proposition 
6.2. The relevant
estimates on $\Gamma$ are those derived in the proof of Proposition 
6.2 in the case $t_1 =
\infty$. The same estimates also apply to the case $t_0 = t_1 < 
\infty$, which is relevant for
Part (1) of this proposition. They are independent of $t_0$. \par

We define the set
\bea
\label{6.77e}
&&{\cal R} = \Big \{ (q, \sigma ) \in {\cal C} (I, X^{k,\ell}) : (q, 
\sigma ) (t_0) = 0 \quad
\hbox{if} \ t_0 < \infty\ , \nn \\
&&\parallel t^{\lambda_0}q;L^{\infty}(I, L^2)\parallel\ \leq Y_0 \ , 
\ \parallel t^{\lambda}
\omega^k q; L^{\infty} (I, L^2)\parallel \ \leq Y \ , \nn \\
&&\mathrel{\mathop {\rm Sup}_{0 \leq
m \leq \ell}} \parallel t^{\lambda_0- \beta (m+1)} \omega^m \nabla 
\sigma ; L^{\infty}(I,
L^2)\parallel \ \leq Z \Big \} \ .\eea
\indent We first show that ${\cal R}$ is stable under $\Gamma$ for 
suitable $Y_0$, $Y$, $Z$ and
for sufficiently large $T$. Let $(q, \sigma) \in {\cal R}$ and $(q' , 
\sigma ') = \Gamma (q,
\sigma )$. Then $(q', \sigma ')$ satisfies the estimates 
(\ref{6.56e}) (\ref{6.64e})
(\ref{6.68e}) (\ref{6.71e}) where $Y'_0$, $Y'$, $Z'_m$ are defined by 
(\ref{6.53e})
(\ref{6.55e}) (\ref{6.63e}) (\ref{6.65e}) and/or their extension to 
$t_0 < \infty$. It is
therefore sufficient to ensure that the RHS of (\ref{6.56e}) (\ref{6.64e})
(\ref{6.68e}) (\ref{6.71e}) are not larger than $Y_0$, $Y$, $Z$, and 
$Z$ respectively. For
that purpose, it is sufficient to choose
\beq
\label{6.78e}
\left \{ \begin{array}{l} Y_0 = 2c_0 \quad , Z = 2e\nu^{-1} (c_2 + 4ac_0) \\
\\
Y = e \Big \{ (4\lambda^{-1} b_0)^k 2c_0 + 4 \lambda^{-1}(c_1 + c_0) \Big \}
\end{array} \right .
\eeq
\noi and to take $T$ sufficiently large in the sense that
\beq
\label{6.79e}
\left \{ \begin{array}{l} a\ Z\ T^{-(1- \beta )} + a^2 \ Y_0 \ 
T^{-\beta} + a \left ( Y_0^3 \
\bar{Y}\right )^{1/2} T^{-(k-1/2+ \lambda )} \leq c_0 \ ,\\ \\
\eta\ T^{\eta} \geq 8 b \quad , \quad T \geq 4eZ \quad , \quad a\ 
T^{\beta} \geq Y_0 \ , \\
\\
E(T) \leq 1 \quad , \quad C_1(T) \leq 2 (c_1 + c_0) \ .
\end{array} \right .
\eeq
\noi The conditions (\ref{6.79e}) are lower bounds on $T$ expressed 
in terms of the parameters
listed in the Proposition, after substitution of (\ref{6.78e}). \par

We next show that $\Gamma$ is a contraction on ${\cal R}$ in the 
norms considered in Lemma
6.2. Let $(q_i, \sigma_i) \in {\cal R}$ and $(q'_i, \sigma '_i) = 
\Gamma (q_i, \sigma_i)$, $i
= 1,2$, and define $(q_{\pm} , \sigma _{\pm})$ and $(q'_{\pm} , 
\sigma '_{\pm})$ as in Lemma
6.2. We define in addition
\bea
\label{6.80e}
&&y_- = \ \parallel q_-\parallel_2 \quad , \quad z_{-m} = \parallel 
\omega^m \nabla
\sigma_-\parallel_2\\
&&Y_- =\ \parallel t^{\lambda_0} y_- ; L^{\infty}(I) \parallel \quad 
, \quad Z_- =
\mathrel{\mathop {\rm Sup}_{0 \leq m \leq \ell_0}}\parallel 
t^{\lambda_0 - \beta (m+1)} z_{-m};
L^{\infty}(I)\parallel \label{6.81e} \eea
\noi and we make similar definitions for the primed quantities. We 
take $\ell_0 = [3/2 - k]_+$
and estimate $y'_-$ and $z'_{-m}$ by (\ref{6.29e}) (\ref{6.30e}), 
taking advantage of the fact
that $m' = m$ in (\ref{6.30e}) for that choice of $\ell_0$. Using the 
fact that $\Gamma$ maps
${\cal R}$ into itself and omitting again overall constants, we obtain
\bea
\label{6.82e}
&&\left | \partial_t y'_-\right | \leq a\ Z_-\ t^{-2-\lambda_0 + 
\beta} + \bar{Y}\  Z_-\
t^{-2-\lambda_0 + \beta (\ell_0 + 1)- \lambda } + a^2 \ Y_-\ t^{-1 - 
\beta - \lambda_0} \nn \\
&&+ a\ \bar{Y}Y_-\ t^{-1-2\lambda_0 + (\lambda_0 -
\lambda )/2k} + \bar{Y}^2\ Y_- \ t^{-1-3\lambda_0 + 2(\lambda_0 - \lambda )/k}
\eea
\noi where $\bar{Y} = Y \vee Y_0$,
\bea
\label{6.83e}
&&\left | \partial_t z'_{-m}\right | \leq \left ( b\ t^{-1-\eta} + Z\
t^{-2} \right ) \left ( z'_{-m} + Z_-\ t^{-\lambda_0 + \beta (m+1)} 
\right ) \nn \\
&&+ a \ Y_-\ t^{-1 - \lambda_0 + \beta (m+1)} + Y_0 Y_-\ 
t^{-1-2\lambda_0 + \beta
(m+5/2)}
\eea
\noi for $0 \leq m \leq \ell_0$. Integrating (\ref{6.82e}) 
(\ref{6.83e}) from $t$ to $t_0$ with
$(y'_-, z'_{-m})(t_0) = 0$ and using again the fact that $\lambda_0 > 
\lambda + k > 1$ and
$\lambda_0 > \beta (\ell + 1) > \beta ((\ell_0 + 2) \vee 5/2)$, we obtain
\bea
\label{6.84e}
&&Y'_- \leq a\ Z_-\ T^{-(1-\beta)} + \bar{Y} \  Z_-\ T^{-1-\lambda + 
\beta (\ell_0 + 1)} \nn \\
&&+ a^2 \ Y_-\ T^{-\beta} + a\ \bar{Y}\ Y_-\ T^{-(k-1/2+ \lambda)} + 
\bar{Y}^2 \ Y_-\ T^{-2(k-1
+ \lambda)}\ ,
\eea
\bea
\label{6.85e}
&&Z'_- \leq \exp \left ( b\ \eta^{-1}\ T^{-\eta} + Z\ T^{-1}\right ) 
\Big \{ \left ( b\
\eta^{-1}\ T^{-\eta} + Z\ T^{-1} \right ) Z_- \nn \\
&&+ \beta^{-1} \left (  a\ Y_- + T^{-\beta } Y_0 \ Y_- \right ) \Big \}\ .
\eea
\indent We now ensure that the map $\Gamma$ is a contraction for the 
norms defined by
(\ref{6.80e}) (\ref{6.81e}) in the form
\beq
\label{6.86e}
\left \{ \begin{array}{l} Y'_- \leq  \left ( c^{-1} \ Z_- + Y_- 
\right )/4 \\  \\
Z'_- \leq \left ( Z_- + c\ Y_- \right )/4  \end{array} \right .
\eeq
\noi which imply
\beq
\label{6.87e}
Z'_- + c\ Y'_- \leq \left ( Z_- + c \ Y_- \right )/2
\eeq
\noi by taking $c = 8 \beta^{-1}a$ and $T$ sufficiently large, 
depending on the parameters
listed in the proposition, in part explicitly and in part through 
$Y_0$, $Y$ and $Z$ defined
by (\ref{6.78e}). (It is only at this point that we need the 
condition $1 + \lambda > \beta
(5/2 - k)$, in order to ensure that the power of $T$ in the second 
term in the RHS of
(\ref{6.84e}) is negative). \par

We have proved that for sufficiently large $T$, the map $\Gamma$ maps 
${\cal R}$ defined by
(\ref{6.77e}) into itself and is a contraction for the norms 
(\ref{6.81e}). By a standard
compactness argument, ${\cal R}$ is closed for the latter norms, and 
therefore $\Gamma$ has a
unique fixed point in ${\cal R}$, which completes the existence part 
of the proof of Parts (1)
and (2). \par

The uniqueness statement of Part (1) is a special case of Proposition 
4.2 part (1), while the
uniqueness statement of Part (2) follows from Proposition 4.2 part 
(3) and from the fact that
$\lambda_0 > 1 > \beta_2$. \par

\noi \underbar{Part (3)}. Let $T < t'_1 < t'_2 < \infty$ and let 
$(q_i, \sigma_i)$, $i = 1,2$,
be the solutions of the ssytem (\ref{2.30e}) obtained in part (1) and 
corresponding to $t_0 =
t_1 = t'_i$ respectively. Those solutions satisfy the estimates 
(\ref{6.49e}) (\ref{6.50e})
for $t \leq t'_i$. Define as before $(q_{\pm}, \sigma_{\pm}) = (1/2) 
(q_1 \pm q_2, \sigma_1
\pm \sigma_2)$. We shall estimate $(q_-, \sigma_-)$ for $t \leq t'_1$ 
in the norms considered
in Lemma 6.2. In order to alleviate the notation, we omit the prime 
on $t_1$, $t_2$ in the
rest of the proof. By (\ref{6.49e}) (\ref{6.50e}), we estimate
  \beq
\label{6.88e}
\left \{ \begin{array}{l} \parallel q_-(t_1)\parallel_2 \ = (1/2) 
\parallel q_2(t_1)\parallel_2
\ \leq (1/2) Y_0\ t_1^{-\lambda_0} \\  \\ \parallel \omega^m \nabla 
\sigma_-(t_1)\parallel_2 \ =
(1/2) \parallel \omega^m \nabla \sigma_2 (t_1)\parallel_2 \ \leq (1/2)
Zt_1^{-\lambda_0 + \beta (m+1)}\end{array} \right .
\eeq
\noi for $0 \leq m \leq \ell$. On the other hand $(q_-, \sigma_-)$ 
satisfies a system closely
related to (\ref{6.32e}) where however $(q'_{\pm} , \sigma '_{\pm}) = 
(q_{\pm} , \sigma
_{\pm})$ and where additional terms appear because of the different 
values $t_1$ and $t_2$
occuring in $B_S$ and $B_L$. More precisely
\bea
\label{6.89e}
&&\partial_t q_- = i(2t^2)^{-1} \Delta q_- + t^{-2} \Big \{ Q(s_+, 
q_-) + Q(\sigma_-, w_+)\Big
\} + it^{-1} B_0 q_- \nn \\
&&+ it^{-1} \Big \{ \left ( B_S^{t_1, \infty}(w_+, w_+) + 
B_S^{t_1}(q_-, q_-)\right ) q_- + 2
B_S^{t_1} (w_+, q_-)w_+ \Big \} \nn \\
&&- i(2t)^{-1} \left ( B_S^{t_2} - B_S^{t_1}\right ) \left ( q_2, q_2 
+ 2W \right ) \left ( q_2 +
W \right )
\eea
\begin{eqnarray*}
\partial_t \sigma_- = t^{-2} \left ( s_+ \cdot \nabla \sigma_- + 
\sigma_- \cdot \nabla s_+
\right ) - 2 t^{-1} \nabla B_L^{t_1} (w_+ ,q_-) \\
+ t^{-1} \nabla \left ( B_L^{t_2} - B_L^{t_1} \right ) \left ( q_2 , 
q_2 + 2W \right ) \ .
\end{eqnarray*}
\noi We first estimate the additional terms in (\ref{6.89e}) as 
compared with (\ref{6.32e}).
 From (\ref{6.36e})-(\ref{6.39e}) we obtain
\bea
\label{6.90e}
&&\parallel \left ( B_S^{t_2} - B_S^{t_1}\right ) \left ( q_2, q_2 + 
2W \right ) \left (q_2 + W
\right ) \parallel_2 \nn \\
&&\leq C \Big \{ \left ( at^{-\beta} + \ \parallel q_2\parallel_3 
\right ) a\ I_0 \left (
\parallel q_2\parallel_2\right ) + a\ I_{-1} \left ( \parallel 
q_2\parallel_2 \ \parallel q_2
\parallel_3 \right ) \nn \\
&&+ \parallel q_2 \parallel_6 \ I_{-1/2} \left ( \parallel 
q_2\parallel_2\ \parallel q_2
\parallel_6 \right ) \Big \}
\eea
\noi where the various $I_m$'s are taken in the interval
$[t_1, t_2]$. From (\ref{6.49e}) and Sobolev inequalities, we obtain
\beq
\label{6.91e}
I_0 \left ( \parallel q_2\parallel_2 \right ) \leq Y_0 \ t^{1/2} 
\int_{t_1}^{t_2} dt' \
t'^{-3/2-\lambda_0} \leq Y_0 \ t^{1/2} \ t_1^{-1/2 - \lambda_0}
\eeq
\bea
\label{6.92e}
I_{-1} \left ( \parallel q_2\parallel_2 \ \parallel q_2\parallel_3 
\right ) &\leq& C\ Y_0 \
\bar{Y}\ t^{-1/2} \int_{t_1}^{t_2} dt' \ 
t'^{-1/2-2\lambda_0+(\lambda_0 - \lambda )/2k} \nn
\\
&\leq& C\ Y_0 \
\bar{Y}\ t^{-1/2} \ t_1^{1/2-2\lambda_0+(\lambda_0 - \lambda )/2k}
\eea
\bea
\label{6.93e}
&&\parallel q_2\parallel_6 I_{-1/2} \left ( \parallel q_2\parallel_2 
\ \parallel q_2\parallel_6
\right ) \leq C\ Y_0 \ \bar{Y}^2\ t^{-\lambda_0+(\lambda_0 - \lambda 
)/k} \nn \\
&&\times \int_{t_1}^{t_2} dt' \ t'^{-1-2\lambda_0+(\lambda_0 - 
\lambda )/k}\leq C\ Y_0 \
\bar{Y}^2\ t^{-\lambda_0+(\lambda_0 - \lambda )/k}\ 
t_1^{-2\lambda_0+(\lambda_0 - \lambda )/k}
\eea
\noi and therefore for $t \leq t_1$
\bea
\label{6.94e}
&&\parallel \left ( B_S^{t_2} - B_S^{t_1}\right ) \left ( q_2, q_2 + 
2W \right ) \left (q_2 + W
\right ) \parallel_2 \ \leq C\ t^{-\lambda_0 + 1/2} \Big \{ a^2 \ 
Y_0\ t_1^{-\beta -1/2}\nn \\
&&+ a\ Y_0 \ \bar{Y} \ t_1^{-\lambda_0 + (\lambda_0 - \lambda 
)/2k-1/2} + Y_0 \ \bar{Y}^2 \
t_1^{-2\lambda_0 + 2(\lambda_0 - \lambda )/k-1/2}\Big \} \ .  \eea
\noi Similarly using (\ref{6.43e}) and
$$I_{-3/2} \left ( \parallel q_2\parallel_2^2 \right ) \leq Y_0^2 \ 
t^{-1} \int_{t_1}^{t_2} dt' \
t'^{-2\lambda_0} \leq Y_0^2 \ t^{-1} \ t_1^{1-2\lambda_0}
$$
\noi we estimate for $t \leq t_1$
\beq
\label{6.95e}
\parallel \omega^{m+2} \left ( B_L^{t_2} - B_L^{t_1}\right ) \left ( 
q_2, q_2 + 2W \right
) \parallel_2 \ \leq C\ \left ( a\ Y_0 + Y_0^2 t^{3\beta /2 - 1} \ 
t_1^{1- \lambda_0} \right
)  t^{\beta (m+1)} \ t_1^{-\lambda_0} \ .  \eeq  \noi We define $y_-$ 
and $z_{-m}$ by
(\ref{6.80e}), we take again $\ell_0 = [3/2 - k]_+$, we choose 
$\lambda '_0$ satisfying
\beq
\label{6.96e}
1 \vee \left ( \lambda_0 - 1/2\right ) \vee \beta \left ( \ell_0 + 1 
\right ) < \lambda '_0 <
\lambda_0 \ ,
\eeq
\noi we define (see (\ref{6.81e}))
\beq
\label{6.97e}
Y_- = \ \parallel t^{\lambda '_0} y_-;L^{\infty}([T, t_1])\parallel \ 
, \ Z_- = \mathrel{\mathop
{\rm Sup}_{0 \leq m \leq \ell_0}}\parallel t^{\lambda '_0 - \beta (m 
+ 1)} z_{-m};L^{\infty}([T,
t_1])\parallel \eeq
\noi and we estimate those quantities in the same way as in the proof 
of Parts (1) and (2).
 From (\ref{6.89e}), we obtain differential inequalities for $y_-$, 
$z_{-m}$, very similar to
(\ref{6.82e}) (\ref{6.83e}) with $y'_- = y_-$, $z'_{-m} = z_{-m}$ 
with however additional terms
estimated by (\ref{6.94e}) (\ref{6.95e}). \par

We integrate those inequalities from $t$ to $t_1$, with initial 
condition at $t_1$ estimated by
(\ref{6.88e}). We then substitute the result in (\ref{6.97e}) and 
omitting an overall
constant, we obtain finally (see (\ref{6.84e}) (\ref{6.85e}))
\bea
\label{6.98e}
&&Y_- \leq a\ Z_-\ T^{-(1-\beta)} + \bar{Y} \ Z_-\  T^{-1-\lambda + 
\beta (\ell_0 + 1)} + a^2\
Y_-\ T^{-\beta} \nn \\
&&+ a\ \bar{Y}\ Y_-\ T^{-(k-1/2+ \lambda)} +\bar{Y}^2 \ Y_-\ 
T^{-2(k-1 + \lambda)} \nn \\
&&+ \Big \{ Y_0 + a^2\ Y_0 \ t_1^{-\beta} + a \ Y_0\ \bar{Y} \ 
t_1^{-(k-1/2 + \lambda )} +
Y_0\ \bar{Y}^2 \ t_1^{-2(k-1+\lambda)} \Big \} t_1^{-(\lambda_0 - 
\lambda '_0)}  \eea
\bea
\label{6.99e}
&&Z_- \leq \exp \left ( b\ \eta^{-1}\ T^{-\eta} + Z\ T^{-1}\right ) 
\Big \{ \left ( b\
\eta^{-1}\ T^{-\eta} + Z\ T^{-1} \right ) Z_- \nn \\
&&+ \beta^{-1} \left (  a\ Y_- + T^{-\beta } Y_0 \ Y_- \right ) + 
\left ( Z + a\ Y_0 + Y_0^2
\ t_1^{-\beta} \right ) t_1^{-(\lambda_0 - \lambda '_0)} \Big \}\ .  \eea
\noi Proceeding as above, we deduce therefrom that for $T$ 
sufficiently large and for a
suitable constant $c$
\beq
\label{6.100e}
Y_- + c Z_- \leq O\left ( t_1^{-(\lambda_0 - \lambda '_0)} \right ) \ .
\eeq
\noi From (\ref{6.100e}) it follows that $(q_{t_0}, \sigma_{t_0})$ 
tends to a limit uniformly
in compact subintervals of $[T, \infty )$ in the norms (\ref{6.80e}). 
By a standard compactness
argument, that limit belongs to ${\cal C} ([T, \infty ), X^{k,\ell})$ 
and satisfies
(\ref{6.49e}) (\ref{6.50e}). One sees easily that the limit satisfies 
the system (\ref{2.30e})
with $t_1 = \infty$, and therefore coincides with the solution 
obtained in Part (2). Actually,
as mentioned before, Part (3) provides an alternative (more 
complicated ) proof of Part (2).\par
\nobreak \hfill $\sq$\par

\mysection{Choice of (W, S) and remainder estimates}
\hspace*{\parindent} In this section, we construct approximate 
solutions $(W,S)$ of the system
(\ref{2.20e}) satisfying the assumptions needed for Propositions 6.2 
and 6.3 and in particular
the remainder estimates (\ref{6.47e}) (\ref{6.48e}), thereby allowing 
for the applicability of
Proposition 6.3, namely for the construction of solutions of the 
system (\ref{2.30e}). More
general $(W,S)$ also suitable for the same purpose, could also be 
constructed by exploiting the
gauge invariance of the system (\ref{2.20e}). \par

We rewrite the remainders as
\beq
\label{7.1e}
R_1(W,S) = U^*(1/t) \partial_t (U(1/t)W) - t^{-2} Q(S,W) - it^{-1} 
(B_0 + B_S (W, W)) W
\eeq
$$R_2(W,S) = \partial_t S - t^{-2} S \cdot \nabla S + t^{-1} \nabla 
B_L(W, W) \ .
\eqno(2.29)\equiv(7.2)$$
\noi We recall that $t_1 = \infty$ in $R_1$, $R_2$, and we omit $t_1$ 
from the notation. We
construct $(W,S)$ by solving the system (\ref{2.20e}) approximately 
by iteration. The $n$-th
iteration should be sufficient to cover the case $\lambda_0 < n$. 
Here we need $\lambda_0 > k >
1$, and we must therefore use at least the second iteration, which 
will allow for $\lambda_0 <
2$. For simplicity, we shall not go any further here. Accordingly we 
take $$W = w_0 + w_1
\qquad , \quad S = s_0 + s_1 \eqno(7.3)$$ \noi where
$$\left \{ \begin{array}{ll} \partial_t \ U(1/t)w_0 = 0 &\qquad w_0 
(\infty ) = w_+ \ ,\\ & \\
\partial_t \ s_0 = - t^{-1} \nabla B_L(w_0, w_0) &\qquad s_0 (1) = 0 
\ , \end{array} \right .
\eqno(7.4)$$
\noi so that
$$\left \{ \begin{array}{l} w_0 = U^*(1/t) w_+ \ ,\\  \\
s_0 (t) = - \displaystyle{\int_1^t} dt' \ t'^{-1} \nabla B_L(w_0(t'), 
w_0(t'))  \end{array}
\right . \eqno(7.5)$$
\noi and
$$\left \{ \begin{array}{ll} \partial_t ( U(1/t)w_1) = t^{-2}\ U(1/t) 
\ Q(s_0, w_0) &\qquad w_1
(\infty ) = 0 \ ,\\ & \\ \partial_t \ s_1 =  t^{-2} s_0 \cdot \nabla 
s_0 - 2t^{-1} \nabla
B_L(w_0, w_1) &\qquad s_1 (\infty) = 0 \ , \end{array} \right . \eqno(7.6)$$
\noi so that
$$\left \{ \begin{array}{l} w_1(t) = - U^*(1/t) 
\displaystyle{\int_t^{\infty}} dt'\ t'^{-2} \
U(1/t')\ Q(s_0(t'),w_0(t')) \\  \\ s_1 (t) = - 
\displaystyle{\int_t^{\infty}} dt' \ t'^{-2}
s_0(t') \cdot \nabla s_0(t') + 2 \displaystyle{\int_t^{\infty}} dt' \ 
t'^{-1} \nabla
B_L(w_0(t'), w_1(t'))  \ .\end{array} \right . \eqno(7.7)$$
\noi The remainders then become
$$\left \{ \begin{array}{l} R_1(W,S) = - t^{-2} \Big \{ Q(s_0, w_1) + 
Q(s_1 , w_0) + Q(s_1,
w_1) \Big \}\\
\qquad \qquad \qquad -it^{-1} (B_0 + B_S(W,W))W \ , \\
\\
R_2(W,S) = - t^{-2} \Big \{ s_0 \cdot \nabla s_1 + s_1 \cdot \nabla 
s_0 + s_1 \cdot \nabla s_1
\Big \} + t^{-1} \nabla B_L(w_1, w_1)\ . \end{array} \right . \eqno(7.8)$$
\noi Note that the term with $B_0 + B_S (W,W)$ in (\ref{7.1e}) is 
regarded as short range and
not included in the definition of $(W,S)$. \par

We now turn to the derivation of the estimates 
(\ref{6.45e})-(\ref{6.48e}). The regularity
properties of $(W,S)$ used in Section 6 follow from similar but 
simpler estimates.  We first
estimate all the terms not containing $B_0$. \\

\noi {\bf Lemma 7.1.} {\it Let $0 < \beta < 1$, $k_+ \geq 3$, $w_+ 
\in H^{k_+}$ and $a_+ =
|w_+|_{k_+}$. Then the following estimates hold~:
$$|w_0|_{k_+} \ \leq a_+ \eqno(7.9)$$

$$\parallel \omega^m \ s_0\parallel_2 \ \leq \left \{ 
\begin{array}{ll} C\ a_+^2 \
\ell n \ t &\qquad \hbox{for}\ 0 \leq m \leq k_+ \\
\\
C \ a_+^2 \ t^{\beta (m - k_+)} &\qquad \hbox{for}\ m > k_+ \ ,
\end{array} \right .   \eqno(7.10)$$

$$|w_1|_{k_+ - 1} \leq C \ a_+^3\ t^{-1} (1 + \ell n\ t) \ , \eqno(7.11)$$

$$\parallel \omega^m \ s_1\parallel_2 \ \leq \left \{ 
\begin{array}{l} C\ a_+^4 \ t^{-1} (1
+ \ell n \ t)^2 \quad \qquad \hbox{for}\ 0 \leq m \leq k_+ - 1 \\ \\
C \ a_+^4 \ t^{-1 + \beta (m +1- k_+)} (1 + \ell n \ t) \\
\qquad \hbox{for}\ k_+ - 1 < m < k_+ -1 +
\beta^{-1} \ , \end{array} \right . \eqno(7.12)$$
 
$$\parallel \omega^m \ R_2(W,S)\parallel_2 \ \leq \left \{ 
\begin{array}{l} C(a_+) \ t^{-3} (1
+ \ell n \ t)^3 \quad \qquad \hbox{for}\ 0 \leq m \leq k_+ - 2 \\ \\
C(a_+) \ t^{-3 + \beta (m +2- k_+)} (1 + \ell n \ t)^2 \\ \qquad 
\hbox{for}\ k_+ - 2 < m < k_+
-2 + \beta^{-1} \ , \end{array} \right . \eqno(7.13)$$

$$\parallel \omega^m (Q(S,w_1) + Q(s_1, w_0))\parallel_2 \ \leq C(a_+)
\ t^{-1} (1 + \ell n \ t)^2 \qquad \hbox{for}\ 0 \leq m \leq k_+ - 
2\eqno(7.14)$$

$$\parallel \omega^m(B_S(W,W)W)\parallel_2 \ \leq C(a_+) \ t^{-\beta 
(k_+ - m + 1)} \parallel
\qquad \hbox{for}\ 0 \leq m \leq k_+ - 1 \ .  \eqno(7.15)$$
  \vskip 3 truemm }

\noi {\bf Proof.} (7.9) is trivial. \par
(7.10). By (\ref{3.10e}), Lemma 3.2 and (\ref{3.9e}) we estimate
$$\hskip - 4 truecm \parallel \omega^m \ s_0 \parallel_2 \ \leq 
\displaystyle{\int_1^t} dt'\
t'^{-1} \parallel \omega^{m+1}\ B_L(w_0 (t'), w_0(t'))\parallel_2$$
$$\leq \left \{ \begin{array}{l} C \displaystyle{\int_1^t} dt' \ 
t'^{-1}\ I_{m} \left (
\parallel \omega^m w_0(t')\parallel_2 \ \parallel w_0(t')\parallel_{\infty}
\right ) \leq C a_+^2 \ell n \ t \\ \hskip 6 truecm \hbox{for} \ 0 
\leq m \leq k_+ \\
\\
C \displaystyle{\int_1^t} dt' \ t'^{-1+ \beta (m - k_+)}\ I_{k_+} \left (
\parallel \omega^{k_+} w_0(t')\parallel_2 \ \parallel w_0(t')\parallel_{\infty}
\right ) \leq C a_+^2 \ t^{\beta (m - k_+)} \\ \hskip 6 truecm 
\hbox{for} \ m > k_+ \ .
\end{array}\right .$$

(7.11). By Lemma 3.2 and (7.10), we estimate
$$\parallel Q(s_0, w_0)\parallel_2 \ \leq C \parallel \nabla s_0 
\parallel_2 \left ( \parallel
\omega^{3/2} w_0 \parallel_2 \ + \ \parallel w_0 \parallel_{\infty} 
\right ) \leq C\ a_+^3
\ell n\ t$$

$$\parallel  \omega^{k_+ -1} Q(s_0, w_0)\parallel_2 \ \leq C \Big \{ 
\parallel \omega^{k_+} s_0
\parallel_2 \left ( \parallel \omega^{3/2} w_0 \parallel_2 \ + \ \parallel w_0
\parallel_{\infty} \right ) $$ $$+ \left ( \parallel \omega^{3/2} s_0 
\parallel_2 \ + \ \parallel
s_0 \parallel_{\infty} \right ) \parallel \omega^{k_+} w_0
\parallel_2 \Big \} \leq C\ a_+^3 \ell n\ t$$
\noi from which (7.11) follows by integration. \par

(7.12). By Lemma 3.2 and (7.10) we estimate
$$\parallel \omega^m ( s_0 \cdot \nabla s_0 )\parallel_2 \ \leq C\parallel
\omega^{m+1}\ s_0\parallel_2 \left ( \parallel \omega^{3/2} s_0 
\parallel_2 \ + \ \parallel s_0
\parallel_{\infty} \right ) $$
$$\leq \left \{ \begin{array}{ll} C a_+^4 (\ell n \ t)^2 &\qquad 
\hbox{for} \ m \leq k_+ -1\\
\\
C a_+^4 \ t^{\beta (m + 1-  k_+)} \ell n\ t&\qquad \hbox{for} \ m > k_+ - 1 \ .
\end{array}\right .\eqno(7.16)$$
\noi On the other hand
$$\hskip - 10 truecm \parallel \omega^{m+1} B_L(w_0, w_1) \parallel_2$$
$$\leq \left \{ \begin{array}{l} CI_m \left ( \parallel \omega^{m} 
w_0 \parallel_2 \
\parallel w_1 \parallel_{\infty} \ + \ \parallel \omega^{m} w_1 
\parallel_2  \ \parallel w_0
\parallel_{\infty}\right ) \leq C a_+^4 t^{-1}(1 + \ell n \ t) \\
\hskip 6.2 truecm \hbox{for} \ m \leq k_+ -1\\
\\
Ct^{\beta (m+1-k_+)}I_{k_+-1} \left ( \parallel \omega^{k_+-1} w_0 
\parallel_2 \
\parallel w_1 \parallel_{\infty} \ + \ \parallel \omega^{k_+-1} w_1 
\parallel_2  \ \parallel w_0
\parallel_{\infty}\right ) \\
\leq C a_+^4 \ t^{\beta (m + 1-  k_+)- 1} (1 + \ell n\
t)\qquad \hbox{for} \ m > k_+ - 1 \ . \end{array}\right .\eqno(7.17)$$

\noi (7.12) now follows from (7.16) and (7.17) by integration 
provided $\beta (m + 1 - k_+) <
1$. \par

(7.13). By Lemma 3.2 again, and by (7.10) (7.12) we estimate
$$\hskip -4 truecm\parallel \omega^m \left ( s_0 \cdot \nabla s_1 + 
s_1 \cdot \nabla s_0 + s_1
\cdot \nabla s_1 \right ) \parallel_2\ \leq$$
$$C \Big \{   \parallel \omega^{m+1} s_0
\parallel_2\   \left (  \parallel \omega^{3/2} s_1 \parallel_{2} + 
\parallel s_1
\parallel_{\infty}\right ) + \ \parallel \omega^{m+1} s_1 \parallel_2 
\left ( \parallel
\omega^{3/2} s_0 \parallel_2 \ + \ \parallel s_0 \parallel_{\infty} \right .$$
$$\hskip -7 truecm\left . + \parallel \omega^{3/2} s_1 \parallel_2 \ 
+ \ \parallel
s_1\parallel_{\infty} \right ) \Big \}$$
$$\leq \left \{ \begin{array}{ll} C(a_+) \ t^{-1} (1 + \ell n \ t)^3 
&\hbox{for} \ m
\leq k_+ - 2\\
\\
C(a_+) \ t^{-1 + \beta (m+2 - k_+)} (1 + \ell n\ t)^2 &\hbox{for} \ 
k_+ - 2 < m < k_+
- 2 + \beta^{-1} \ .
\end{array} \right . \eqno(7.18)$$

\noi On the other hand
$$\hskip - 10 truecm \parallel \omega^{m+1} B_L(w_1, w_1) \parallel_2$$
$$\leq \left \{ \begin{array}{l} CI_m \left ( \parallel \omega^{m} 
w_1 \parallel_2 \
\parallel w_1 \parallel_{\infty} \right ) \leq C a_+^6 t^{-2}(1 + 
\ell n \ t)^2 \qquad
\hbox{for} \ m \leq k_+ -1\\ \\
Ct^{\beta (m+1-k_+)}I_{k_+-1} \left ( \parallel \omega^{k_+-1} w_1 
\parallel_2 \
\parallel w_1 \parallel_{\infty} \right )
\leq C a_+^6 \ t^{\beta (m + 1-  k_+)- 2} (1 + \ell n\
t)^2\\
\hskip 9.3 truecm  \hbox{for} \ m > k_+ - 1 \ . \end{array}\right 
.\eqno(7.19)$$

(7.13) now follows from (7.8) (7.18) (7.19). \par
(7.14). By Lemma 3.2 again, and by (7.10) (7.11) (7.12) we estimate
$$\hskip -7 truecm\parallel \omega^m \left ( Q(s_0, w_1) + Q(s_1, w_0) + Q(s_1, w_1) \right )
\parallel_2$$  $$\hskip - 1 truecm \leq C \Big \{  \parallel \omega^{m+1} s_0 \parallel_2\left (
\parallel \omega^{3/2} w_1\parallel_2 \ + \ \parallel 
w_1\parallel_{\infty} \right ) + \
\parallel \omega^{m+1} w_1 \parallel_2 \left ( \parallel \omega^{3/2} 
s_0 \parallel_2 \right
.$$  $$\left. + \ \parallel s_0 \parallel_{\infty} \ + \parallel 
\omega^{3/2} s_1 \parallel_2
\ + \ \parallel s_1\parallel_{\infty} \right ) + \ \parallel 
\omega^{m+1} s_1 \parallel_2 \left
( \parallel \omega^{3/2} w_0 \parallel_2 \ + \ \parallel w_0 
\parallel_{\infty}\right . $$
$$\hskip - 1 truecm \left .  +\ \parallel \omega^{3/2} w_1 
\parallel_2 \ + \ \parallel w_1
\parallel_{\infty} \right ) + \parallel \omega^{m+1} w_0 \parallel_2 
\left ( \parallel
\omega^{3/2} s_1 \parallel_2 \ + \ \parallel s_1 \parallel_{\infty} 
\right ) \Big \}$$
$$\hskip - 6 truecm\leq Ca_+^5 \ t^{-1} (1 + \ell n \ t)^2 (1 +a^2 \ 
t^{-1} (1 + \ell n\
t))\eqno(7.20)$$
\noi from which (7.14) follows. \par

(7.15). As previously, we estimate for $0 \leq m \leq k_+$
$$\parallel \omega^m B_S(W,W) \parallel_2\ \leq \ \parallel \omega^m 
B_S(w_0,w_0) \parallel_2\
+ 2 \parallel \omega^m B_S(w_0,w_1) \parallel_2\ + \parallel \omega^m 
B_S(w_1,w_1) \parallel_2$$
$$\leq Ct^{\beta (m - k_+)} \Big \{ t^{-\beta} I_{k_+} \left ( 
\parallel \omega^{k_+} w_0
\parallel_2 \ \parallel w_0 \parallel_{\infty}\right ) + I_{k_+-1} 
\left ( \parallel
\omega^{k_+-1} w_0 \parallel_2\ \parallel w_1 \parallel_{\infty} \right .$$
$$\hskip - 6 truecm + \left . \parallel
\omega^{k_+-1} w_1 \parallel_2\left ( \parallel w_0 
\parallel_{\infty} \ + \ \parallel w_1
\parallel_{\infty}\right ) \right ) \Big \}$$
$$\hskip - 2 truecm \leq C t^{\beta (m - k_+)} \Big \{ t^{-\beta}\ 
a_+^2 + t^{-1} \ a_+^4 (1 +
\ell n\ t) + t^{-2} \ a_+^6(1 + \ell n\ t)^2 \Big \} \eqno(7.21)$$
\noi by (7.11). Therefore
$$\parallel B_S(W,W)W \parallel_2\ \leq \ \parallel B_S(W,W) \parallel_2\
\parallel W\parallel_{\infty} \ \leq C(a_+) \ t^{-\beta (k_+ + 1)}$$
$$ \parallel \omega^{k_+-1} (B_S(W,W)W) \parallel_2 \ \leq C \Big \{ 
\parallel \omega^{k_+-1}
B_S(W,W \parallel_2 \ \parallel W \parallel_{\infty}$$
$$ + \ \parallel B_S(W,W) \parallel_{\infty} \ \parallel 
\omega^{k_+-1} W \parallel_2 \Big \} \
\leq C (a_+) \ t^{-2\beta} \eqno(7.22)$$
\noi which yields (7.15) by interpolation.\par \nobreak \hfill
$\sq$\par

We next estimate the terms in $R_1$ containing $B_0$.\\

\noi {\bf Lemma 7.2} {\it Let $0 < \beta < 1$. Let $1/2 < \lambda_0 < 
2$ and $k_+ \geq
2\lambda_0 \vee 3$. Let $B_0$ satisfy the estimates (\ref{3.17e}) for 
$0 \leq m \leq 2$, let
$w_+ \in H^{k_+}$ and assume that $B_0$ and $w_+$ satisfy the 
estimate (\ref{5.35e}) for all
multi-indices $\alpha_1$, $\alpha_2$ with $0 \leq |\alpha_1| \leq 2$ 
and $0 \leq |\alpha_2 | <
2 \lambda_0$. Then the following estimate holds for all $m$, $0 \leq 
m \leq 2$, and all $t \geq
1$.}
$$\parallel \omega^m ( B_0 \ W)\parallel_2 \ \leq C \ t^{-\lambda_0 + 
m} \ . \eqno(7.23)$$
\vskip 3 truemm

\noi {\bf Proof.} The contribution of $w_0$ to (7.23) is estimated by 
Lemma 5.1 with $\bar{m} =
2$. In order to estimate the contribution of $w_1$, we decompose $w_1 
= w'_1 + w''_1$ where
\begin{eqnarray*}
&&w'_1 = - \int_t^{\infty} dt' \ t'^{-2} \ Q(s_0(t'), w_+)\\
&&w''_1(t) = (1 - U^*(1/t)) \int_t^{\infty} dt'\ t'^{-2} U(1/t') Q(s_0 (t'),
w_0 (t'))\\
&&+ \int_t^{\infty} dt'\ t'^{-2} \Big \{ \left (1 - U(1/t')) 
Q(s_0(t'), w_0(t')) +
Q(s_0(t'), (1 - U(1/t'))w_+ \right )\Big \} \ .
\end{eqnarray*}
\noi We first consider
$$B_0 (t) w'_1(t) = - \int_t^{\infty} dt'\ t'^{-2} \Big \{ s_0 (t') 
\cdot B_0 (t) \nabla w_+ +
(1/2) (\nabla \cdot s_0)(t') B_0(t) w_+ \Big \}$$
\noi We estimate

$$\parallel B_0(t) w'_1(t) \parallel_2 \ \leq \int_t^{\infty} dt' \ 
t'^{-2} \Big \{ \parallel
s_0 (t') \parallel_{\infty} \ \parallel B_0(t)\nabla w_+ \parallel_2$$
$$+ \ \parallel (\nabla \cdot s_0 )(t') \parallel_{\infty} \ 
\parallel B_0 (t) w_+ \parallel_2
\Big \} \leq C \ t^{-\lambda_0 + 1/2} \int_t^{\infty} dt' \ t'^{-2} 
\ell n \ t'$$
$$ \leq C \ t^{-\lambda_0 - 1/2} (1 + \ell  n\ t) \eqno(7.24)$$
\noi by (5.35) and (7.10). \par

Similarly, we estimate
$$\parallel \Delta (B_0 (t) w'_1(t))\parallel_2 \ \leq C 
\int_t^{\infty} dt' \ t'^{-2} \Big \{
\parallel \Delta s_0 \parallel_3 \ \parallel B_0 \nabla w_+ \parallel_6$$
$$+ \ \parallel \nabla s_0 \parallel_{\infty} \ \parallel \nabla (B_0 
\nabla  w_+) \parallel_2
\ + \ \parallel s_0 \parallel_{\infty} \ \parallel \Delta (B_0 \nabla 
w_+) \parallel_2 \ + \ \parallel \Delta
\nabla \cdot s_0 \parallel_2 \ \parallel B_0 w_+ \parallel_{\infty}$$
$$+ \ \parallel \nabla^2 s_0 \parallel_6 \ \parallel \nabla (B_0 
w_+)\parallel_3 \ + \ \parallel
\nabla \cdot  s_0 \parallel_{\infty} \ \parallel \Delta (B_0 w_+) 
\parallel_2 \Big \}
\eqno(7.25)$$
\noi where $s_0 = s_0(t')$ and $B_0 = B_0(t)$, and therefore by (7.10)
$$\parallel \Delta (B_0\ w'_1)\parallel_2 \ \leq C\ a_+^2 \ t^{-1}(1 
+ \ell n\ t)
\Big\{ \left ( \parallel B_0 \parallel_{\infty} \ + \ \parallel 
\nabla B_0 \parallel_3
\right ) a_+$$
$$+ \ \parallel (\Delta B_0)  \nabla w_+ \parallel_2 + \parallel 
(\Delta B_0) w_+ \parallel_2
\Big \}$$
$$\leq C\ a_+^2 \left ( a_+ \ b_0 \ t^{-1/3} + b_1 \ t^{-\lambda_0 + 
3/2} \right ) (1 + \ell n\
t) \eqno(7.26)$$
\noi by (3.17) and (5.35). \par

We next estimate the contribution of $w''_1$. By the same estimates 
as for $w_1$ (see the
proof of (7.11)) we obtain
$$|w''_1|_{k_+-1} \leq C\ a_+^3 \ t^{-1} (1 + \ell n\ t)$$
$$|w''_1|_{k_+-3} \leq C\ a_+^3 \ t^{-2} (1 + \ell n\ t)$$
\noi where we have used the fact that the factors $(1 - 
U^{(*)}(1/t))$ can be replaced by
$t^{-1}\Delta$ for the purpose of the second estimate, and therefore
$$\parallel \omega^m w''_1 \parallel_2 \ \leq C\ a_+^3\ t^{-2+ m/2} 
(1 + \ell n\ t) \quad
\hbox{for} \ 0 \leq m \leq k_+ - 1  \eqno(7.27)$$
\noi by interpolation. By Lemma 3.2 and (3.17) we then obtain
$$\parallel \omega^m (B_0 \ w''_1) \parallel_2 \ \leq C \left ( 
\parallel \omega^m B_0
\parallel_{\infty} \ \parallel w''_1 \parallel_2 \ + \ \parallel 
B_0\parallel_{\infty} \
\parallel \omega^m w''_1 \parallel_2 \right )$$
$$ \leq C \ b_0 \ a_+^3\ t^{-2 + m} \ (1 + \ell n \ t) \ . \eqno(7.28)$$
\noi Collecting (7.24) (7.26) (7.28) and the estimates of $B_0 w_0$ 
coming from Lemma 5.1
yields (7.23) for $0 \leq m \leq 2$. \par \nobreak \hfill
$\sq$\par

We can now collect Proposition 6.3 and Lemmas 7.1 and 7.2 to obtain 
the main technical result
on the Cauchy problem for the auxiliary system in the difference form 
(2.30). We again
keep the assumptions on $B_0$ in the implicit form of the estimates 
(3.17) and (\ref{5.35e}),
which can however be replaced by sufficient conditions on $(w_+, A_+, 
\dot{A}_+)$ by the use
of Lemmas 3.5 and 5.2. \\

\noi {\bf Proposition 7.1.} {\it Let $1 < k \leq \ell$ and $\ell > 
3/2$. Let $\beta$,
$\lambda_0$ and $\lambda$ satisfy
$$0 < \beta < 2/3 \quad , \quad \lambda > 0\quad , \quad \lambda + k 
< \lambda_0 <
2\quad , \quad \lambda_0 > \beta (\ell + 1)\ . \eqno(7.29)$$
\noi Let $k_+$ satisfy
$$k_+ \geq  k + 2 \quad , \quad k_+ \geq 2 \lambda_0\quad , \quad 
\beta (k_+ + 1) \geq
\lambda_0 \quad , \quad \beta (\ell + 3 - k_+) < 1\ . \eqno(7.30)$$
\noi Let $w_+ \in H^{k_+}$, let $B_0$ satisfy the estimates 
(\ref{3.17e}) for $0 \leq m \leq
k$ and let $(w_+, B_0)$ satisfy the estimates (\ref{5.35e}) for all 
multi-indices $\alpha_1$,
$\alpha_2$ with $0 \leq |\alpha_1| \leq 2$ and $0 \leq |\alpha_1| < 
2\lambda_0$. Let
$(W,S)$ be defined by (7.3) (7.5) (7.7). Then \par

(1) $(W,S)$ satisfy the estimates (6.45) (6.46) (6.47) (6.48), with 
$0 < \eta < 1 - 3 \beta
/2$ in (\ref{6.46e}). \par

(2) All the statements of Proposition 6.3 hold. }\\

\noi {\bf Proof.} It follows from Lemmas 7.1 and 7.2 that all the 
assumptions of Proposition
6.3, and in particular the estimates (6.45)-(6.48), are 
satisfied.\par \nobreak \hfill
$\sq$\par

\mysection{Wave operators and asymptotics for (u, A)}
\hspace*{\parindent} In this section we complete the construction of 
the wave operators for
the system (1.1) (1.2) and we derive asymptotic properties of 
solutions in their range. The
construction relies in an essential way on Proposition 7.1. So far we have worked with the
system (\ref{2.20e}) for $(w, s)$ and the first task is t 
o reconstruct the phase $\varphi$.
Corresponding to $S = s_0 + s_1$, we define $\phi = \varphi_0 + 
\varphi_1$ where
$$\varphi_0 = - \int_1^t dt'\ t'^{-1} \ B_L^{\infty} (w_0(t'), w_0
(t')) \eqno(8.1)$$
$$\varphi_1 = - \int_t^{\infty} dt' (2t'^2)^{-1} |s_0(t')|^2 + 2 
\int_t^{\infty} dt' \ t'^{-1}
\ B_L^{\infty} (w_0(t') , w_1(t')) \eqno(8.2)
$$
\noi so that $s_0 = \nabla \varphi_0$ and $s_1 = \nabla \varphi_1$. \par

Let now $(q, \sigma )$ be the solution of the system (\ref{2.30e}) 
constructed in Proposition
6.3 part (2) and let $(w, s) = (W,S) + (q, \sigma )$. We define
$$\psi = - \int_t^{\infty} dt' (2t'^2)^{-1} \left ( \sigma \cdot 
(\sigma + 2S) + s_1 \cdot (
s_1 + 2s_0)\right ) (t')$$
$$+ \int_t^{\infty} dt' \ t'^{-1} \left ( B_L^{\infty} (q, q) + 2
B_L^{\infty}(W,q) + B_L^{\infty} (w_1, w_1) \right ) (t') \eqno(8.3)$$
\noi which is taylored to ensure that $\nabla \psi = \sigma$, given 
the fact that $s_0$, $s_1$
and $\sigma$ are gradients. The integral converges in $\dot{H}^1$, as 
follows from
(\ref{6.49e}) (\ref{6.50e}) and from the estimate (see the proof of 
(\ref{6.4e}))
$$\partial_t \parallel \sigma \parallel_2 \ \leq t^{-2} \parallel 
\nabla \sigma \parallel_2
\left ( \parallel s \parallel_{\infty} \ + \ \parallel \nabla S 
\parallel_3 \right ) + t^{-1}\
a\  I_0 (\parallel q \parallel_2)$$
$$\hskip - 4 truecm + t^{-1 + 3 \beta /2} \ I_{-3/2} \left ( 
\parallel q \parallel_2^2 \right )
+ \ \parallel R_2(W,S) \parallel_2$$
$$\leq C \left ( t^{-2 - \lambda_0 + \beta} (1 + \ell n \ t) + 
t^{-1-\lambda_0} +
t^{-1-2\lambda_0 + 3 \beta /2} + t^{-3} (1 + \ell n\ t)^3\right )$$
$$\hskip - 9 truecm \leq C\ t^{-1 - \lambda_0} \ . \eqno(8.4)$$
\noi Furthermore, this implies that
$$\parallel \nabla \psi \parallel_2 \ = \ \parallel \sigma 
\parallel_2 \ \leq C\
t^{-\lambda_0} \ . \eqno(8.5)$$
\noi Finally we define $\varphi = \phi + \psi$ so that $\nabla 
\varphi = s$, and $(w,
\varphi )$ solves the system (\ref{2.18e}). For more details on the 
reconstruction of
$\varphi$ from $s$, we refer to Section 7 of \cite{6r}. \par

We can now define the wave operators for the system (1.1) (1.2) as 
follows. We start from the
asymptotic state $(u_+, A_+, \dot{A}_+)$ for $(u, A)$. We define $w_+ 
= Fu_+$, we define
$(W,S)$ by (7.3) (7.5) (7.7) and $B_0$ by (\ref{2.3e}) (\ref{2.13e}), namely
$$A_0 = \dot{K}(t) \ A_+ + K(t) \ \dot{A}_+ = t^{-1} \ D_0\ B_0 \ .$$
\noi We next solve the system (\ref{2.30e}) with $t_1 = \infty$ and 
with initial time $t_0 =
\infty$ for $(q, \sigma )$ by Proposition 6.3, part (2), we define 
$(w, s) = (W,S) + (q ,
\sigma )$ and we reconstruct $\varphi$ from $s$ as explained above, 
namely $\varphi =
\varphi_0 + \varphi_1 + \psi$ with $\varphi_0$, $\varphi_1$ and 
$\psi$ defined by (8.1) (8.2)
(8.3). We finally substitute $(w, \varphi )$ thereby obtained into 
(\ref{2.11e}) (\ref{2.2e}),
thereby obtaining a solution $(u, A)$ of the system (1.1) (1.2). The 
wave operator is defined
as the map $\Omega : (u_+, A_+, \dot{A}_+) \to (u, A)$. \par

In order to state the regularity properties of $u$ that follow in a 
natural way from the
previous construction, we introduce appropriate function spaces. In 
addition to the operators
$M = M(t)$ and $D = D(t)$ defined by (2.8) (2.9), we introduce the operator
$$J = J(t) = x + it \ \nabla \ , \eqno(8.6)$$
\noi the generator of Galilei transformations. The operators $M$, 
$D$, $J$ satisfy the
commutation relation
$$i\ M\ D \ \nabla = J\ M\ D \ . \eqno(8.7)$$
\noi For any interval $I \subset [1, \infty )$ and any $k \geq 0$, we 
define the space
$${\cal X}^k(I) = \Big \{ u : D^*M^* u \in {\cal C} (I, H^k) \Big \}$$
$$\hskip 1.9 truecm = \Big \{ u : <J(t)>^k \ u \in {\cal C} (I, L^2) 
\Big \} \eqno(8.8)
$$
\noi where $< \lambda > = (1 + \lambda^2)^{1/2}$ for any real number 
or self-adjoint operator
$\lambda$ and where the second equality follows from (8.7). \par

We now collect the information obtained for the solutions of the 
system (1.1) (1.2) and state
the main result of this paper as follows.\\

\noi {\bf Proposition 8.1.} {\it Let $1 < k \leq \ell$, $\ell > 3/2$. 
Let $\beta$, $\lambda_0$
and $\lambda$ satisfy
$$0 < \beta < 2/3 \quad , \quad \lambda > 0 \quad , \quad \lambda + k 
< \lambda_0 < 2 \quad ,
\quad \lambda_0 > \beta (\ell + 1) \ . \eqno(7.29)$$
\noi Let $k_+$ satisfy
$$k_+ \geq  k + 2 \quad , \quad k_+ \geq 2 \lambda_0\quad , \quad 
\beta (k_+ + 1) \geq
\lambda_0 \quad , \quad \beta (\ell + 3 - k_+) < 1\ . \eqno(7.30)$$
\noi Let $u_+ \in F H^{k_+}$, let $w_+ = F u_+$ and $a_+ = 
|w_+|_{k_+}$. Let $(A_+, \dot{A}_+)
\in H^k \oplus H^{k-1}$. Let $A_0$ defined by (\ref{2.3e}) satisfy 
the estimates
$$\parallel \omega^m\ A_0(t) \parallel_r\ \leq b_0\ t^{2/r-1} 
\eqno(3.15)\equiv(8.9)$$
\noi for $0 \leq m \leq k$, $2 \leq r \leq \infty$ and all $t \geq 
1$, and the estimates
$$\parallel \left ( \partial^{\alpha_1}A_0 \right ) \left ( \left ( 
\partial^{\alpha_2} w_+
\right ) (x/t)\right ) \parallel_2 \ \leq b_1\ t^{-\lambda_0 + (1 + 
|\alpha_2|)/2}
\eqno(8.10)$$
\noi for all multi-indices $\alpha_1$, $\alpha_2$ with $0 \leq 
|\alpha_1| \leq 2$ and $0 \leq
|\alpha_2| < 2 \lambda_0$. Let $(W, S)$ be defined by (7.3) (7.5) 
(7.7) and let $\phi =
\varphi_0 + \varphi_1$ be defined by (8.1) (8.2). Then \par

(1) There exists $T$, $1 \leq T < \infty$ and there exists a unique 
solution $(u, A)$ of the
system (1.1) (2.2) with $u \in {\cal X}^k([T, \infty ))$, $(A, 
\partial_t A) \in {\cal C}
([T, \infty ), H^k \oplus H^{k-1})$ such that $u$ can be represented as
$$u = MD \exp (-i \varphi ) w \eqno(2.11)\sim(8.11)$$
\noi where $\varphi = \phi + \psi$ with $\psi$ defined by (8.3) and 
$(q, \sigma ) = (w, s) -
(W,S)$, and where $(w, s)$ is a solution of the system (\ref{2.20e}) 
with $t_1 = \infty$ such
that $(w, s) \in {\cal C}([T, \infty ) , X^{k,\ell})$ and such that
$$\parallel q(t) \parallel_2 \ \equiv \ \parallel w(t) - W(t) 
\parallel_2 \ \leq C
\ t^{-\lambda_0} \eqno(8.12)$$
$$\parallel \omega^k q(t) \parallel_2 \ \equiv \ \parallel 
\omega^k(w(t) - W(t))\parallel_2 \
\leq C\ t^{-\lambda} \eqno(8.13)$$
$$\parallel \omega^m \sigma (t) \parallel_2 \ \equiv \ \parallel 
\omega^m(s(t) -
S(t))\parallel_2 \ \leq C\ t^{-\lambda_0+\beta m} \quad \hbox{for}\ 0 
\leq m \leq
\ell + 1 \ .\eqno(8.14)$$
\indent That solution is obtained as $(u, A) = \Omega (u_+, A_+, 
\dot{A}_+)$ as defined above.
\par

(2) The solution $(u, A)$ satisfies the following estimates~:
$$\parallel u(t) - M(t) \ D(t) \exp (- i \phi (t)) W(t) \parallel_2 \ 
\leq C(a_+, b_0, b_1)
t^{-\lambda_0} \eqno(8.15)$$
$$\parallel |J(t)|^k \left ( \exp (i \phi (t, x/t)) u(t) - M(t) \ 
D(t) \ W(t) \right )
\parallel_2 \ \leq C(a_+, b_0, b_1)t^{-\lambda} \eqno(8.16)$$

$$\parallel u(t)- M(t) \ D(t) \exp (i \phi (t)), W(t)  \parallel_r \ 
\leq C(a_+, b_0,
b_1)t^{-\lambda_0 +(\lambda_0 - \lambda ) \delta (r)/k} \eqno(8.17)$$
\noi for $0 \leq \delta (r) = (3/2 - 3/r) \leq [k \wedge 3/2]$. \par

Define in addition
$$A_2 = A - A_0 - A_1^{\infty} (|DW|^2) \ . \eqno(8.18)$$
\noi Then the following estimates hold~:
$$\parallel A_2(t)\parallel_2 \ \leq C(a_+, b_0, b_1) \ t^{-\lambda_0 
+ 1/2} \ . \eqno(8.19)$$
\noi Furthermore, for $3/2 < k (<2)$~:
$$\parallel \nabla A_2(t)\parallel_2\ \leq C(a_+, b_0, b_1) \ 
t^{-\lambda_0-1/2} \eqno(8.20)$$
$$\parallel \omega^k \nabla A_2(t)\parallel_2 \ \leq \ C(a_+, b_0, 
b_1) \ t^{-\lambda -k-1/2}
\ , \eqno(8.21)$$
\noi while for $(1 <) k < 3/2$~:
$$\parallel \nabla A_2(t)\parallel_2 \ \leq C(a_+, b_0, b_1) \left ( 
t^{-\lambda_0 -
1/2} + t^{-2\lambda_0 - 1/2 + (\lambda_0 - \lambda ) 3/2k} \right 
)\eqno(8.22)$$
$$\parallel \omega^{2k-1/2} A_2(t)\parallel_2 \ \leq C(a_+, b_0, b_1) 
t^{-\lambda -2k+1} \left
(  t^{-\lambda} + t^{k-3/2} \right ) \ . \eqno(8.23)$$ \noi A similar 
result holds for $k =
3/2$ with a $t^{\varepsilon}$ loss in the decay. } \\

\noi {\bf Proof.} \underbar{Part (1)} follows from Propositions 6.3 
part (2) and from
Proposition 7.1, supplemented with the reconstruction of $\varphi$ 
described above in this
section. In particular the estimate (8.10) is nothing but the 
estimate (\ref{5.35e}) expressed
in terms of $A_0$ instead of $B_0$ while the estimates (8.12) (8.13) 
(8.14) are essentially
(\ref{6.49e}) (\ref{6.50e}) supplemented with (8.4) (8.5). \par

\noi \underbar{Part (2)}. We first prove the estimates on $u$. From 
(8.11) with $\varphi =
\phi + \psi$ and from (8.7), it follows that
$$\parallel |J|^m\left ( \exp (i\ D_0\ \phi )u - MDW \right 
)\parallel_2\ = \parallel
\omega^m\left ( w\ e^{-i\psi} - W \right ) \parallel_2 \eqno(8.24)$$
\noi For $m = 0$, we estimate
\begin{eqnarray*}
\parallel w \exp (- i \psi ) - W \parallel_2  &\leq&  \parallel w\left ( \exp
(-i \psi ) - 1 \right ) \parallel_2 \ + \ \parallel  w - W\parallel_2\\
&\leq& \ \parallel
w\parallel_3 \ \parallel \psi\parallel_6 \ + \ \parallel q 
\parallel_2 \ \leq C\
t^{-\lambda_0}
\end{eqnarray*}
  \noi by (8.5), a Sobolev inequality and (8.12). This proves (8.15). For
$m = k$, we estimate by Lemma 3.2
\begin{eqnarray*}
&&\parallel \omega^k\left ( \exp (-i \psi ) w - W \right ) \parallel_2 \ \leq C
\Big \{ \parallel \omega^k \left ( \exp (-i \psi ) - 1 \right ) 
\parallel_3\ \parallel w
\parallel_6 \\
&&+ \ \parallel \exp (-i \psi ) - 1 \parallel_{\infty} \ \parallel 
\omega^k w \parallel_2 \ +
\ \parallel \omega^k (w - W)\parallel_2 \\
&&\leq C \ \parallel \omega^{k-1/2} \sigma \parallel_2 \ \exp \left ( 
C \parallel \psi
\parallel_{\infty} \right ) \parallel \nabla w \parallel_2 \ + \left 
( \parallel \sigma
\parallel_2 \ \parallel \nabla \sigma \parallel_2 \right )^{1/2} \ 
\parallel \omega^k w
\parallel_2\\
&&+ \ \parallel \omega^k q\parallel_2 \ \leq C \left ( t^{-\lambda_0 
+ \beta (k - 1/2)} +
t^{-\lambda_0 + \beta /2} + t^{-\lambda} \right ) \leq C\ t^{-\lambda}
\end{eqnarray*}

\noi by Lemma 3.3, by the Sobolev inequality
$$\parallel \psi \parallel_{\infty} \ \leq C \left ( \parallel \sigma 
\parallel_2 \ \parallel
\nabla \sigma \parallel_2 \right )^{1/2}$$
\noi and by (8.13) (8.14). This proves (8.16). \par

The estimate (8.17) follows immediately from (8.15) (8.16) and from 
the inequality
$$\parallel f\parallel_r\ = t^{-\delta (r)} \ \parallel D^* M^* f 
\parallel_r \ \leq  C\
t^{-\delta (r)} \parallel \omega^{\delta (r)} D^*M^*f \parallel_2$$
$$= C\ t^{- \delta (r)} \ \parallel |J(t)|^{\delta (r)} f \parallel_2 \ .$$
\indent We next prove the estimates on $A$. It follows from the 
definitions (\ref{2.2e})
(\ref{2.3e}) (\ref{2.4e}) (8.18) and from (\ref{2.13e}) (\ref{2.14e}) that
$$A_2 = t^{-1} \ D_0 \ B_1^{\infty} (q, q + 2W) \ . \eqno(8.25)$$
\noi It is therefore sufficient to estimate $B_1^{\infty}(q, q + 
2W)$. We omit the superscript
$\infty$ for brevity. We first estimate by (\ref{3.10e})
$$\parallel B_1(q, q + 2W) \parallel_2\ \leq C\ I_{-1} \left ( 
\parallel \omega^{-1} (q(q +
2W))\parallel_2 \right )$$ $$\leq C\ I_{-1} \left ( \parallel q 
\parallel_2 \ \parallel q + 2W
\parallel_3 \right ) \leq  C\ t^{-\lambda_0} \eqno(8.26)$$
\noi by Sobolev inequalities and by (8.12), since $q + 2W$ is bounded 
in $H^k$ and a fortiori
in $L^3$. This proves (8.19). \par

For $k > 3/2$, we estimate similarly
$$\parallel \nabla B_1(q, q + 2W) \parallel_2\ \leq  I_{0} \left ( 
\parallel q \parallel_2\
\parallel q + 2W \parallel_{\infty} \right ) \leq C\ t^{-\lambda_0} 
\eqno(8.27)$$
\noi by (\ref{3.10e}) and (8.12), since $q + 2W$ is bounded in 
$L^{\infty}$ in that case.
Furthermore, by (\ref{3.10e}), Lemma 3.2 and Sobolev inequalities
$$\parallel \omega^{k+1} \ B_1(q, q + 2W) \parallel_2\ \leq C\ I_k 
\left ( \parallel
\omega^{k} q\parallel_2 \left ( \parallel q\parallel_{\infty} \ + \ \parallel
W\parallel_{\infty} \right ) \right .$$
$$+ \left . \parallel \nabla q\parallel_2 \ \parallel \omega^{k+ 1/2}
W\parallel_2\right ) \leq C\ t^{-\lambda}\eqno(8.28)$$
\noi by (8.12) (8.13). The last two inequalities imply (8.20) and 
(8.21) respectively. \par

For $k < 3/2$, we must estimate $B_1(q,q)$ and $B_1(q, W)$ separately 
because $q$ is no longer
controlled in $L^{\infty}$. We estimate as before
$$\parallel \nabla B_1(q, W) \parallel_2 \ \leq I_0 \left ( \parallel 
q \parallel_2 \
\parallel W \parallel_{\infty} \right ) \leq C\ t^{-\lambda_0}$$
\noi by (8.12), while
$$\parallel \nabla B_1(q, q) \parallel_2 \ \leq I_0 \left ( \parallel 
q \parallel_4^2 \
\right ) \leq C\  t^{-2\lambda_0 + (\lambda_0 - \lambda )3/2k}$$
\noi by (8.12) (8.13), which together imply (8.22). \par

We next estimate by (\ref{3.12e}) and (8.13)
$$\parallel \omega^{2k-1/2} \ B_1(q, q) \parallel_2 \ \leq C\ 
I_{2k-3/2} \left ( \parallel
\omega^k q \parallel_2^2 \right ) \leq C\ t^{-2\lambda} \eqno(8.29)$$
\noi while by (\ref{3.10e}) and Lemma 3.2
$$\parallel \omega^{2k-1/2}\ B_1(q, W)\parallel_2 \ \leq C\ 
I_{2k-3/2} \left ( \parallel
\omega^{2k- 3/2} q \parallel_2 \ \parallel W\parallel_{\infty} \right .$$
$$\left .+ \parallel q\parallel_r \ \parallel \omega^{2k-3/2} W 
\parallel_{3/\delta } \right
)$$
$$\leq C\ I_{2k - 3/2} \left ( \parallel \omega^{2k-3/2} q\parallel_2 
\left ( \parallel W
\parallel_{\infty} \ + \ \parallel \omega^{3/2} W \parallel_2 \right 
) \right ) $$
\noi by Sobolev inequalities, with $1/2 < \delta = \delta (r) = 2k - 
3/2 < 3/2$,
$$\cdots \leq C\ t^{-\lambda - (\lambda_0 - \lambda ) (3/2k - 1)} 
\leq C\ t^{-\lambda - 3/2 +
k} \eqno(8.30)$$
\noi by interpolation between (8.12) and (8.13). Now (8.23) follows 
from (8.29) and (8.30).\par \nobreak \hfill
$\sq$\par

We conclude this section with some remarks on variations which can be 
made or attempted in
the formulation of Proposition 8.1. \\

\noi {\bf Remark 8.1.} We have stated the assumptions on $(A_+, 
\dot{A}_+)$ in an implicit
way in the form of conditions on the solution $A_0$ of the free wave 
equation that they
generate. Sufficient conditions for (8.9) and (8.10) to hold directly 
expressed in terms of
$(A_+, \dot{A}_+)$ and possibly $w_+$ can be found in Lemma 3.5 and 
Lemma 5.2, but those
conditions are far from being optimal (especially those of Lemma 5.2). \\

\noi {\bf Remark 8.2.} We have stated the result for the partly 
resolved system (1.1) (2.2)
instead of the original system (1.1) (1.2). This includes the fact 
that the initial time for
(1.2) is $t_1 = \infty$ from the outset and makes it unnecessary to 
specify the uniqueness
class for $A$, since (1.2) is solved explicitly by (2.2). On the 
other hand, the available
regularity for $A$ is stronger than stated, as follows from the 
assumption (8.9) on $A_0$,
from simple estimates on $A_1^{\infty}(|DW|^2)$ and from the 
remainder estimates
(8.19)-(8.23). \\

\noi {\bf Remark 8.3.} The asymptotic behaviour of the scalar field 
$A$ differs in an
important way from that of a solution of the free wave equation. 
Omitting the remainder $A_2$
which tends to zero in a sense made precise by the estimates 
(8.19)-(8.23), one is left with
$$A \sim A_0 + A_1^{\infty} (|DW|^2)$$
\noi and the last term is easily estimated by
$$\parallel \omega^m A_1^{\infty} (|DW|^2)\parallel_2 \ \leq C(a_+) \ 
t^{-m+ 1/2}$$
\noi for $0 \leq m \leq k_+$. Replacing $W$ by $w_+$ as a first 
approximation in that term,
one obtains
$$A_1^{\infty} (|Dw_+|^2) = t^{-1}\ D_0 \ B_1^{\infty} (w_+, w_+)$$
\noi with $B_1^{\infty}(w_+, w_+)$ constant in time. This yields a 
contribution to $A$ which
spreads by dilation by $t$ and decays as $t^{-1}$ in $L^{\infty}$ 
norm. That contribution
can in no obvious sense be regarded as small as compared with $A_0$. \\

\noi {\bf Remark 8.4.} One might be tempted to look for simpler 
asymptotic forms for $u$
and for $A$ by replacing for instance $W$ by $w_+$ in (8.15) (8.16) 
(8.18) and/or by
omitting a few factors $U^{(*)}(1/t)$ in (7.5) (7.7). This however 
would introduce errors at
least $O(t^{-1})$ and spoil the $t^{-\lambda_0}$ decay in (8.12) 
(8.15) (8.19) (8.20) (8.22).\\

\noi{\large \bf Acknowledgements.} We are grateful to Professor Yves 
Meyer for enlightening
conversations.
\newpage
\mysection*{Appendix A}
\hspace*{\parindent} In this appendix, we prove Warnings 4.1 and 4.2. \\

\noi {\bf Proof of Warning 4.1.} One sees easily that (\ref{4.57e}) 
with $y(1) = y_0 > 0$ has
a unique maximal increasing solution $y \in {\cal C}^1([1, 
T^*),{I\hskip-1truemm R}^+)$ for
some $T^* > 1$. We shall argue by contradiction by showing that if 
$T^*$ is sufficiently
large, then $y(t)$ is infinite for some $t < T^*$. By integration, 
(\ref{4.57e}) with $y(1)
= y_0$ is converted into the integral equation
$$y(t) = y_0 \exp \Big \{ (m - \beta_1)^{-1} \ t^{-\beta_1} 
\int_{1/t}^1 d\nu \ y(\nu t)
\left ( \nu^{-1-m} - \nu^{-1-\beta_1} \right ) \Big \} \ . \eqno({\rm A.1})$$
\noi We shall prove by induction that (A.1) implies a sequence of 
lower bounds $y(t) \geq
a_n\ t^{\alpha_n}$ with $\alpha_n$ rapidly growing and $a_n$ not too 
small. We start with
$a_0 = y_0$, $\alpha_0 = 0$. Substituting that lower bound into (A.1) yields
$$y \geq y_0 \exp \{ y_0 \ h(t) \}$$
\noi where
$$h(t) = (m - \beta_1)^{-1} \ t^{-\beta_1} \int_{1/t}^1 d\nu \left ( 
\nu^{-1-m} - \nu^{-1 -
\beta_1} \right )$$
$$ = (m - \beta_1)^{-1} \left ( m^{-1} \ t^{m-\beta_1} - \beta_1^{-1} 
\right ) + m^{-1} \
\beta_1^{-1} \ t^{-\beta_1}$$
\noi so that $y \geq a_1 t^{\alpha_1}$ provided
$$\ell n \ a_1/y_0 \leq (m - \beta_1)^{-1} \left ( - \alpha_1 \ \ell 
n \ \tau + y_0
( m^{-1} \tau - \beta_1^{-1} ) \right ) \eqno({\rm A.2})$$
\noi where $\tau = t^{m- \beta_1} \geq 1$. The minimum of the RHS is 
attained for $\tau = m
\alpha_1/y_0$, which we take $> 1$, and we can then take
$$a_1 = y_0 \exp \Big \{ (m - \beta_1)^{-1} \ \beta_1^{-1} \ y_0 \Big 
\} (e\ y_0/m\
\alpha_1)^{\alpha_1/(m - \beta_1)} \ .\eqno({\rm A.3})$$
\noi Here $\alpha_1$ is an arbitrary fixed parameter, which we take large. In
particular we impose $\alpha_1 > (m^{-1} y_0 \vee 2 \beta_1)$. \par

At the following steps of the iteration, it will be sufficient to replace
(A.1) by the lower bound obtained by letting $m$ decrease to $\beta_1$,
namely
$$y(t) \geq y_0 \exp \Big \{ t^{-\beta_1} \int_{1/t}^1 d \nu \ y(\nu t) \
\nu^{-1-\beta_1} ÷ |\ell n \ \nu | \Big \}\eqno({\rm A.4})$$
\noi or equivalently
$$y(t) \geq y_0 \exp \Big \{ \int_{1}^t d t'\  \ y(t') \
t'^{-1-\beta_1} \ \ell n (t/t') \Big \}\ . \eqno({\rm A.5})$$
\noi We now describe the determination of $(a_{n+1}, \alpha_{n+1}) = 
(a' , \alpha ')$ from
$(a_n, \alpha_n) = (a, \alpha )$. Substituting the induction 
assumption into (A.5), we obtain
for $\alpha > \beta_1$ (a condition that will be ensured below)
\begin{eqnarray*}
&&y \geq y_0 \exp \Big \{ a \int_1^t dt' \ t'^{-1-\beta_1 + \alpha} \ \ell n
(t/t')\Big \}\\
&&= y_0 \exp \Big \{ a(\alpha - \beta_1)^{-2} (\tau - 1 - \ell n\ \tau )
\Big\}
\end{eqnarray*} \noi where $\tau = t^{\alpha - \beta_1} > 1$. This 
implies $y \geq a't^{\alpha
'}$ provided  $$\ell n\ a'/y_0 \leq a(\alpha - \beta_1)^{-2} (\tau - 
1 - \ell n \
\tau ) - \alpha ' (\alpha - \beta_1)^{-1} \ell n\ \tau$$
$$\hskip - 2 truecm = \theta \left ( \widetilde{a} (\tau - 1) - 
(\widetilde{a} + 1) \ell n \
\tau \right )
  \eqno({\rm A.6})$$
\noi where
$$\theta = \alpha '/(\alpha - \beta_1) \quad , \quad \widetilde{a}  = 
a/ \alpha ' (\alpha -
\beta_1) \ .  \eqno({\rm A.7})$$
\noi The minimum over $\tau$ of the last member of (A.6) is attained 
for $\widetilde{a}\tau =
\widetilde{a} + 1$, and it suffices to impose
$$\ell n (a'/y_0) \leq \theta (1 - (\widetilde{a} + 1) \ell n (\widetilde{a} +
1)/\widetilde{a}) \equiv \theta (\ell n \ \widetilde{a} - 
f(\widetilde{a}))$$ \noi which allows
us to take $a' = y_0 \widetilde{a}^{\theta}$ provided
$$f(\widetilde{a}) \equiv (\widetilde{a} + 1) \ell n (\widetilde{a} + 
1) - \widetilde{a} \ \ell
n \ \widetilde{a} - 1 \leq 0$$
\noi a condition which is easily seen to hold for $\widetilde{a} \leq 
1/2$. \par

Finally we can take
$$\alpha ' = \theta (\alpha - \beta_1 ) \qquad , \quad a' = y_0 \left 
( a/\theta (\alpha -
\beta_1)^2 \right )^{\theta}  \eqno({\rm A.8})$$
\noi provided
$$a \leq \theta (\alpha - \beta_1)^2/2 \ . \eqno({\rm A.9})$$
\noi So far $\theta$ is a free parameter. For definiteness we choose 
$\theta = 2$, so that
after coming back to the original notation, (A.8) (A.9) become
$$\alpha_{n+1} = 2 (\alpha_n - \beta_1) \ , \eqno({\rm A.10})$$
$$ a_{n+1} = y_0 \ a_n^2 /4 (\alpha_n - \beta_1)^4 \ ,  \eqno({\rm A.11})$$
$$a_n \leq (\alpha_n - \beta_1)^2 \ .  \eqno({\rm A.12})$$
\noi (A.10) is readily solved by
$$\alpha_n = 2 \beta_1 + 2^{n-1} (\alpha_1 - 2\beta_1) \ .$$
\noi (A.12) is harmless and holds for all $n$ if it holds for $n = 1$ 
and if $y_0 \leq
4(\alpha_1 - \beta_1)^2$, which can be arranged by taking $\alpha_1$ 
sufficiently large.
(A.11) can be rewritten as
$${a_{n+1}\ y_0 \over 64 (\alpha_{n+1} - \beta_1)^4} = {y_0^2\ a_n^2 
\over 64^2 (\alpha_n -
\beta_1)^8} \left ( {2(\alpha_n - \beta_1) \over \alpha_{n+1} - 
\beta_1} \right )^4 \geq
\left \{ {a_n \ y_0 \over 64 (\alpha_n - \beta_1)^4} \right \}^2 
\eqno({\rm A.13})$$
\noi by
(A.10). Let now $t \geq 1$ and define
$$u_n = a_n \ t^{\alpha_n - 2 \beta_1} \ y_0 /64 (\alpha_n - \beta_1)^4 \ .$$
\noi It follows from (A.10) (A.13) that $u_{n+1} \geq u_n^2$ and
in particular $u_n \geq 1$ for all $n$ if $u_1 \geq 1$, namely if $t$ 
is sufficiently large in
the sense that
$$t^{\alpha_1 -2\beta_1} \geq (a_1 y_0)^{-1} \ 64(\alpha_1 - 
\beta_1)^4 \ .\eqno({\rm A.14})$$
\noi For such $t$, the condition $u_n \geq 1$ can be rewritten as
$$
y(t) \geq a_n\ t^{\alpha_n}\geq t^{2\beta_1}\ y_0^{-1} \ 64(\alpha_n 
- \beta_1)^4$$
$$\geq 4 t^{2\beta_1} \ y_0^{-1} \ 2^{4n} (\alpha_1 - 2 \beta_1)^4 \ 
. \eqno({\rm
A.15})$$
\noi Since the last member of (A.15) tends to infinity with $n$, such 
a $t$ cannot be smaller
than $T^*$, which proves finite time blow up. \par \nobreak \hfill
$\sq$\par
\noi {\bf Remark A1.} Since the RHS of (4.56) and (4.57) is 
decreasing in $\beta_1$ and
increasing in $m$, blow up in finite time for $(\beta_1, m)$ implies 
blow up in finite time for
$(\beta '_1, m)$ with $\beta '_1 \leq \beta_1$ and for $(\beta_1, 
m')$ with $m' \geq m$, while
the opposite situation prevails as regards the existence of global 
solutions. Actually it is
easy to see that (4.56) or (4.57) admits global solutions for small 
data if $\beta_1 > 0$ and
$m \leq \beta_1$. When coming back to the original equation (4.55), 
the condition of small
data becomes a condition of large $t_0$. \\

\noi {\bf Proof of Warning 4.2.} We want to prove finite time blow up 
for (4.61) with $y(t_0) =
y_0 > 0$. Omitting the second term in the RHS and integrating the 
remaining inequality, we
obtain
$$y \geq \left ( y_0^{1/k} + t - t_0 \right )^k \geq (t - t_0)^k \ . 
\eqno({\rm A.16})$$
\noi We next keep (A.16), omit the first term in the RHS of (4.61) 
and change $t$ to $t +
t_0$. It is then sufficient to prove blow up for
$$\left \{ \begin{array}{l} y \geq t^k \\
\\
  \partial_t y \geq
(t+t_0)^{-1-\beta_1}\ y^3 \ .\end{array}\right . \eqno({\rm A.17})$$
\noi For that purpose, we show by induction that $y$ satisfies
$$y(t) \geq y_n(t) \geq a_n \ t^{\alpha_n} 
(t+t_0)^{-(1+\beta_1)\gamma_n} \ ,  \eqno({\rm
A.18})$$
\noi starting with $a_0 = 1$, $\alpha_0 = k$ and $\gamma_0 = 0$ given 
by (A.17). We obtain
\begin{eqnarray*}
&&y_{n+1} = \int_0^t dt'(t_0 + t')^{-1-\beta_1} \ y_n^3(t')\\
&&\geq a^3_n \int_0^t dt' (t_0 + t')^{-(1 + \beta_1)(3\gamma_n + 1)} 
\ t'^{3\alpha_n} \\
&&\geq a_n^3 (t_0 + t)^{-(1 + \beta_1)(3\gamma_n+1)} \ t^{3 \alpha_n 
+ 1}(3 \alpha_n + 1)^{-1}
\end{eqnarray*}
\noi thereby ensuring (A.18) at the level $n+1$ if we choose
$$\alpha_{n+1} = 3\alpha_n + 1 \qquad , \quad \gamma_{n+1} = 
3\gamma_n + 1 \ ,   \eqno({\rm
A.19})$$
$$a_{n+1} = a_n^3/(3 \alpha_n + 1) \ .  \eqno({\rm A.20})$$
\noi (A.19) is readily solved by
$$\alpha_{n} = 3^n (k + 1/2) - 1/2 \qquad ,\quad \gamma_{n} = (3^n - 1)/2
\eqno({\rm A.21})$$
\noi so that
$$a_{n+1} \geq  a_n^3(k + 1/2)^{-1} \ 3^{-(n+1)}   \eqno({\rm A.22})$$
\noi or equivalently
$$b_{n+1} \geq b_n^3   \eqno({\rm A.23})$$
\noi where
$$b_{n} = a_n 3^{-n/2-3/4} \ (k + 1/2)^{-1/2} \ .  \eqno({\rm A.24})$$
\noi Let now $t > 0$ and
$$u_{n} = b_n \ t^{\alpha_n + 1/2} \ (t_0 + t)^{-(1 + 
\beta_1)(\gamma_n + 1/2)}\ .
\eqno({\rm A.25})$$
\noi It follows from (A.19) (A.23) that $u_{n+1} \geq u_n^3$ and in 
particular that $u_n \geq
1$ for all $n$ if $u_0 \geq 1$. The condition $u_0 \geq 1$ reduces to
$$t^{2k+1} (t_0 + t)^{-(1 + \beta_1)} \geq 3^{3/2} (k + 1/2) 
\eqno({\rm A.26})$$
\noi and holds for $t$ sufficiently large if $2k > \beta_1$. For such 
a $t$, by (A.18)
$$y \geq a_n\ t^{\alpha_n} (t_0 + t)^{-(1 + \beta_1)\gamma_n} \geq (k 
+ 1/2)^{1/2} \ 3^{n/2 +
3/4}\ t^{-1/2} (t_0 + t)^{(1 + \beta_1)/2} \ .  \eqno({\rm A.27})$$
\noi Since the last
member of (A.27) tends to infinity with $n$, such a $t$ cannot be 
smaller than the maximal time
$T^*$ of existence of the solution $y$ of (\ref{4.61e}), which proves 
finite time blow up.\par
\nobreak \hfill $\sq$\par

\end{document}